\documentclass[a4paper,12pt]{article}

\usepackage[top=2.5cm, bottom=2.5cm, left=2.5cm, right=2.5cm]{geometry}

\usepackage{tikz}
\usetikzlibrary{arrows.meta, decorations.markings}
\definecolor{darkgreen}{rgb}{0,0.5,0}

\usepackage{array}

\usepackage{amsfonts,amstext,amscd,bezier,amsthm,amssymb}
\usepackage[centertags]{amsmath}

\usepackage{mathtools}
\usepackage[integrals]{wasysym}
\usepackage{stmaryrd}

\usepackage[mathcal]{euscript}

\usepackage[utf8]{inputenc}
\usepackage[T1]{fontenc}

\usepackage[english]{babel}

\usepackage{cite}

\usepackage{mathptmx}

\usepackage[inline]{enumitem}

\usepackage[linktocpage, colorlinks, linkcolor=red, citecolor=blue, urlcolor=violet, bookmarks, bookmarksnumbered, pdfstartview={XYZ null null 1.00}]{hyperref}

\theoremstyle{plain}
\newtheorem{theo}{Theorem}[section]
\newtheorem{lemm}[theo]{Lemma}

\newtheorem{prop}[theo]{Proposition}
\theoremstyle{definition}
\newtheorem{defi}[theo]{Definition}
\newtheorem{remk}[theo]{Remark}
\newtheorem{expl}[theo]{Example}
\newenvironment{prf}[1][Proof]%
  {\begin{list}{}{%
    \setlength{\topsep}{0pt}%
    \setlength{\leftmargin}{0pt}%
    \setlength{\rightmargin}{0pt}%
    \setlength{\listparindent}{\parindent}%
    \setlength{\itemindent}{0pt}%
    \setlength{\parsep}{\parskip}}%
    \item[]{\bf #1. }}%
  {\hspace*{\fill} $\blacksquare$ \end{list} \medskip}
\newenvironment{prfsection}%
  {\begin{list}{}{%
    \setlength{\topsep}{0pt}%
    \setlength{\leftmargin}{0pt}%
    \setlength{\rightmargin}{0pt}%
    \setlength{\listparindent}{\parindent}%
    \setlength{\itemindent}{0pt}%
    \setlength{\parsep}{\parskip}}%
    \item[]}%
  {\hspace*{\fill} $\blacksquare$ \end{list} \medskip}
	
\DeclareMathOperator{\Ker}{Ker}
\DeclareMathOperator{\rank}{rk}
\DeclareMathOperator{\range}{Ran}
\DeclareMathOperator{\id}{Id}

\newcommand{\suchthat}{\;|\:}
\newcommand{\midsuchthat}{\;\middle|\:}
\newcommand{\transp}{{\mathrm{T}}}

\newcommand{\norm}[1]{\left\lVert #1\right\lVert}
\newcommand{\abs}[1]{\left\lvert #1\right \rvert}
\newcommand{\floor}[1]{\left\lfloor #1 \right\rfloor}
\newcommand{\ceil}[1]{\left\lceil #1 \right\rceil}

\newcommand{\scalprod}[2]{\left\langle #1, #2 \right\rangle}

\numberwithin{figure}{section}
\numberwithin{equation}{section}

\newcommand{\Epsilon}{\mathrm E}

\begin{document}

\setlist[enumerate, 1]{label={\textnormal{(\alph*)}}, ref={(\alph*)}, leftmargin=*}
\setlist[enumerate, 2]{label={\textnormal{(\roman*)}}, ref={(\roman*)}, leftmargin=*}

\title{Approximate and exact controllability of linear difference equations\thanks{This work is supported by a public grant overseen by the French National Research Agency (ANR) as part of the ``Investissement d'Avenir'' program, through the iCODE project funded by the IDEX Paris-Saclay, ANR-11-IDEX-0003-02. The second author was supported by a public grant as part of the Investissement d'avenir project, reference ANR-10-CAMP-0151-02-FMJH.}}
\author{
Yacine Chitour\thanks{Laboratoire des Signaux et Syst\`emes, Sup\'elec, and Universit\'e Paris Sud, Orsay, France.} ,
Guilherme Mazanti\thanks{Laboratoire de Mathématiques d'Orsay, Univ. Paris-Sud, CNRS, Université Paris-Saclay, 91405 Orsay, France} ,
Mario Sigalotti\thanks{Inria team CAGE, and CMAP, \'Ecole Polytechnique, Palaiseau, France.}
}
\maketitle

\begin{abstract}
In this paper, we study approximate and exact controllability of the linear difference equation $x(t) = \sum_{j=1}^N A_j x(t - \Lambda_j) + B u(t)$ in $L^2$, with $x(t) \in \mathbb C^d$ and $u(t) \in \mathbb C^m$, using as a basic tool a representation formula for its solution in terms of the initial condition, the control $u$, and some suitable matrix coefficients. When $\Lambda_1, \dotsc, \Lambda_N$ are commensurable, approximate and exact controllability are equivalent and can be characterized by a Kalman criterion. This paper focuses on providing characterizations of approximate and exact controllability without the commensurability assumption. In the case of two-dimensional systems with two delays, we obtain an explicit characterization of approximate and exact controllability in terms of the parameters of the problem. In the general setting, we prove that approximate controllability from zero to constant states is equivalent to approximate controllability in $L^2$. The corresponding result for exact controllability is true at least for two-dimensional systems with two delays.
\end{abstract}

\tableofcontents

\paragraph*{Notations}

In this paper, we denote by $\mathbb N$ and $\mathbb N^\ast$ the sets of nonnegative and positive integers, respectively. For $a, b \in \mathbb R$, we write the set of all integers between $a$ and $b$ as $\llbracket a, b \rrbracket = [a, b] \cap \mathbb Z$, with the convention that $[a, b] = \emptyset$ if $a > b$. For $\Lambda \in \mathbb R^N$, we use $\Lambda_{\min}$ and $\Lambda_{\max}$ to denote the smallest and the largest components of $\Lambda$, respectively. For $\xi \in \mathbb R$, the symbol $\floor{\xi}$ is used to the denote the integer part of $\xi$, i.e., the unique integer such that $\xi - 1 < \floor{\xi} \leq \xi$, $\ceil{\xi}$ denotes the unique integer such that $\xi \leq \ceil{\xi} < \xi + 1$, and we set $\{\xi\} = \xi - \floor{\xi}$. For $z \in \mathbb C$, the complex conjugate of $z$ is denoted by $\overline z$. We write $\overline X$ for the closure of the subset $X$ of a topological space. By convention, we set the sum over an empty set to be equal to zero, $\inf\emptyset = +\infty$, and $\sup\emptyset = -\infty$. The characteristic function of a set $A \subset \mathbb R$ is denoted by $\chi_A$.

The set of $d \times m$ matrices with coefficients in $K \subset \mathbb C$ is denoted by $\mathcal M_{d, m}(K)$, or simply by $\mathcal M_d(K)$ when $m = d$. The identity matrix in $\mathcal M_d(\mathbb C)$ is denoted by $\id_d$, the zero matrix in $\mathcal M_{d, m}(\mathbb C)$ is denoted by $0_{d, m}$, or simply by $0$ when its dimensions are clear from the context, and the transpose of a matrix $A \in \mathcal M_{d, m}(K)$ is denoted by $A^\transp$. We write $\mathrm{GL}_d(\mathbb C)$ for the general linear group of order $d$ over $\mathbb C$. The vectors $e_1, \dotsc, e_d$ denote the canonical basis of $\mathbb C^d$. For $p \in [1, +\infty]$, $\abs{\cdot}_{p}$ indicates both the $\ell^p$-norm in $\mathbb C^d$ and the corresponding induced matrix norm in $\mathcal M_{d, m}(\mathbb C)$. We denote the usual scalar product of two vectors $x, y \in \mathbb R^d$ by $x \cdot y$. The range of a matrix $M \in \mathcal M_{d, m}(\mathbb C)$ is denoted by $\range M$, and $\rank M$ denotes the dimension of $\range M$.

For $(A, B) \in \mathcal M_{d}(\mathbb C) \times \mathcal M_{d, m}(\mathbb C)$, the \emph{controllability matrix} of $(A, B)$ is denoted by $\mathcal C(A, B)$, and we recall that
\[
\mathcal C(A, B) = \begin{pmatrix}
B & A B & A^2 B & \cdots & A^{d-1} B \\
\end{pmatrix} \in \mathcal M_{d, d m}(\mathbb C).
\]
We also recall that a pair $(A, B) \in \mathcal M_{d}(\mathbb C) \times \mathcal M_{d, m}(\mathbb C)$ is said to be \emph{controllable} if $\rank \mathcal C(A, B) = d$.

The inner product of a Hilbert space $\mathsf H$ is denoted by $\scalprod{\cdot}{\cdot}_{\mathsf H}$ and is assumed to be anti-linear in the first variable and linear in the second one. The corresponding norm is denoted by $\norm{\cdot}_{\mathsf H}$, and the index $\mathsf H$ is dropped from these notations when the Hilbert space under consideration is clear from the context. For two Hilbert spaces $\mathsf H_1, \mathsf H_2$, the Banach space of all bounded operators from $\mathsf H_1$ to $\mathsf H_2$ is denoted by $\mathcal L(\mathsf H_1, \mathsf H_2)$, with its usual induced norm $\norm{\cdot}_{\mathcal L(\mathsf H_1, \mathsf H_2)}$. The adjoint of an operator $E \in \mathcal L(\mathsf H_1, \mathsf H_2)$ is denoted by $E^\ast$. When $\mathsf H_1 = \mathsf H_2 = \mathsf H$, we write simply $\mathcal L(\mathsf H)$ for $\mathcal L(\mathsf H, \mathsf H)$. The range of an operator $E \in \mathcal L(\mathsf H_1, \mathsf H_2)$ is denoted by $\range E$.

\section{Introduction}
\label{SecIntro}

This paper studies the controllability of the difference equation
\begin{equation}
\label{MainSyst}
\qquad x(t) = \sum_{j=1}^N A_j x(t - \Lambda_j) + B u(t),
\end{equation}
where $x(t) \in \mathbb C^d$ is the state, $u(t) \in \mathbb C^m$ is the control input, $N, d, m \in \mathbb N^\ast$, $\Lambda = (\Lambda_1, \dotsc, \Lambda_N) \in (0, +\infty)^N$ is the vector of positive delays, $A = (A_1, \dotsc, A_N) \in \mathcal M_d(\mathbb C)^N$, and $B \in \mathcal M_{d, m}(\mathbb C)$.

The study of the autonomous difference equation
\begin{equation}
\label{SystAutonomous}
x(t) = \sum_{j=1}^N A_j x(t - \Lambda_j)
\end{equation}
has a long history and its analysis through spectral methods has led to important stability criteria, such as those in \cite{Avellar1980Zeros} and \cite[Chapter 9]{Hale1993Introduction} (see also \cite{Melvin1974Stability, Cruz1970Stability, Henry1974Linear, Datko1977Linear, Hale1985Stability, Michiels2009Strong} and references therein). A major motivation for analyzing the stability of \eqref{SystAutonomous} is that it is deeply related to properties of more general neutral functional differential equations of the form
\begin{equation}
\label{NFDE}
\frac{d}{dt}\left(x(t) - \sum_{j=1}^N A_j x(t - \Lambda_j)\right) = f(x_t)
\end{equation}
where $x_t: [-r, 0] \to \mathbb C^d$ is given by $x_t(s) = x(t + s)$, $r \geq \Lambda_{\max}$, and $f$ is some function defined on a certain space (typically $\mathcal C^k([-r, 0], \mathbb C^d)$ or $W^{k, p}((-r, 0), \mathbb C^d)$); see, e.g., \cite{Cruz1970Stability, Datko1977Linear, Hale1985Stability, Ngoc2015Exponential}, \cite[Section 9.7]{Hale1993Introduction}. Another important motivation is that, using d'Alembert decomposition, some hyperbolic PDEs can be transformed by the method of characteristics into differential or difference equations with delays \cite{Coron2015Dissipative, Coron2008Dissipative, Cooke1968Differential, Fridman2010Bounds, Slemrod1971Nonexistence, Kloss2012Flow}, possibly with time-varying matrices $A_j$ \cite{Chitour2017Persistently, Chitour2016Stability}.

Several works in the literature have studied the control and the stabilization of neutral functional differential equations under the form \eqref{NFDE}. In particular, stabilization by linear feedback laws was addressed in \cite{Pandolfi1976Stabilization, Hale2002Strong, OConnor1983Stabilization}, with a Hautus-type condition for the stabilizability of \eqref{MainSyst} provided in \cite{Hale2002Strong}.

Due to the infinite-dimensional nature of the dynamics of difference equations and neutral functional differential equations, several different notions of controllability can be used, such as approximate, exact, spectral, or relative controllability \cite{Salamon1984Control, Chyung1970Controllability, Diblik2008Controllability, Pospisil2015Relative, Mazanti2017Relative}. Relative controllability was originally introduced in the study of control systems with delays in the control input \cite{Chyung1970Controllability} and consists in controlling the value of $x(T) \in \mathbb C^d$ at some prescribed time $T$. In the context of difference equations under the form \eqref{MainSyst}, it was characterized in some particular situations with integer delays in \cite{Diblik2008Controllability, Pospisil2015Relative}, with a complete characterization on the general case provided in \cite{Mazanti2017Relative}.

We consider in this paper the approximate and exact controllability of \eqref{MainSyst} in the function space $L^2((-\Lambda_{\max}, 0), \mathbb C^d)$. Such a problem is largely absent from the literature, with the notable exception of \cite{Salamon1984Control, OConnor1983Function}, where some controllability notions for neutral functional differential equations under the form \eqref{NFDE} are characterized in terms of corresponding observability properties, such as unique continuation principles, using duality arguments reminiscent of the Hilbert Uniqueness Method introduced later in \cite{Lions1986Controlabilite, Lions1988Exact}.

The above controllability problems have easy answers in some simple situations. Indeed, in the case of a single delay, approximate and exact controllability are equivalent to the standard Kalman controllability criterion for the pair $(A_1, B)$, i.e., the controllability of the finite-dimensional discrete-time system $x_{n+1} = A_1 x_n + B u_n$. More generally, when all delays are commensurable, i.e., integer multiples of a common positive real number, we reduce the problem to the single-delay case by the classical augmented state space technique (see, for instance, \cite[Chapter 4]{Franklin1997Digital}). The Kalman criterion can be interpreted as an explicit test for controllability since it yields a complex-valued function $F$ of the parameters of the problem, polynomial with respect to the coefficients of the matrices, such that controllability of a system is equivalent to $F$ not taking the value zero for that system.

%The above controllability problems have easy answers in some simple situations. Indeed, in the case of a single delay, both controllability notions are equivalent to the controllability of the finite-dimensional discrete-time system $x_{n+1} = A_1 x_n + B u_n$, i.e., to the standard Kalman controllability criterion for the pair $(A_1, B)$. More generally, when all delays are commensurable, i.e., integer multiples of a common positive real number, one reduces the problem to the single-delay case by a classical augmentation of the state space (for details, see Section \ref{SecCommensurable}). The Kalman criterion can be interpreted as an explicit test for controllability since it yields a complex-valued function $F$ of the parameters of the problem, polynomial with respect to the coefficients of the matrices, such that controllability of a system is equivalent to $F$ not taking the value zero for that system.

We are not aware of any result of this type in the incommensurable case, even though the problem seems natural and of primary importance if one is interested in linear controlled difference equations. We show in this paper that such an explicit test can be obtained at least in the first non-trivial incommensurable case, namely two-dimensional systems with two delays and a scalar input (Theorem \ref{MainTheo22}). Note that approximate and exact controllability are no more equivalent but we still characterize explicitly both of them.

Let us now describe the line of arguments we use to derive our results. The approximate controllability in the case of incommensurable delays is reduced to the existence of nonzero functions invariant with respect to a suitable irrational translation modulo $1$. The ergodicity of the latter yields a necessary condition for approximate controllability, which is also shown to be sufficient. As regards exact controllability, the strategy consists in approximating the original system by a sequence of systems $(\Sigma_n)_{n \in \mathbb N}$ with commensurable delays, and, for every $n \in \mathbb N$, the controllability of $\Sigma_n$ is equivalent to the invertibility of a Toeplitz matrix $M_n$, whose size tends to infinity. The heart of the argument boils down to bounding the norm of $M_n^{-1}$ uniformly with respect to $n$.

For more delays or in higher dimension, the existence of explicit controllability tests remains open. Characterizing approximate controllability using our techniques would amount to single out a tractable discrete dynamical system, generalizing the above-mentioned translation modulo $1$. Concerning exact controllability, the difficulty is that the above matrices $M_n$ are now block-Toeplitz. We believe that the general case is not an easy problem and additional techniques may be needed, for instance arguments based on the Laplace transform.

We also prove an additional result stating that approximate controllability from zero to constant states implies approximate controllability in $L^2$, and the same holds true for exact controllability at least for two-dimensional systems with two delays and a scalar input. The interest of this result lies in the fact that reachability of a finite-dimensional space is sufficient to deduce the reachability of the full $L^2$ space.

Throughout the paper, we rely on a basic tool for the controllability analysis of \eqref{MainSyst}, namely a suitable representation formula, describing a solution at time $t$ in terms of its initial condition, the control input, and some matrix-valued coefficients computed recursively (see Proposition~\ref{PropExplicit}). Such a formula, already proved in \cite{Mazanti2017Relative}, generalizes the ones obtained in \cite[Theorems 3.3 and 3.6]{Chitour2017Persistently} for the stability analysis of a system of transport equations on a network under intermittent damping, and the one obtained in \cite[Proposition~3.14]{Chitour2016Stability}, used for providing stability criteria for a non-autonomous version of \eqref{SystAutonomous}.

The plan of the paper goes as follows. In Section \ref{SecPrevResults} we discuss the well-posedness of \eqref{MainSyst}, present the explicit representation formula for its solutions, provide the definitions of $L^2$ approximate and exact controllability, and recall some of their elementary properties. Section \ref{SecCommensurable} considers the case of systems with commensurable delays, for which the usual technique of state augmentation is available. We prove that such a technique and our approach based on the representation formula from Section \ref{SecControlDelayWellPosed} both yield the same Kalman-like controllability criterion. The main results are provided in Sections \ref{SecTwoTwo} and \ref{SecContrConst}. Section \ref{SecTwoTwo} provides the complete algebraic characterization of approximate and exact controllability of \eqref{MainSyst} in dimension $2$ with two delays and a scalar input. Finally, Section \ref{SecContrConst} contains the results regarding controllability from zero to constant states. Some technical proofs are deferred to the appendix.

All the results in this paper also hold, with the same proofs, if one assumes $A = (A_1, \dotsc, A_N)$ to be in $\mathcal M_d(\mathbb R)^N$
and $B$ in $\mathcal M_{d, m}(\mathbb R)$, with the state $x(t)$ in $\mathbb R^d$ and the control $u(t)$ in $\mathbb R^m$. We choose complex-valued matrices, states, and controls for \eqref{MainSyst} in this paper following the approach of \cite{Chitour2016Stability}, which is mainly motivated by the fact that classical spectral conditions for difference equations such as those from \cite{Avellar1980Zeros, Henry1974Linear, Hale2002Strong} and \cite[Chapter 9]{Hale1993Introduction} are more naturally expressed in such a framework.

\section{Definitions and preliminary results}
\label{SecPrevResults}

In this section we provide the definitions of solutions of \eqref{MainSyst} and approximate and exact controllability in $L^2$, and recall the explicit representation formula for solutions of \eqref{MainSyst} and some elementary properties of $L^2$ controllability.

\subsection{Well-posedness and explicit representation of solutions}
\label{SecControlDelayWellPosed}

\begin{defi}
\label{DefiSolution}
Let $A = (A_1, \dotsc, A_N) \in \mathcal M_d(\mathbb C)^N$, $B \in \mathcal M_{d, m}(\mathbb C)$, $\Lambda = (\Lambda_1, \dotsc, \Lambda_N) \in (0, +\infty)^N$, $T > 0$, $x_0: \left[-\Lambda_{\max}, 0\right) \to \mathbb C^d$, and $u: [0, T] \to \mathbb C^m$. We say that $x: \left[-\Lambda_{\max}, T\right] \to \mathbb C^d$ is a \emph{solution} of \eqref{MainSyst} with initial condition $x_0$ and control $u$ if it satisfies \eqref{MainSyst} for every $t \in [0, T]$ and $x(t) = x_0(t)$ for $t \in \left[-\Lambda_{\max}, 0\right)$. In this case, for $t \in [0, T]$, we define $x_t: \left[-\Lambda_{\max}, 0\right) \to \mathbb C^d$ by $x_t = x(t + \cdot)|_{\left[-\Lambda_{\max}, 0\right)}$.
\end{defi}

This notion of solution, already used in \cite{Mazanti2017Relative} and similar to the one used in \cite{Chitour2016Stability}, requires no regularity on $x_0$, $u$, or $x$. Nonetheless, such a weak framework is enough to guarantee existence and uniqueness of solutions.

\begin{prop}
\label{PropControlDelayExistUnique}
Let $A = (A_1, \dotsc, A_N) \in \mathcal M_d(\mathbb C)^N$, $B \in \mathcal M_{d, m}(\mathbb C)$, $\Lambda = (\Lambda_1, \dotsc, \Lambda_N) \in (0, +\infty)^N$, $T > 0$, $x_0: \left[-\Lambda_{\max}, 0\right) \to \mathbb C^d$, and $u: [0, T] \to \mathbb C^m$. Then \eqref{MainSyst} admits a unique solution $x: \left[-\Lambda_{\max},\allowbreak T\right]\allowbreak \to \mathbb C^d$ with initial condition $x_0$ and control $u$.
\end{prop}

Proposition~\ref{PropControlDelayExistUnique} can be easily proved from \eqref{MainSyst}, which is already an explicit representation formula for the solution in terms of the initial condition and the control when $t < \Lambda_{\min}$. Its proof can be found in \cite[Proposition~2.2]{Mazanti2017Relative} and is very similar to that of \cite[Proposition~3.2]{Chitour2016Stability}.

We also recall that, as in \cite[Remark 3.4]{Chitour2016Stability} and \cite[Remark 2.3]{Mazanti2017Relative}, if $x_0, \widetilde x_0: \left[-\Lambda_{\max}, 0\right) \to \mathbb C^d$ and $u, \widetilde u: [0, T] \to \mathbb C^m$ are such that $x_0 = \widetilde x_0$ and $u = \widetilde u$ almost everywhere on their respective domains, then the solutions $x, \widetilde x: \left[-\Lambda_{\max}, T\right] \to \mathbb C^d$ of \eqref{MainSyst} associated respectively with $x_0$, $u$, and $\widetilde x_0$, $\widetilde u$, satisfy $x = \widetilde x$ almost everywhere on $\left[-\Lambda_{\max}, T\right]$. In particular, one still obtains existence and uniqueness of solutions of \eqref{MainSyst} for initial conditions in $L^p((-\Lambda_{\max}, 0), \mathbb C^d)$ and controls in $L^p((0, T), \mathbb C^m)$ for some $p \in [1, +\infty]$, and, in this case, solutions $x$ of \eqref{MainSyst} satisfy $x \in L^p(\left(-\Lambda_{\max}, T\right), \mathbb C^d)$, and hence $x_t \in L^p((-\Lambda_{\max}, 0), \mathbb C^d)$ for every $t \in [0, T]$.

In order to provide an explicit representation for the solutions of \eqref{MainSyst}, we first provide a recursive definition of the matrix coefficients $\Xi_{\mathbf n}$ appearing in such a representation.

\begin{defi}
\label{ControlDelayDefiXi}
For $A = (A_1, \dotsc, A_N) \in \mathcal M_d(\mathbb C)^N$ and $\mathbf n \in \mathbb Z^N$, we define the matrix $\Xi_{\mathbf n} \in \mathcal M_d(\mathbb C)$ inductively by
\begin{equation}
\label{EqControlDelayDefiXi}
\Xi_{\mathbf n} = 
\begin{dcases*}
0, & if $\mathbf n \in \mathbb Z^N \setminus \mathbb N^N$, \\
\id_d, & if $\mathbf n = 0$, \\
\sum_{k=1}^N A_k \Xi_{\mathbf n - e_k}, & if $\mathbf n \in \mathbb N^N \setminus \{0\}$. \\
\end{dcases*}
\end{equation}
\end{defi}

The explicit representation for the solutions of \eqref{MainSyst} used throughout the present paper is the one from \cite[Proposition~2.7]{Mazanti2017Relative}, which we state below.

\begin{prop}
\label{PropExplicit}
Let $A = (A_1, \dotsc, A_N) \in \mathcal M_d(\mathbb C)^N$, $B \in \mathcal M_{d, m}(\mathbb C)$, $\Lambda = (\Lambda_1, \dotsc, \Lambda_N) \in (0, +\infty)^N$, $T > 0$, $x_0: \left[-\Lambda_{\max}, 0\right) \to \mathbb C^d$, and $u: [0, T] \to \mathbb C^m$. The corresponding solution $x: \left[-\Lambda_{\max}, T\right] \to \mathbb C^d$ of \eqref{MainSyst} is given for $t \in [0, T]$ by
\begin{equation}
\label{ExplicitSol}
x(t) = \sum_{\substack{(\mathbf n, j) \in \mathbb N^N \times \llbracket 1, N\rrbracket \\ -\Lambda_j \leq t - \Lambda \cdot \mathbf n < 0}} \Xi_{\mathbf n - e_j} A_j x_0(t - \Lambda \cdot \mathbf n) + \sum_{\substack{\mathbf n \in \mathbb N^N \\ \Lambda \cdot \mathbf n \leq t}} \Xi_{\mathbf n} B u(t - \Lambda \cdot \mathbf n).
\end{equation}
\end{prop}

\begin{remk}
\label{RemkSemigroup}
Let $p \in [1, +\infty]$. For $t \geq 0$, we define $\Upsilon(t) \in \mathcal L(L^p((-\Lambda_{\max},\allowbreak 0),\allowbreak \mathbb C^d))$ by
\[
(\Upsilon(t)x_0)(s) = \sum_{\substack{(\mathbf n, j) \in \mathbb N^N \times \llbracket 1, N\rrbracket \\ -\Lambda_j \leq t + s - \Lambda \cdot \mathbf n < 0}} \Xi_{\mathbf n - e_j} A_j x_0(t + s - \Lambda \cdot \mathbf n).
\]
The operator $\Upsilon(t)$ maps an initial condition $x_0$ to the state $x_t = \left.x(t + \cdot)\right|_{(-\Lambda_{\max}, 0)}$, where $x$ is the solution of \eqref{MainSyst} at time $t$ with initial condition $x_0$ and control $0$. Using the fact that translations define continuous operators in $L^p$ when $p < \infty$, one proves that the family $\{\Upsilon(t)\}_{t \geq 0}$ is a strongly continuous semigroup in $L^p((-\Lambda_{\max}, 0), \mathbb C^d)$ for $p \in [1, +\infty)$ (see, e.g., \cite[Proposition~3.5]{Chitour2016Stability}).
\end{remk}

\subsection{Approximate and exact controllability in \texorpdfstring{$L^2$}{L2}}
\label{SecControlL2}

We now define the main notions we consider in this paper, namely the approximate and exact controllability of the state $x_t = \left.x(t + \cdot)\right|_{[-\Lambda_{\max}, 0)}$ of \eqref{MainSyst} in the function space $L^2((-\Lambda_{\max}, 0),\allowbreak \mathbb C^d)$. We start with the notations that will be used throughout the rest of the paper.

\begin{defi}
\label{DefiContrExactApprox}
Let $T \in (0, +\infty)$. We define the Hilbert spaces $\mathsf X$ and $\mathsf Y_T$ by $\mathsf X = L^2((-\Lambda_{\max}, 0),\allowbreak \mathbb C^d)$ and $\mathsf Y_T = L^2((0, T),\allowbreak \mathbb C^m)$ endowed with their usual inner products and associated norms.
\begin{enumerate}
\item We say that \eqref{MainSyst} is \emph{approximately controllable in time $T$} if, for every $x_0, y \in \mathsf X$ and $\varepsilon > 0$, there exists $u \in \mathsf Y_T$ such that the solution $x$ of \eqref{MainSyst} with initial condition $x_0$ and control $u$ satisfies $\norm{x_T - y}_{\mathsf X} < \varepsilon$.

\item We say that \eqref{MainSyst} is \emph{exactly controllable in time $T$} if, for every $x_0, y \in \mathsf X$, there exists $u \in \mathsf Y_T$ such that the solution $x$ of \eqref{MainSyst} with initial condition $x_0$ and control $u$ satisfies $x_T = y$.

\item We define the \emph{end-point operator} $E(T) \in \mathcal L(\mathsf Y_T, \mathsf X)$ by
\begin{equation}
\label{EqDefET}
(E(T)u)(t) = \sum_{\substack{\mathbf n \in \mathbb N^N \\ \Lambda \cdot \mathbf n \leq T + t}} \Xi_{\mathbf n} B u(T + t - \Lambda \cdot \mathbf n).
\end{equation}
\end{enumerate}
\end{defi}

Approximate or exact controllability in time $T$ implies the same kind of controllability for every time $T^\prime \geq T$, since one can take a control $u$ equal to zero in the interval $(0, T^\prime - T)$ and control the system from $T^\prime - T$ until $T^\prime$.

It follows immediately from Proposition~\ref{PropExplicit} that, for every $T > 0$, $x_0 \in \mathsf X$, and $u \in \mathsf Y_T$, the corresponding solution $x$ of \eqref{MainSyst} satisfies
\begin{equation}
\label{EqSolOperators}
x_T = \Upsilon(T) x_0 + E(T) u,
\end{equation}
where $\{\Upsilon(t)\}_{t \geq 0}$ is the semigroup defined in Remark \ref{RemkSemigroup}. Equation \eqref{EqSolOperators} allows one to immediately obtain the following classical characterization of approximate and exact controllability in terms of the operator $E(T)$ (cf.\ \cite[Lemma~2.46]{Coron2007Control}).

\begin{prop}
\label{PropETControl}
Let $T \in (0, +\infty)$.
\begin{enumerate}
\item\label{PropETControlApprox} System \eqref{MainSyst} is approximately controllable in time $T$ if and only if $\range E(T)$ is dense in $\mathsf X$.
\item\label{PropETControlExact} System \eqref{MainSyst} is exactly controllable in time $T$ if and only if $E(T)$ is surjective.
\end{enumerate}
\end{prop}

We recall in the next proposition the classical characterizations of approximate and exact controllability in terms of the adjoint operator $E(T)^\ast$, whose proofs can be found, e.g., in \cite[Section 2.3.2]{Coron2007Control}.

\begin{prop}
\label{PropControl}
Let $T \in (0, +\infty)$.
\begin{enumerate}
\item System \eqref{MainSyst} is approximately controllable in time $T$ if and only if $E(T)^\ast$ is injective, i.e., for every $x \in \mathsf X$,
\begin{equation}
\label{UniqueContinuation}
E(T)^\ast x = 0 \implies x = 0.
\end{equation}

\item\label{ControlExactETAst} System \eqref{MainSyst} is exactly controllable in time $T$ if and only if there exists $c > 0$ such that, for every $x \in \mathsf X$,
\begin{equation}
\label{ObservabilityInequality}
\norm{E(T)^\ast x}_{\mathsf Y_T}^2 \geq c \norm{x}_{\mathsf X}^2.
\end{equation}
\end{enumerate}
\end{prop}

Properties \eqref{UniqueContinuation} and \eqref{ObservabilityInequality} are called \emph{unique continuation property} and \emph{observability inequality}, respectively. In order to apply Proposition~\ref{PropControl}, we provide in the next lemma an explicit formula for $E(T)^\ast$, which can be immediately obtained from the definition of adjoint operator.

\begin{lemm}
\label{LemmETAdj}
Let $T \in (0, +\infty)$. The adjoint operator $E(T)^\ast \in \mathcal L(\mathsf X, \mathsf Y_T)$ is given by
\begin{equation}
\label{EqETAst}
(E(T)^\ast x)(t) = \sum_{\substack{\mathbf n \in \mathbb N^N \\ -\Lambda_{\max} \leq t - T + \Lambda \cdot \mathbf n < 0}} B^\ast \Xi_{\mathbf n}^\ast x(t - T + \Lambda \cdot \mathbf n).
\end{equation}
\end{lemm}

\begin{remk}
Exact controllability is preserved under small perturbations of $(A, B)$. This follows from Proposition~\ref{PropControl}\ref{ControlExactETAst} and the continuity of $E(T)^\ast$ with respect to the operator norm (which clearly results from \eqref{EqETAst}). However, exact controllability is not preserved for small perturbations of $\Lambda$ (cf.\ Theorem~\ref{MainTheo22}\ref{MainTheoA1A2Contr}\ref{MainTheoExact}). As regards approximate controllability, it is not preserved for small perturbations of $(A, B, \Lambda)$ (cf.\ Theorem~\ref{MainTheo22}\ref{MainTheoA1A2Contr}\ref{MainTheoApprox}, where $(A, B, \Lambda)$ is chosen such that the set $\mathcal S$ defined in that theorem is infinite).
\end{remk}

A useful result for studying approximate and exact controllability is the following lem\-ma, which states that such properties are preserved under linear change of coordinates, linear feedback, and changes of the time scale.

\begin{lemm}
\label{LemmControlDelayInvariant}
Let $T > 0$, $\lambda > 0$, $K_j \in \mathcal M_{m, d}(\mathbb C)$ for $j \in \llbracket 1, N\rrbracket$, $P \in \mathrm{GL}_d(\mathbb C)$, and consider the system
\begin{equation}
\label{MainSystChange}
x(t) = \sum_{j=1}^N P(A_j + B K_j) P^{-1} x\left(t - \frac{\Lambda_j}{\lambda}\right) + P B u(t).
\end{equation}
Then
\begin{enumerate}
\item\label{ApproxContrlIff} \eqref{MainSyst} is approximately controllable in time $T$ if and only if \eqref{MainSystChange} is approximately controllable in time $\frac{T}{\lambda}$;

\item\label{ExactContrlIff} \eqref{MainSyst} is exactly controllable in time $T$ if and only if \eqref{MainSystChange} is exactly controllable in time $\frac{T}{\lambda}$.
\end{enumerate}
\end{lemm}

\begin{prf}
Let us prove \ref{ApproxContrlIff}, the proof of \ref{ExactContrlIff} being similar. Assume that \eqref{MainSyst} is approximately controllable in time $T$ and take $x_0, y \in L^2((-\Lambda_{\max}/\lambda, 0),\allowbreak \mathbb C^d)$ and $\varepsilon > 0$. Let $\widetilde x_0, \widetilde{y} \in L^2((-\Lambda_{\max}, 0),\allowbreak \mathbb C^d)$ be given by $\widetilde x_0(t) = P^{-1} x_0(t/\lambda)$ and $\widetilde{y}(t) = P^{-1} y(t/\lambda)$. Since \eqref{MainSyst} is approximately controllable in time $T$, there exists $\widetilde u \in L^2((0, T), \mathbb C^m)$ such that the solution $\widetilde x$ of \eqref{MainSyst} with initial condition $\widetilde x_0$ and control $\widetilde u$ satisfies $\norm{\widetilde x_T - \widetilde{y}}_{\mathsf X} < \frac{\varepsilon \sqrt{\lambda}}{\abs{P}_{2}}$. Let $u \in L^2((0, T/\lambda), \mathbb C^m)$ and $x \in L^2((-\Lambda_{\max}/\lambda, T/\lambda), \mathbb C^d)$ be given by
\[
u(t) = \widetilde u(\lambda t) - \sum_{j=1}^N K_j \widetilde x(\lambda t - \Lambda_j), \qquad x(t) = P \widetilde x(\lambda t).
\]
A straightforward computation shows that $x$ is the solution of \eqref{MainSystChange} with initial condition $x_0$ and control $u$, and that $x_{T/\lambda}(t) = P \widetilde x_T(\lambda t)$ for $t \in (-\Lambda_{\max}/\lambda, 0)$. Hence $\norm{x_{T/\lambda} - y}_{L^2((-\Lambda_{\max}/\lambda, 0), \mathbb C^d)} \allowbreak < \varepsilon$, and thus \eqref{MainSystChange} is approximately controllable in time $\frac{T}{\lambda}$. The converse is proved in a similar way.
\end{prf}

\begin{remk}
\label{RemkGraphic}
One can provide a graphical representation for the operators $E(T)$ and $E(T)^\ast$ as follows. In a plane with coordinates $(\xi, \zeta)$, we draw in the domain $[0, T) \times [-\Lambda_{\max}, 0)$, for $\mathbf n \in \mathbb N^N$, the line segment $\sigma_{\mathbf n}$ defined by the equation $\zeta = \xi - T + \Lambda \cdot \mathbf n$ (see Figure \ref{FigET}). We associate with the line segment $\sigma_{\mathbf n}$ the matrix coefficient $\Xi_{\mathbf n} B$.

\begin{figure}[ht]
\centering
\begin{tikzpicture}[scale=2]

\draw[thick, arrows={-Stealth[scale=1.5]}] (-0.25, 0) -- (5.63906, 0) node[right] {$\xi$};
\draw[thick, arrows={-Stealth[scale=1.5]}] (0, -2.25000) -- (0, 0.25) node[left] {$\zeta$};

\draw[thick, dashed, color=darkgreen] (5.38906, -2.25000) -- (5.38906, 0.125) node[above, color=black] {$T$};
\draw[thick, dashed, color=darkgreen] (5.63906, -2.00000) -- (-0.125, -2.00000) node[left, color=black] {$-\Lambda_{\max}$};

\draw (5.38906, 0.00000) -- node[font=\scriptsize, midway, below, rotate=45] {$B$} (3.38906, -2.00000);
\draw (4.24746, 0.00000) -- node[font=\scriptsize, midway, below, rotate=45] {$\Xi_{(0, 0, 1)} B$} (2.24746, -2.00000);
\draw (3.10587, 0.00000) -- (1.10587, -2.00000);
\draw (1.96428, 0.00000) -- (0.00000, -1.96428);
\draw (0.82269, 0.00000) -- (0.00000, -0.82269);
\draw (3.77102, 0.00000) -- node[font=\scriptsize, midway, below, rotate=45] {$\Xi_{(0, 1, 0)} B$} (1.77102, -2.00000);
\draw (2.62943, 0.00000) -- (0.62943, -2.00000);
\draw (1.48784, 0.00000) -- (0.00000, -1.48784);
\draw (0.34624, 0.00000) -- (0.00000, -0.34624);
\draw (2.15299, 0.00000) -- (0.15299, -2.00000);
\draw (1.01140, 0.00000) -- (0.00000, -1.01140);
\draw (0.53495, 0.00000) -- (0.00000, -0.53495);
\draw (3.38906, 0.00000) -- node[font=\scriptsize, midway, below] {$\cdots$} (1.38906, -2.00000);
\draw (2.24746, 0.00000) -- (0.24746, -2.00000);
\draw (1.10587, 0.00000) -- (0.00000, -1.10587);
\draw (1.77102, 0.00000) -- (0.00000, -1.77102);
\draw (0.62943, 0.00000) -- (0.00000, -0.62943);
\draw (0.15299, 0.00000) -- (0.00000, -0.15299);
\draw (1.38906, 0.00000) -- (0.00000, -1.38906);
\draw (0.24746, 0.00000) -- (0.00000, -0.24746);

\draw[thick, dashed, color=red] (1.16076, -2.25000) -- (1.16076, 0.125) node[above, color=black] {$t$};
\draw[thick, dashed, color=blue] (5.63906, -1.76535) -- (-0.125, -1.76535) node[left, color=black] {$s$};
\end{tikzpicture}
\caption{Graphical representation for $E(T)$ and $E(T)^\ast$ in the case $N = 3$, $\Lambda_1 = 2$, $\Lambda_2 = \frac{\sqrt{5}+1}{2}$, $\Lambda_3 = \pi - 2$, and $T = e^2 - 2$. The matrix coefficients associated with the line segments $\sigma_{\mathbf n}$ are given in the picture for $\mathbf n = (0, 0, 0)$, $\mathbf n = (0, 0, 1)$, and $\mathbf n = (0, 1, 0)$.}
\label{FigET}
\end{figure}
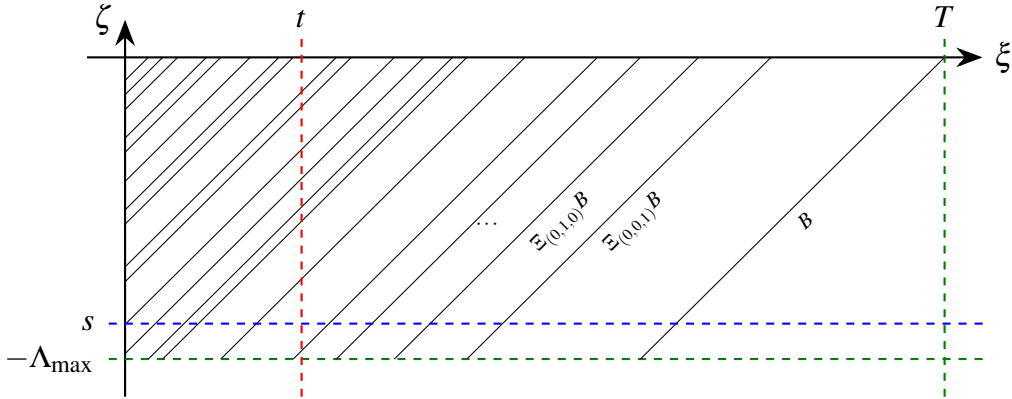

For $u \in \mathsf Y_T$, \eqref{EqDefET} can be interpreted as follows. For $s \in [-\Lambda_{\max}, 0)$, we draw the horizontal line $\zeta = s$. Each intersection between this line and a line segment $\sigma_{\mathbf n}$ gives one term in the sum for $(E(T)u)(s)$. This term consists of the matrix coefficient corresponding to the line $\sigma_{\mathbf n}$ multiplied by $u$ evaluated at the $\xi$-coordinate of the intersection point.

Similarly, for $x \in \mathsf X$, \eqref{EqETAst} can be interpreted as follows. For $t \in [0, T)$, we draw the vertical line $\xi = t$. As before, each intersection between this line and a line segment $\sigma_{\mathbf n}$ gives one term in the sum for $(E(T)^\ast x)(t)$. This term consists of the Hermitian transpose of the matrix coefficient corresponding to the line $\sigma_{\mathbf n}$ multiplied by $x$ evaluated at the $\zeta$-coordinate of the intersection point.
\end{remk}

\section{Controllability of systems with commensurable delays}
\label{SecCommensurable}

We consider in this section the problem of characterizing approximate and exact controllability of \eqref{MainSyst} in the case where the delays $\Lambda_1, \dotsc, \Lambda_N$ are commensurable. A classical procedure is to perform an augmentation of the state of the system to obtain an equivalent system with a single delay, whose controllability can be easily characterized using Kalman criterion for discrete-time linear control systems. For the sake of completeness, we detail such an approach in Lemma~\ref{LemmAAugm} and Proposition~\ref{PropCommensurable}. An important limitation of this technique is that it cannot be generalized to the case where $\Lambda_1, \dotsc, \Lambda_N$ are not assumed to be commensurable.

Thanks to Proposition~\ref{PropETControl}, another possible approach to the controllability of \eqref{MainSyst}, which will be extended to the case of incommensurable delays in Section \ref{SecTwoTwo}, is to consider the range of the operator $E(T)$. Following this approach, we characterize the operator $E(T)$ in Lemma~\ref{LemmMatrixET} in order to obtain a controllability criterion for \eqref{MainSyst} in Proposition~\ref{PropCommET}. It turns out that, in both criteria, controllability is equivalent to a full-rank condition on the \emph{same} matrix, as we prove in the main result of this section, Theorem~\ref{TheoCommens}.

\subsection{Kalman criterion based on state augmentation}

Let us first consider the augmentation of the state of \eqref{MainSyst}. The next lemma, whose proof is straightforward, provides the construction of the augmented state and the difference equation it satisfies.

\begin{lemm}
\label{LemmAAugm}
Let $T \in (0, +\infty)$, $u: [0, T] \to \mathbb C^m$, and suppose that $(\Lambda_1, \dotsc, \Lambda_N) = \lambda (k_1, \dotsc, k_N)$ with $\lambda > 0$ and $k_1, \dotsc, \allowbreak k_N \in \mathbb N^\ast$. Let $K = \max_{j \in \llbracket 1, N\rrbracket} k_j$.
\begin{enumerate}
\item If $x: [-\Lambda_{\max}, T] \to \mathbb C^d$ is the solution of \eqref{MainSyst} with initial condition $x_0: [-\Lambda_{\max}, 0) \to \mathbb C^d$, then the function $X: [-\lambda, T) \to \mathbb C^{K d}$ defined by
\begin{equation}
\label{AugmX}
X(t) = \begin{pmatrix}
x(t) \\
x(t - \lambda) \\
x(t - 2\lambda) \\
\vdots \\
x(t - (K-1)\lambda) \\
\end{pmatrix}
\end{equation}
satisfies
\begin{equation}
\label{AugmSyst}
X(t) = \widehat A X(t - \lambda) + \widehat B u(t),
\end{equation}
with $\widehat A$ and $\widehat B$ given by
\begin{equation}
\label{AugmMat}
\begin{aligned}
\widehat A & = \begin{pmatrix}
\widehat A_1 & \widehat A_2 & \widehat A_3 & \cdots & \widehat A_K \\
\id_d & 0 & 0 & \cdots & 0 \\
0 & \id_d & 0 & \cdots & 0 \\
\vdots & \vdots & \ddots & \ddots & \vdots \\
0 & 0 & \cdots & \id_d & 0 \\
\end{pmatrix} \in \mathcal M_{Kd}(\mathbb C), \quad \widehat B = \begin{pmatrix}
B \\
0 \\
0 \\
\vdots \\
0 \\
\end{pmatrix} \in \mathcal M_{K d, m}(\mathbb C),\\
\widehat A_k & = \sum_{\substack{j = 1 \\ k_j = k}}^N A_j \quad \text{ for } k \in \llbracket 1, K\rrbracket\text{ (in particular, $\widehat A_k = 0$ if $k_j \neq k$ for all $j \in \llbracket 1, N\rrbracket$)},
\end{aligned}
\end{equation}
and with initial condition $X_0: [-\lambda, 0) \to \mathbb C^{K d}$ given by
\begin{equation}
\label{AugmInitial}
X_0(t) = \begin{pmatrix}
x_0(t) \\
x_0(t - \lambda) \\
x_0(t - 2\lambda) \\
\vdots \\
x_0(t - (K-1)\lambda) \\
\end{pmatrix}.
\end{equation}

\item If $X: [-\lambda, T] \to \mathbb C^{K d}$ is the solution of \eqref{AugmSyst} with initial condition $X_0: [-\lambda, 0) \to \mathbb C^{K d}$, with $\widehat A$ and $\widehat B$ given by \eqref{AugmMat}, then the function $x: [-\Lambda_{\max}, T] \to \mathbb C^d$ defined by
\[
x(t) = 
\begin{dcases*}
\widehat C X(t), & if $t \in [0, T]$, \\
x_0(t), & if $t \in [-\Lambda_{\max}, 0)$,
\end{dcases*}
\]
is the solution of \eqref{MainSyst} with initial condition $x_0: [-\Lambda_{\max}, 0) \to \mathbb C^d$, where the matrix $\widehat C \in \mathcal M_{d, K d}(\mathbb C)$ is given by $\widehat C = \begin{pmatrix}\id_d & 0_{d, (K-1)d}\end{pmatrix}$ and $x_0$ is the unique function satisfying \eqref{AugmInitial} for every $t \in [-\lambda, 0)$.
\end{enumerate}
\end{lemm}

\begin{remk}
Lemma~\ref{LemmAAugm} considers solutions of \eqref{MainSyst} and \eqref{AugmSyst} in the sense of Definition \ref{DefiSolution}, i.e., with no regularity assumptions. However, one immediately obtains from \eqref{AugmX} that, for every $t \in [0, T]$, $x_t \in \mathsf X$ if and only if $X_t \in L^2((-\lambda, 0), \mathbb C^{Kd})$, and in this case $\norm{x_t}_{\mathsf X} = \norm{X_t}_{L^2((-\lambda, 0), \mathbb C^{Kd})}$.
\end{remk}

As an immediate consequence of Lemma~\ref{LemmAAugm}, we obtain the following criterion.

\begin{prop}
\label{PropCommensurable}
Let $T \in (0, +\infty)$ and suppose that $(\Lambda_1, \dotsc, \Lambda_N) = \lambda (k_1, \dotsc, k_N)$ with $\lambda > 0$ and $k_1, \dotsc, k_N \in \mathbb N^\ast$. Let $K = \max_{j \in \llbracket 1, N\rrbracket} k_j$ and define $\widehat A$ and $\widehat B$ from $A_1, \dotsc, A_N, B$ as in \eqref{AugmMat}. Then the following assertions are equivalent.
\begin{enumerate}
\item\label{CtrlCommensurableApprox} System \eqref{MainSyst} is approximately controllable in time $T$;
\item\label{CtrlCommensurableExact} System \eqref{MainSyst} is exactly controllable in time $T$;
\item\label{CtrlCommensurableKalman} $T \geq (\kappa+1) \lambda$, where $\kappa = \inf\left\{n \in \mathbb N \midsuchthat \rank\begin{pmatrix}\widehat B & \widehat A \widehat B & \widehat A^2 \widehat B & \cdots & \widehat A^n \widehat B\end{pmatrix} = K d\right\} \in \mathbb N \cup \{\infty\}$.
\end{enumerate}
\end{prop}

\begin{prf}
Notice first that the solution $X: [-\lambda, T] \to \mathbb C^{Kd}$ of \eqref{AugmSyst} with initial condition $X_0: [-\lambda, 0) \to \mathbb C^{Kd}$ and control $u: [0, T] \to \mathbb C^m$ is given by
\begin{equation}
\label{ExplicitSolAugmX}
X(t) = \widehat A^{1 + \floor{t/\lambda}} X_0\left(t - \left(1 + \floor{\frac{t}{\lambda}}\right) \lambda\right) + \sum_{n=0}^{\floor{t/\lambda}} \widehat A^n \widehat B u(t - n\lambda).
\end{equation}

We will prove that \ref{CtrlCommensurableExact} $\implies$ \ref{CtrlCommensurableApprox} $\implies$ \ref{CtrlCommensurableKalman} $\implies$ \ref{CtrlCommensurableExact}. The first implication is trivial due to the definitions of approximate and exact controllability. Suppose now that \ref{CtrlCommensurableApprox} holds, let $M = \floor{\frac{T}{\lambda}}$, $\rho = (M+1)\lambda - T > 0$, take $w \in \mathbb C^{Kd}$ and $\varepsilon > 0$, and write $w = \left(w_1^\transp, \dotsc, w_K^\transp\right)^\transp$ with $w_1, \dotsc, w_K \in \mathbb C^d$. Let $y \in \mathsf X$ be defined by the relations $y(t) = w_j$ for $t \in [-j\lambda, -(j-1)\lambda)$, $j \in \llbracket 1, K\rrbracket$. By \ref{CtrlCommensurableApprox}, there exists $u \in \mathsf Y_T$ such that the solution $x$ of \eqref{MainSyst} with zero initial condition and control $u$ satisfies $\norm{x_T - y}_{\mathsf X} < \rho\varepsilon$. Defining $X \in L^2((-\lambda, T), \mathbb C^{Kd})$ by \eqref{AugmX}, we obtain that $\norm{X_T - w}_{L^2((-\lambda, 0), \mathbb C^{Kd})} < \rho\varepsilon$. Using Lemma~\ref{LemmAAugm} and \eqref{ExplicitSolAugmX}, we obtain that
\[
\int_{T - \lambda}^{M\lambda} \abs{\sum_{n=0}^{M-1} \widehat A^n \widehat B u(t - n\lambda) - w}_{2}^2 dt \leq \int_{T - \lambda}^T \abs{\sum_{n=0}^{\floor{t/\lambda}} \widehat A^n \widehat B u(t - n\lambda) - w}_{2}^2 dt < \rho\varepsilon,
\]
and, in particular, there exists a set of positive measure $J \subset (T - \lambda, M\lambda)$ such that
\[
\abs{\sum_{n=0}^{M-1} \widehat A^n \widehat B u(t - n\lambda) - w}_{2}^2 < \varepsilon
\]
for $t \in J$. Hence, we have shown that, for every $w \in \mathbb C^{K d}$ and $\varepsilon > 0$, there exist $u_0, \dotsc, u_{M-1} \in \mathbb C^m$ such that $\abs{\sum_{n=0}^{M-1} \widehat A^n \widehat B u_n - w}_2^2 < \varepsilon$, which in particular implies that $M \geq 1$. This proves that the range of the matrix $\begin{pmatrix}\widehat B & \widehat A \widehat B & \widehat A^2 B & \cdots & \widehat A^{M-1} \widehat B\end{pmatrix} \in \mathcal M_{Kd, M m}(\mathbb C)$ is dense in $\mathbb C^{Kd}$, and hence is equal to $\mathbb C^{Kd}$, yielding $\kappa \leq M - 1$ by definition of $\kappa$. Thus $T \geq M \lambda \geq (\kappa + 1)\lambda$, which proves \ref{CtrlCommensurableKalman}.

Assume now that \ref{CtrlCommensurableKalman} holds. In particular, since $T < +\infty$, one has $\kappa \in \mathbb N$. We will prove the exact controllability of \eqref{MainSyst} in time $T_0 = (\kappa + 1)\lambda$, which implies its exact controllability in time $T$. Let $x_0, y \in \mathsf X$. Define $X_0, Y \in L^2((-\lambda, 0), \mathbb C^{Kd})$ from $x_0, y$ respectively as in \eqref{AugmInitial}. Let $C = \begin{pmatrix}\widehat B & \widehat A \widehat B & \cdots & \widehat A^\kappa \widehat B\end{pmatrix} \in \mathcal M_{Kd, (\kappa + 1) m}(\mathbb C)$, which, by \ref{CtrlCommensurableKalman}, has full rank, and thus admits a right inverse $C^\# \in \mathcal M_{(\kappa + 1) m, Kd}(\mathbb C)$. Let $u \in \mathsf Y_{T_0}$ be the unique function defined by the relation
\[
\begin{pmatrix}
u(t + (\kappa + 1)\lambda) \\
u(t + \kappa \lambda) \\
\vdots \\
u(t + \lambda)
\end{pmatrix} = C^\# \left(Y(t) - \widehat A^{\kappa + 1} X_0(t)\right) \qquad \text{ for almost every } t \in (-\lambda, 0).
\]
A straightforward computation shows, together with \eqref{ExplicitSolAugmX}, that the unique solution $X$ of \eqref{AugmSyst} with initial condition $X_0$ and control $u$ satisfies $X_{T_0} = Y$, and hence, by Lemma~\ref{LemmAAugm}, the unique solution of \eqref{MainSyst} with initial condition $x_0$ and control $u$ satisfies $x_{T_0} = y$, which proves \ref{CtrlCommensurableExact}.
\end{prf}

\begin{remk}
A first important consequence of Proposition~\ref{PropCommensurable} is that approximate and exact controllability are equivalent for systems with commensurable delays. As it follows from the results in Section \ref{SecTwoTwo}, this is no longer true when the commensurability hypothesis does not hold.
\end{remk}

\begin{remk}
\label{RemkCayleyHamilton}
It follows from Cayley--Hamilton theorem that $\kappa$ from Proposition~\ref{PropCommensurable} is either infinite or belongs to $\llbracket 0, Kd - 1\rrbracket$. In particular, \ref{CtrlCommensurableKalman} is satisfied for some $T \in (0, +\infty)$ if and only if the controllability matrix $\mathcal C(\widehat A, \widehat B) \in \mathcal M_{Kd, Kdm}(\mathbb C)$ has full rank. Moreover, condition \ref{CtrlCommensurableKalman} is satisfied for some $T \in (0, +\infty)$ if and only if it is satisfied for every $T \in [(\kappa + 1) \lambda, +\infty)$, and thus (approximate or exact) controllability in time $T \geq (\kappa + 1) \lambda$ is equivalent to (the same kind of) controllability in time $T = (\kappa + 1)\lambda$.
\end{remk}

\begin{remk}
\label{RemkKappaSingleInput}
When $m = 1$, it follows from the definition of $\kappa$ that $\kappa \geq K d - 1$ and thus, from Remark \ref{RemkCayleyHamilton}, $\kappa \in \{K d - 1, +\infty\}$. It follows that a system with a single input is either (approximately and exactly) controllable in time $T = d \Lambda_{\max}$ or not controllable in any time $T \in (0, +\infty)$.
\end{remk}

\begin{expl}
To illustrate the result from Proposition \ref{PropCommensurable} which relies on the state augmentation from Lemma \ref{LemmAAugm}, we provide the following example. Let $N = 2$ and $\Lambda = \lambda(1, 2)$ with $\lambda > 0$. Then \eqref{MainSyst} reads
\[x(t) = A_1 x(t - \lambda) + A_2 x(t - 2 \lambda) + B u(t),\]
and $K = 2$. The augmented matrices from \eqref{AugmMat} are given by
\[
\widehat A = 
\begin{pmatrix}
A_1 & A_2 \\
\id_d & 0 \\
\end{pmatrix}, \qquad \widehat B = 
\begin{pmatrix}
B \\ 0 \\
\end{pmatrix}.
\]
We now choose $d = 2$ and $A_1 = A_2 = \begin{pmatrix}0 & 1 \\ 0 & 0\\\end{pmatrix}$, $B = \begin{pmatrix}0 \\ 1 \\\end{pmatrix}$. Then
\[
\widehat A = 
\begin{pmatrix}
0 & 1 & 0 & 1 \\
0 & 0 & 0 & 0 \\
1 & 0 & 0 & 0 \\
0 & 1 & 0 & 0 \\
\end{pmatrix}, \qquad \widehat B = 
\begin{pmatrix}
0 \\
1 \\
0 \\
0 \\
\end{pmatrix}.
\]
It is easy to see that the condition from Proposition \ref{PropCommensurable}\ref{CtrlCommensurableKalman} is satisfied with $\kappa = K d - 1 = 3$ as soon as $T \geq 4 \lambda$. This value of $\kappa$ is in accordance with Remark \ref{RemkKappaSingleInput}.
\end{expl}

\subsection{\texorpdfstring{Controllability analysis through the range of $E(T)$}{Controllability analysis through the range of E(T)}}

We now turn to the characterization of the controllability of \eqref{MainSyst} using the operator $E(T)$ from \eqref{EqDefET} instead of the augmented system from Lemma~\ref{LemmAAugm}.

\begin{defi}
\label{DefiR1R2}
Let $T \in (0, +\infty)$ and suppose that $(\Lambda_1, \dotsc, \Lambda_N) = \lambda (k_1, \dotsc, k_N)$ with $\lambda > 0$ and $k_1, \dotsc, k_N \in \mathbb N^\ast$. Let $K = \max_{j \in \llbracket 1, N\rrbracket} k_j$, $M = \floor{\frac{T}{\lambda}}$, and $\delta = T - \lambda M \in [0, \lambda)$. We define $R_1 \in \mathcal L\left(\mathsf X, L^2((-\lambda, 0), \mathbb C^d)^K\right)$ and $R_2 \in \mathcal L\left(\mathsf Y_T, L^2((-\lambda, 0), \mathbb C^m)^M \times L^2((-\delta, 0), \mathbb C^m)\right)$ by
\[
\begin{aligned}
(R_1 x(t))_n & = x(t - (n-1)\lambda), & \quad & \text{ for } t \in (-\lambda, 0) \text{ and } n \in \llbracket 1, K\rrbracket, \\
(R_2 u(t))_n & = u(t + T - (n-1)\lambda), & & \text{ for } \left\{
\begin{aligned}
t \in (-\lambda, 0) & \text{ if } n \in \llbracket 1, M\rrbracket, \\
t \in (-\delta, 0) & \text{ if } n = M+1.
\end{aligned}\right.
\end{aligned}
\]
\end{defi}

It follows immediately from the definitions of $R_1$ and $R_2$ that these operators are unitary transformations. The operator $R_1$ allows to represent a function defined on $(-\Lambda_{\max}, 0)$ as a vector of $K$ functions defined on $(-\lambda, 0)$. The operator $R_2$ acts similarly on functions defined on $(0, T)$, with the interval of length $\delta < \lambda$ corresponding to the fact that $T$ is not necessarily an integer multiple of $\lambda$. In the next result, these transformations are used to provide a representation of $E(T)$ in terms of a block-Toeplitz matrix $C$ and a matrix $\Epsilon$.

\begin{lemm}
\label{LemmMatrixET}
Let $T \in (0, +\infty)$ and suppose that $(\Lambda_1, \dotsc, \Lambda_N) = \lambda (k_1, \dotsc, k_N)$ with $\lambda > 0$ and $k = (k_1, \dotsc,\allowbreak k_N) \in (\mathbb N^\ast)^N$. Let $K$, $M$, $\delta$, $R_1$, and $R_2$ be as in Definition \ref{DefiR1R2}. Then, for every $u \in L^2((-\lambda, 0),\allowbreak \mathbb C^m)^M \times L^2((-\delta, 0), \mathbb C^m)$,
\[
R_1 E(T) R_2^{-1} u = C P_1 u + \Epsilon P_2 u,
\]
where $P_1 \in \mathcal L\left(L^2((-\lambda, 0), \mathbb C^m)^M \times L^2((-\delta, 0), \mathbb C^m), L^2((-\lambda, 0), \mathbb C^m)^M\right)$ is the projection in the first $M$ coordinates, $P_2 \in \mathcal L\left(L^2((-\lambda, 0), \mathbb C^m)^M \times L^2((-\delta, 0), \mathbb C^m), L^2((-\lambda, 0), \mathbb C^m)\right)$ is the projection in the last coordinate composed with an extension by zero in the interval $(-\lambda, -\delta)$, and $C \in \mathcal M_{K d, M m}(\mathbb C), \Epsilon \in \mathcal M_{K d, m}(\mathbb C)$ are given by
\begin{equation}
\label{DefCEps}
\begin{aligned}
C & = \left(C_{j\ell}\right)_{j \in \llbracket 1, K\rrbracket, \ell \in \llbracket 1, M\rrbracket}, & \qquad & & C_{j\ell} & = \sum_{\substack{\mathbf n \in \mathbb N^N \\ k \cdot \mathbf n = \ell - j}} \Xi_{\mathbf n} B & & \text{ for } j \in \llbracket 1, K\rrbracket,\; \ell \in \llbracket 1, M\rrbracket, \\
\Epsilon & = \left(\Epsilon_j\right)_{j \in \llbracket 1, K\rrbracket}, & & & \Epsilon_j & = \sum_{\substack{\mathbf n \in \mathbb N^N \\ k \cdot \mathbf n = M+1 - j}} \Xi_{\mathbf n} B & & \text{ for } j \in \llbracket 1, K\rrbracket.
\end{aligned}
\end{equation}
\end{lemm}

\begin{prf}
Let $u \in \mathsf Y_T$ and extend $u$ by zero in the interval $(-\infty, 0)$. From \eqref{EqDefET} and Definition \ref{DefiR1R2}, we have that, for $j \in \llbracket 1, K\rrbracket$ and $t \in (-\lambda, 0)$,
\begin{align*}
\left(R_1 E(T) u(t)\right)_{j} & = \sum_{\substack{\mathbf n \in \mathbb N^N \\ \Lambda \cdot \mathbf n \leq T + t - (j - 1)\lambda}} \Xi_{\mathbf n} B u(t + T - \Lambda \cdot \mathbf n - (j - 1)\lambda) \displaybreak[0] \\
 & = \sum_{\substack{\mathbf n \in \mathbb N^N \\ k \cdot \mathbf n \leq \frac{T + t}{\lambda} - (j - 1)}} \Xi_{\mathbf n} B u(t + T - (k \cdot \mathbf n + j - 1)\lambda) \displaybreak[0] \\
 & = \sum_{\ell=1}^M \sum_{\substack{\mathbf n \in \mathbb N^N \\ k \cdot \mathbf n = \ell - j}} \Xi_{\mathbf n} B u(t + T - (\ell - 1)\lambda) + \sum_{\substack{\mathbf n \in \mathbb N^N \\ k \cdot \mathbf n = M+1 - j}} \Xi_{\mathbf n} B u(t + T - M \lambda) \displaybreak[0] \\
 & = \sum_{\ell=1}^M C_{j\ell} \left(P_1 R_2 u(t)\right)_{\ell} + \Epsilon_j \left(P_2 R_2 u(t)\right),
\end{align*}
which gives the required result.
\end{prf}

\begin{remk}
\label{RemkGraphCommens}
One can use the graphical representation of $E(T)$ from Remark \ref{RemkGraphic} to construct the matrices $C$ and $\Epsilon$ from Lemma~\ref{LemmMatrixET}. Indeed, when $(\Lambda_1, \dotsc, \Lambda_N) = \lambda (k_1, \dotsc, k_N)$ for some $\lambda > 0$ and $k_1, \dotsc, k_N \in \mathbb N^\ast$, one can consider a grid in $[0, T) \times [-\Lambda_{\max}, 0)$ defined by the horizontal lines $\zeta = -j\lambda$, $j \in \llbracket 1, K\rrbracket$, and by the vertical lines $\xi = T - (\ell - 1)\lambda$, $\ell \in \llbracket 1, M+1\rrbracket$, where $K = \max_{j \in \llbracket 1, N\rrbracket} k_j$ and $M = \floor{\frac{T}{\lambda}}$. This grid contains square cells $\mathsf S_{j\ell} = (T - \ell\lambda, T - (\ell - 1)\lambda) \times (-j\lambda, -(j-1)\lambda)$ for $j \in \llbracket 1, K\rrbracket$, $\ell \in \llbracket 1, M+1\rrbracket$, and rectangular cells $\mathsf R_j = (0, T - M\lambda) \times (-j\lambda, -(j-1)\lambda)$, the latter being empty when $T$ is an integer multiple of $\lambda$ (see Figure \ref{FigETMatrix}).

\begin{figure}[ht]
\centering
\begin{tikzpicture}[scale=2]

%\fill[red!15] (4.08727, 0) -- (2.28727, 0) -- (2.28727, -1.00000) -- (4.08727, -1.00000);
%\fill[blue!15] (0.08727, 0) -- (0, 0) -- (0, -1.20000) -- (0.08727, -1.20000);

\draw[thick, arrows={-Stealth[scale=1.5]}] (-0.25, 0) -- (4.33727, 0) node[right] {$\xi$};
\draw[thick, arrows={-Stealth[scale=1.5]}] (0, -2.25000) -- (0, 0.25) node[left] {$\zeta$};

\draw[color=gray] (4.08727, -2.00000) -- (4.08727, 0) node[above, color=black] {$T$};
\draw[color=gray] (3.88727, -2.00000) -- (3.88727, 0);
\draw[color=gray] (3.68727, -2.00000) -- (3.68727, 0);
\draw[color=gray] (3.48727, -2.00000) -- (3.48727, 0);
\draw[color=gray] (3.28727, -2.00000) -- (3.28727, 0);
\draw[color=gray] (3.08727, -2.00000) -- (3.08727, 0);
\draw[color=gray] (2.88727, -2.00000) -- (2.88727, 0);
\draw[color=gray] (2.68727, -2.00000) -- (2.68727, 0);
\draw[color=gray] (2.48727, -2.00000) -- (2.48727, 0);
\draw[color=gray] (2.28727, -2.00000) -- (2.28727, 0);
\draw[color=gray] (2.08727, -2.00000) -- (2.08727, 0);
\draw[color=gray] (1.88727, -2.00000) -- (1.88727, 0);
\draw[color=gray] (1.68727, -2.00000) -- (1.68727, 0);
\draw[color=gray] (1.48727, -2.00000) -- (1.48727, 0);
\draw[color=gray] (1.28727, -2.00000) -- (1.28727, 0);
\draw[color=gray] (1.08727, -2.00000) -- (1.08727, 0);
\draw[color=gray] (0.88727, -2.00000) -- (0.88727, 0);
\draw[color=gray] (0.68727, -2.00000) -- (0.68727, 0);
\draw[color=gray] (0.48727, -2.00000) -- (0.48727, 0);
\draw[color=gray] (0.28727, -2.00000) -- (0.28727, 0);
\draw[color=gray] (0.08727, -2.00000) -- (0.08727, 0);
\draw[color=gray] (4.08727, -0.00000) -- (0, -0.00000);
\draw[color=gray] (4.08727, -0.20000) -- (0, -0.20000);
\draw[color=gray] (4.08727, -0.40000) -- (0, -0.40000);
\draw[color=gray] (4.08727, -0.60000) -- (0, -0.60000);
\draw[color=gray] (4.08727, -0.80000) -- (0, -0.80000);
\draw[color=gray] (4.08727, -1.00000) -- (0, -1.00000);
\draw[color=gray] (4.08727, -1.20000) -- (0, -1.20000);
\draw[color=gray] (4.08727, -1.40000) -- (0, -1.40000);
\draw[color=gray] (4.08727, -1.60000) -- (0, -1.60000);
\draw[color=gray] (4.08727, -1.80000) -- (0, -1.80000);
\draw[color=gray] (4.08727, -2.00000) -- (0, -2.00000) node[left, color=black] {$-\Lambda_{\max}$};

\draw (4.08727, 0.00000) -- (2.08727, -2.00000);
\draw (2.68727, 0.00000) -- (0.68727, -2.00000);
\draw (1.28727, 0.00000) -- (0.00000, -1.28727);
\draw (3.48727, 0.00000) -- (1.48727, -2.00000);
\draw (2.08727, 0.00000) -- (0.08727, -2.00000);
\draw (0.68727, 0.00000) -- (0.00000, -0.68727);
\draw (2.88727, 0.00000) -- (0.88727, -2.00000);
\draw (1.48727, 0.00000) -- (0.00000, -1.48727);
\draw (0.08727, 0.00000) -- (0.00000, -0.08727);
\draw (2.28727, 0.00000) -- (0.28727, -2.00000);
\draw (0.88727, 0.00000) -- (0.00000, -0.88727);
\draw (1.68727, 0.00000) -- (0.00000, -1.68727);
\draw (0.28727, 0.00000) -- (0.00000, -0.28727);
\draw (1.08727, 0.00000) -- (0.00000, -1.08727);
\draw (0.48727, 0.00000) -- (0.00000, -0.48727);
\draw (2.08727, 0.00000) -- (0.08727, -2.00000);
\draw (0.68727, 0.00000) -- (0.00000, -0.68727);
\draw (1.48727, 0.00000) -- (0.00000, -1.48727);
\draw (0.08727, 0.00000) -- (0.00000, -0.08727);
\draw (0.88727, 0.00000) -- (0.00000, -0.88727);
\draw (0.28727, 0.00000) -- (0.00000, -0.28727);
\draw (0.08727, 0.00000) -- (0.00000, -0.08727);

\end{tikzpicture}
\caption{Graphical representation for $E(T)$ in the case $N = 3$, $\Lambda = \left(1, \frac{7}{10}, \frac{3}{10}\right)$, $\lambda = \frac{1}{10}$, and $T \in (2, 2 + \lambda)$.}
\label{FigETMatrix}
\end{figure}
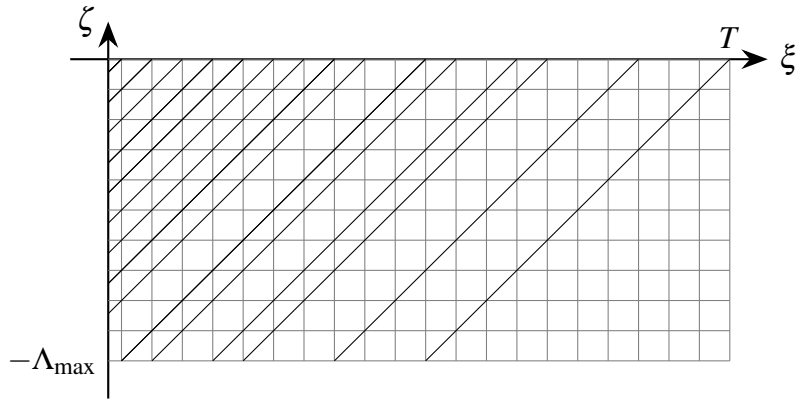

Consider the line segments $\sigma_{\mathbf n}$ from Remark \ref{RemkGraphic}. Due to the commensurability of the delays $\Lambda_1, \dotsc, \Lambda_N$, the intersection between each line segment $\sigma_{\mathbf n}$ and a square $\mathsf S_{j\ell}$ is either empty or equal to the diagonal of the square from its bottom-left to its top-right edge, and, similarly, the intersection between each $\sigma_{\mathbf n}$ and a rectangle $\mathsf R_j$ is either empty or equal to a line segment starting at the top-right edge of the rectangle. The matrix $C = \left(C_{j\ell}\right)_{j \in \llbracket 1, K\rrbracket,\; \ell \in \llbracket 1, M\rrbracket}$ can thus be constructed as follows. For $j \in \llbracket 1, K\rrbracket$ and $\ell \in \llbracket 1, M\rrbracket$, the matrix $C_{j\ell}$ is the sum over all $\mathbf n \in \mathbb N^N$ such that $\sigma_{\mathbf n}$ intersects the square $\mathsf S_{j\ell}$ of the matrix coefficients corresponding to $\sigma_{\mathbf n}$. Notice, in particular, that $C$ is a block-Toeplitz matrix, which is clear from its definition in \eqref{DefCEps}. Similarly, $\Epsilon = \left(\Epsilon_j\right)_{j \in \llbracket 1, K\rrbracket}$ is constructed by defining, for $j \in \llbracket 1, K\rrbracket$, $\Epsilon_j$ as the sum over all $\mathbf n \in \mathbb N^N$ such that $\sigma_{\mathbf n}$ intersects the rectangle $\mathsf R_{j}$ of the matrix coefficients corresponding to $\sigma_{\mathbf n}$. In the case $N = 3$, $\Lambda = \left(1, \frac{7}{10}, \frac{3}{10}\right)$, $\lambda = \frac{1}{10}$, and $T \in (2, 2 + \lambda)$, represented in Figure \ref{FigETMatrix}, the first $5d$ lines and $9m$ columns of the matrix $C$ are
\begin{equation*}
C = 
\begin{pmatrix}
B & 0 & 0 & \Xi_{(0, 0, 1)} B & 0 & 0 & \Xi_{(0, 0, 2)} B & \Xi_{(0, 1, 0)} B & 0 & \cdots \\
0 & B & 0 & 0 & \Xi_{(0, 0, 1)} B & 0 & 0 & \Xi_{(0, 0, 2)} B & \Xi_{(0, 1, 0)} B & \ddots \\
0 & 0 & B & 0 & 0 & \Xi_{(0, 0, 1)} B & 0 & 0 & \Xi_{(0, 0, 2)} B & \ddots \\
0 & 0 & 0 & B & 0 & 0 & \Xi_{(0, 0, 1)} B & 0 & 0 & \ddots \\
0 & 0 & 0 & 0 & B & 0 & 0 & \Xi_{(0, 0, 1)} B & 0 & \ddots \\
\vdots & \ddots & \ddots & \ddots & \ddots & \ddots & \ddots & \ddots & \ddots & \ddots \\
\end{pmatrix},
\end{equation*}
and the first $6d$ lines of $\Epsilon$ are
\[
\Epsilon = 
\begin{pmatrix}
\left(\Xi_{(2, 0, 0)} + \Xi_{(1, 1, 1)} + \Xi_{(0, 2, 2)}\right) B \\
\left(\Xi_{(1, 0, 3)} + \Xi_{(0, 1, 4)}\right) B \\
\Xi_{(0, 0, 6)} B \\
\left(\Xi_{(1, 1, 0)} + \Xi_{(0, 2, 1)}\right) B \\
\left(\Xi_{(1, 0, 2)} + \Xi_{(0, 1, 3)}\right) B \\
\Xi_{(0, 0, 5)} B \\
\vdots \\
\end{pmatrix}.
\]
\end{remk}

We now provide a controllability criterion for \eqref{MainSyst} in terms of the rank of $C$.

\begin{prop}
\label{PropCommET}
Let $T \in (0, +\infty)$ and suppose that $(\Lambda_1, \dotsc, \Lambda_N) = \lambda (k_1, \dotsc, k_N)$ with $\lambda > 0$ and $k_1, \dotsc, k_N \in \mathbb N^\ast$. Let $K$, $M$, and $C \in \mathcal M_{Kd, Mm}(\mathbb C)$ be as in Lemma~\ref{LemmMatrixET}. Then the following assertions are equivalent.
\begin{enumerate}
\item\label{CtrlCommETApprox} System \eqref{MainSyst} is approximately controllable in time $T$;
\item\label{CtrlCommETExact} System \eqref{MainSyst} is exactly controllable in time $T$;
\item\label{CtrlCommETKalman} The matrix $C$ has full rank.
\end{enumerate}
\end{prop}

\begin{prf}
The equivalence of \ref{CtrlCommETApprox} and \ref{CtrlCommETExact} has been proved in Proposition \ref{PropCommensurable}. Suppose that \ref{CtrlCommETExact} holds, which means, from Proposition~\ref{PropETControl}\ref{PropETControlExact}, that $E(T)$ is surjective. Since $R_1$ and $R_2$ are unitary transformations, Lemma~\ref{LemmMatrixET} shows that the operator $C P_1 + \Epsilon P_2: L^2((-\lambda, 0), \mathbb C^m)^M \times L^2((-\delta, 0), \mathbb C^m) \to L^2((-\lambda, 0), \mathbb C^d)^K$ is also surjective. Define the operator $\Pi \in \mathcal L\bigl(L^2((-\lambda, 0),\allowbreak \mathbb C^d)^K, L^2((-\lambda, -\delta), \mathbb C^d)^K\bigr)$ as the restriction to the non-empty interval $(-\lambda, -\delta)$, which is surjective. Thus $\Pi(C P_1 + \Epsilon P_2)$ is surjective, and one has, from the definition of $\Pi$ and $P_2$, that $\Pi \Epsilon P_2 = 0$, which shows that $\Pi C P_1$ is surjective. On the other hand, $(\Pi C P_1 u(t))_j = \sum_{\ell=1}^M C_{j\ell} u_\ell(t)$ for every $u \in L^2((-\lambda, 0), \mathbb C^m)^M \times L^2((-\delta, 0), \mathbb C^m)$, $j \in \llbracket 1, M\rrbracket$, and $t \in (-\lambda, -\delta)$, and hence $C$ has full rank, which proves \ref{CtrlCommETKalman}.

Suppose now that \ref{CtrlCommETKalman} holds. Notice that the matrix $C$ can be canonically identified with an operator, still denoted by $C$, in $\mathcal L\left(L^2((-\lambda, 0), \mathbb C^m)^M, L^2((-\lambda, 0), \mathbb C^d)^K\right)$, and such an operator is surjective. Defining $Q \in \mathcal L(L^2((-\lambda, 0), \mathbb C^m)^M, L^2((-\lambda, 0), \mathbb C^m)^M \times L^2((-\delta, 0), \mathbb C^m))$ by $Q u = (u, 0)$ for $u \in L^2((-\lambda, 0), \mathbb C^m)^M$, one has $C = (C P_1 + \Epsilon P_2) Q$, and thus $C P_1 + \Epsilon P_2$ is surjective, which yields, by Lemma~\ref{LemmMatrixET} and the fact that $R_1$ and $R_2$ are unitary transformations, that $E(T)$ is surjective. Thus, by Proposition~\ref{PropETControl}\ref{PropETControlExact}, \eqref{MainSyst} is exactly controllable in time $T$.
\end{prf}

\subsection{Comparison between Propositions \ref{PropCommensurable} and \ref{PropCommET}}

Propositions \ref{PropCommensurable} and \ref{PropCommET} provide two criteria for the controllability of \eqref{MainSyst} for commensurable delays $\Lambda_1, \dotsc, \Lambda_N$. The first one is obtained by the usual augmentation of the state and corresponds to a Kalman condition on the augmented matrices $\widehat A$ and $\widehat B$ from \eqref{AugmMat}, whereas the second one uses the characterizations of controllability in terms of the operator $E(T)$ from Proposition~\ref{PropETControl} in order to provide a criterion in terms of the matrix $C$ constructed from the matrix coefficients $\Xi_{\mathbf n} B$. It follows clearly from Propositions \ref{PropCommensurable} and \ref{PropCommET} that $C$ has full rank if and only if the matrix
\[
\begin{pmatrix}
\widehat B & \widehat A \widehat B & \widehat A^2 \widehat B & \cdots & \widehat A^{\floor{\frac{T}{\lambda}}-1} \widehat B
\end{pmatrix}
\]
has full rank. The main result of this section is that the two matrices coincide.

\begin{theo}
\label{TheoCommens}
Let $T \in (0, +\infty)$ and assume that $(\Lambda_1, \dotsc, \Lambda_N) = \lambda (k_1, \dotsc, k_N)$ with $\lambda > 0$ and $k_1, \dotsc, \allowbreak k_N \in \mathbb N^\ast$. Let $K$, $\widehat A$, $\widehat B$ be as in Proposition~\ref{PropCommensurable} and $M$, $C$ as in Proposition~\ref{PropCommET}. Then
\[C = \begin{pmatrix}\widehat B & \widehat A \widehat B & \widehat A^2 \widehat B & \cdots & \widehat A^{M-1} \widehat B\end{pmatrix}.\]
\end{theo}

\begin{prf}
For $j \in \llbracket 1, K\rrbracket$ and $\ell \in \llbracket 1, M\rrbracket$, let $C_{j\ell}$ be defined as in \eqref{DefCEps} and set $C_\ell = \left(C_{j\ell}\right)_{j \in \llbracket 1, K\rrbracket} \in \mathcal M_{Kd, m}(\mathbb C)$. We will prove the theorem by showing that $C_1 = \widehat B$ and that $C_{\ell + 1} = \widehat A C_\ell$ for $\ell \in \llbracket 1, M-1\rrbracket$. Let $k = (k_1, \dotsc, k_N)$.

By \eqref{DefCEps}, $C_{j1} = \sum_{\substack{\mathbf n \in \mathbb N^N \\ k \cdot \mathbf n = 1 - j}} \Xi_{\mathbf n} B$ for $j \in \llbracket 1, K\rrbracket$, and thus, since $\Xi_{\mathbf n} = 0$ for $\mathbf n \in \mathbb Z^N \setminus \mathbb N^N$, we obtain that $C_{j1} = 0$ for $j \in \llbracket 2, K\rrbracket$ and $C_{11} = \Xi_{0} B = B$, which shows that $C_1 = \widehat B$.

Let $\ell \in \llbracket 1, M-1\rrbracket$. For $j \in \llbracket 2, K\rrbracket$, we have $C_{j, \ell + 1} = \sum_{\substack{\mathbf n \in \mathbb N^N \\ k \cdot \mathbf n = \ell + 1 - j}} \Xi_{\mathbf n} B = C_{j-1, \ell} = \left(\widehat A C_{\ell}\right)_j$. Moreover, it follows from \eqref{EqControlDelayDefiXi} that
\begin{align*}
C_{1, \ell + 1} & = \sum_{\substack{\mathbf n \in \mathbb N^N \\ k \cdot \mathbf n = \ell}} \Xi_{\mathbf n} B = \sum_{\substack{\mathbf n \in \mathbb N^N \\ k \cdot \mathbf n = \ell}} \sum_{j=1}^N A_j \Xi_{\mathbf n - e_j} B {} = \sum_{\substack{\mathbf n \in \mathbb N^N \\ k \cdot \mathbf n = \ell}} \sum_{m=1}^K \sum_{\substack{j=1 \\ k_j = m}}^N A_j \Xi_{\mathbf n - e_j} B \displaybreak[0] \\
 & = \sum_{m=1}^K \sum_{\substack{\mathbf n \in \mathbb N^N \\ k \cdot \mathbf n = \ell}} \sum_{\substack{j=1 \\ k_j = m}}^N A_j \Xi_{\mathbf n - e_j} B = \sum_{m=1}^K \sum_{\substack{\mathbf n^\prime \in \mathbb N^N \\ k \cdot \mathbf n^\prime = \ell - m}} \sum_{\substack{j=1 \\ k_j = m}}^N A_j \Xi_{\mathbf n^\prime} B = \sum_{m=1}^K \widehat A_m C_{m\ell} = \left(\widehat A C_{\ell}\right)_1,
\end{align*}
where $\widehat A_m$ is defined as in \eqref{AugmMat}. Hence $\widehat A C_\ell = C_{\ell + 1}$, as required.
\end{prf}

\begin{remk}
Lemma~\ref{LemmMatrixET} shows that, when $\Lambda_1, \dotsc, \Lambda_N$ are commensurable, the operator $E(T)$ can be represented by the matrices $\Epsilon$ and $C$, and Proposition~\ref{PropCommET} shows that the controllability of \eqref{MainSyst} is encoded only in the matrix $C$. The representation of $E(T)$ by the matrix $C$ is also highlighted in Remark \ref{RemkGraphCommens}. Hence, the fact that $C$ coincides with the Kalman matrix $\begin{pmatrix}\widehat B & \widehat A \widehat B & \cdots & \widehat A^{M-1} \widehat B\end{pmatrix}$ for the augmented system \eqref{AugmSyst} shows that $E(T)$ generalizes the Kalman matrix for difference equations without the commensurability hypothesis on the delays.
\end{remk}

\begin{remk}
\label{RemkCn}
The main idea used here, namely the representation of $E(T)$ by the matrix $C$ in the commensurable case, is useful for the strategy we adopt in Section \ref{SecTwoTwo} to address the general case of incommensurable delays. Indeed, we characterize in Section \ref{SecTwoTwo} approximate and exact controllability through an operator $S$ which can be seen as a ``representation'' of $E(T)$ (see Definition \ref{DefinitionS}, Lemma~\ref{LemmS}, and Remark \ref{RemkGraphicSSAst}), and our strategy consists in approximating the delay vector $\Lambda$ by a sequence of commensurable delays $(\Lambda_n)_{n \in \mathbb N}$ and studying the asymptotic behavior of a corresponding sequence of matrices $(M_n)_{n \in \mathbb N}$, these matrices representing the operator $S$ in the same way as $C$ represents the operator $E(T)$.
\end{remk}

\section{Controllability of two-dimensional systems with two delays}
\label{SecTwoTwo}

In this section we investigate the controllability of \eqref{MainSyst} when the delays are not commensurable. The extension from the commensurable case is nontrivial, since the technique of state augmentation from Lemma~\ref{LemmAAugm} cannot be applied anymore and a deeper analysis of the operator $E(T)$ is necessary. In this section, we carry out such an analysis in the particular case $N = d = 2$ and $m = 1$, obtaining necessary and sufficient conditions for approximate and exact controllability. This simple-looking low-dimensional case already presents several non-trivial features that illustrate the difficulties stemming from the non-commensurability of the delays, including the fact that, contrarily to Propositions \ref{PropCommensurable} and \ref{PropCommET}, approximate and exact controllability are no longer equivalent.

Consider the difference equation
\begin{equation}
\label{Control22}
x(t) = A_1 x(t - \Lambda_1) + A_2 x(t - \Lambda_2) + B u(t),
\end{equation}
where $x(t) \in \mathbb C^2$, $u(t) \in \mathbb C$, $A_1, A_2 \in \mathcal M_2(\mathbb C)$, and $B \in \mathcal M_{2, 1}(\mathbb C)$, the latter set being canonically identified with $\mathbb C^2$. Without loss of generality, we assume that $\Lambda_1 > \Lambda_2$ and $B \neq 0$. The main result of this section is the following controllability criterion.

\begin{theo}
\label{MainTheo22}
Let $A_1, A_2 \in \mathcal M_2(\mathbb C)$, $B \in \mathcal M_{2, 1}(\mathbb C)$, $T \in (0, +\infty)$, and $(\Lambda_1, \Lambda_2) \in (0, +\infty)^2$ with $\Lambda_1 > \Lambda_2$ and $B \neq 0$.
\begin{enumerate}
\item\label{MainTheoA1BNonContr} If $\range A_1 \subset \range B$ or both pairs $(A_1, B)$, $(A_2, B)$ are not controllable, then \eqref{Control22} is neither approximately nor exactly controllable in time $T$.
\item\label{MainTheoA2BNonContr} If $\range A_1 \not\subset \range B$ and exactly one of the pairs $(A_1, B)$, $(A_2, B)$ is controllable, then the following are equivalent.
\begin{enumerate}
\item System \eqref{Control22} is approximately controllable in time $T$.
\item System \eqref{Control22} is exactly controllable in time $T$.
\item $T \geq 2 \Lambda_{1}$.
\end{enumerate}
\item\label{MainTheoA1A2Contr} If $(A_1, B)$ and $(A_2, B)$ are controllable, let $B^\perp \in \mathbb C^2$ be the unique vector such that $\det(B,\allowbreak B^\perp) \allowbreak = 1$ and $B^\transp B^\perp = 0$. Set
\begin{equation}
\label{DefAlphaBeta}
\beta = \frac{\det \mathcal C(A_1, B)}{\det \mathcal C(A_2, B)}, \qquad \alpha = \det \begin{pmatrix}B & (A_1 - \beta A_2) B^\perp\end{pmatrix}.
\end{equation}
Let $\mathcal S \subset \mathbb C$ be the set of all possible complex values of the expression $\beta + \alpha^{1 - \frac{\Lambda_{2}}{\Lambda_{1}}}$.
\begin{enumerate}
\item\label{MainTheoApprox} System \eqref{Control22} is approximately controllable in time $T$ if and only if $T \geq 2 \Lambda_{1}$ and $0 \notin \mathcal S$.
\item\label{MainTheoExact} System \eqref{Control22} is exactly controllable in time $T$ if and only if $T \geq 2 \Lambda_{1}$ and $0 \notin \overline{\mathcal S}$.
\end{enumerate}
\end{enumerate}
\end{theo}

\begin{remk}
The set $\mathcal S$ from case \ref{MainTheoA1A2Contr} is
\[
\mathcal S = \left\{\beta + \abs{\alpha}^{1 - \frac{\Lambda_2}{\Lambda_1}} e^{i (\theta + 2 k \pi)\left(1 - \frac{\Lambda_2}{\Lambda_1}\right)} \midsuchthat k \in \mathbb Z\right\},
\]
where $\theta \in (-\pi, \pi]$ is such that $\alpha = \abs{\alpha} e^{i\theta}$. Notice that $\mathcal S$ is a subset of the circle centered in $\beta$ with radius $\abs{\alpha}^{1 - \frac{\Lambda_2}{\Lambda_1}}$ (which reduces to a point when $\alpha = 0$). When $\frac{\Lambda_2}{\Lambda_1} \in \mathbb Q$, $\mathcal S$ is finite, $\overline{\mathcal S} = \mathcal S$, and one recovers the equivalence between exact and approximate controllability in time $T$ from Proposition 3.10. When $\frac{\Lambda_2}{\Lambda_1} \notin \mathbb Q$, $\mathcal S$ is a countable dense subset of the circle.
\end{remk}

\begin{remk}
In case \ref{MainTheoA1A2Contr}, approximate and exact controllability are characterized by the position of $0$ with respect of the subset $\mathcal S$ of $\mathbb C$, which is completely defined by $(A, B, \Lambda)$. It would be a striking result to generalize this fact to other values of $N$ and $d$. In this context, we believe that the strategy of our argument, as briefly described in Remark \ref{RemkCn}, is only suited for the case considered here, due to the difficulties in adapting to a more general case the reductions to normal forms from Remark \ref{RemkReduction}, the construction of the operator $S$ from Definition \ref{DefinitionS}, and the spectral study of the matrix $M$ from the appendix.
\end{remk}

The remainder of this section is dedicated to the proof of Theorem~\ref{MainTheo22}.

\subsection{Reduction to normal forms}
\label{SecReduction}

We start by characterizing the complex numbers $\alpha, \beta$ defined in \eqref{DefAlphaBeta}.

\begin{lemm}
\label{LemmAlphaEigenvalue}
Let $A_1, A_2 \in \mathcal M_2(\mathbb C)$, $B \in \mathcal M_{2, 1}(\mathbb C)$, assume that $(A_1, B)$ and $(A_2, B)$ are controllable, and let $\alpha, \beta$ be given by \eqref{DefAlphaBeta}. Let
\begin{equation}
\label{DefR}
R = \begin{pmatrix}
0 & 1 \\
-1 & 0 \\
\end{pmatrix}.
\end{equation}
Then $(A_1 - \beta A_2, B)$ is not controllable, $B$ is a right eigenvector of $A_1 - \beta A_2$, and $\alpha$ is an eigenvalue of $A_1 - \beta A_2$ associated with the left eigenvector $B^\transp R$.
\end{lemm}

\begin{prf}
By definition of $\beta$, one has $\det\mathcal C(A_1 - \beta A_2, B) = \det\begin{pmatrix} B & (A_1 - \beta A_2) B\\\end{pmatrix} = \det\begin{pmatrix} B & A_1 B\\\end{pmatrix}\allowbreak - \beta \det\begin{pmatrix}B & A_2 B\end{pmatrix} = 0$, and thus $(A_1 - \beta A_2, B)$ is not controllable. Moreover, since one has $\det\begin{pmatrix} B & (A_1 - \beta A_2) B\\\end{pmatrix} \allowbreak = 0$, the vectors $(A_1 - \beta A_2) B$ and $B$ are colinear, and thus $(A_1 - \beta A_2) B = \lambda B$ for some $\lambda \in \mathbb C$. Finally, notice that, for every $X, Y \in \mathcal M_{2, 1}(\mathbb C)$, $\det\begin{pmatrix}X & Y\\\end{pmatrix} = X^\transp R Y$, and thus, by definition of $\alpha$,
\begin{equation}
\label{Alpha}
\alpha = B^\transp R (A_1 - \beta A_2) B^\perp.
\end{equation}
Moreover, one has $B^\transp R B = \det\begin{pmatrix}B & B \\\end{pmatrix} = 0$ and $B^\transp R (A_1 - \beta A_2) B = \lambda B^\transp R B = 0$, which shows in particular that $B^\transp R (A_1 - \beta A_2) B = \alpha B^\transp R B$. Together with \eqref{Alpha}, this gives $B^\transp R (A_1 - \beta A_2) (a B^\perp \allowbreak + b B) = \alpha B^\transp R (a B^\perp + b B)$ for every $a, b \in \mathbb C$. Since $\{B, B^\perp\}$ is a basis of $\mathbb C^2$, this yields
\[B^\transp R (A_1 - \beta A_2) = \alpha B^\transp R,\]
and thus $B^\transp R$ is a left eigenvector of $A_1 - \beta A_2$ associated with the eigenvalue $\alpha$.
\end{prf}

We next show, thanks to the characterization of $\alpha, \beta$ from Lemma~\ref{LemmAlphaEigenvalue}, that $\alpha$ and $\beta$ are invariant under linear change of variables and linear feedbacks.

\begin{lemm}
\label{LemmAlphaBetaInvariant}
Let $A_1, A_2 \in \mathcal M_2(\mathbb C)$, $B \in \mathcal M_{2, 1}(\mathbb C)$, $P \in \mathrm{GL}_2(\mathbb C)$, $K_1, K_2 \in \mathcal M_{1, 2}(\mathbb C)$, and set
\[
\widetilde B = P B, \qquad \widetilde A_j = P (A_j + B K_j) P^{-1} \quad \text{ for } j \in \{1, 2\}.
\]
Suppose that $(A_1, B)$ and $(A_2, B)$ are controllable. Let $\alpha, \beta \in \mathbb C$ be defined by \eqref{DefAlphaBeta} and define $\widetilde\alpha, \widetilde\beta \in \mathbb C$ by
\[
\widetilde\beta = \frac{\det \mathcal C(\widetilde A_1, \widetilde B)}{\det \mathcal C(\widetilde A_2, \widetilde B)}, \qquad \widetilde\alpha = \det \begin{pmatrix}\widetilde B & (\widetilde A_1 - \widetilde \beta \widetilde A_2) \widetilde B^\perp\end{pmatrix},
\]
where $\widetilde B^\perp \in \mathbb C^2$ is the unique vector such that $\det(\widetilde B, \widetilde B^\perp) = 1$ and $\widetilde B^\transp \widetilde B^\perp = 0$. Then $\widetilde\alpha = \alpha$ and $\widetilde\beta = \beta$.
\end{lemm}

\begin{prf}
Since $\mathcal C(\widetilde A_j, \widetilde B) = P \mathcal C(A_j + B K_j, B)$ and $\det \mathcal C(A_j + B K_j, B) = \det \mathcal C(A_j, B)$ for $j \in \{1, 2\}$, one immediately obtains from the definitions of $\beta$ and $\widetilde\beta$ that $\widetilde\beta = \beta$. Let $R$ be given by \eqref{DefR}. By Lemma~\ref{LemmAlphaEigenvalue}, $\alpha$ is an eigenvalue of $A_1 - \beta A_2$ associated with the left eigenvector $B^\transp R$ and $\widetilde\alpha$ is an eigenvalue of $\widetilde A_1 - \widetilde\beta \widetilde A_2$ associated with the left eigenvector $\widetilde B^\transp R$. Using that $(P B)^\transp R (P B) = \det\begin{pmatrix}PB & PB\\\end{pmatrix} = 0$ and that $P^\transp R P = (\det P) R$, we get
\[
\begin{split}
\widetilde B^\transp R (\widetilde A_1 - \widetilde\beta \widetilde A_2) & = B^\transp P^\transp R P \left((A_1 - \beta A_2) + B(K_1 - \beta K_2)\right) P^{-1} = (\det P) B^\transp R (A_1 - \beta A_2) P^{-1} \\
 & = \alpha (\det P) B^\transp R P^{-1} = \alpha B^\transp P^\transp R P P^{-1} = \alpha \widetilde B^\transp R,
\end{split}
\]
which shows that $\widetilde\alpha = \alpha$.
\end{prf}

\begin{remk}
\label{RemkReduction}
It follows from Lemmas \ref{LemmControlDelayInvariant} and \ref{LemmAlphaBetaInvariant} that, in order to prove Theorem~\ref{MainTheo22}, it suffices to prove it for
\begin{equation}
\label{NormalForm}
A_j = \begin{pmatrix}
a_{j1} & a_{j2} \\
0 & 0 \\
\end{pmatrix} \text{ for } j \in \{1, 2\}, \qquad B = \begin{pmatrix}
0 \\
1 \\
\end{pmatrix}, \qquad (\Lambda_1, \Lambda_2) = (1, L)
\end{equation}
with $a_{jk} \in \mathbb C$ for $j, k \in \{1, 2\}$ and $L \in (0, 1)$. Indeed, given $A_1, A_2 \in \mathcal M_2(\mathbb C)$, $B \in \mathcal M_{2, 1}(\mathbb C)$, and $\Lambda_1, \Lambda_2 \in (0, +\infty)$ with $\Lambda_1 > \Lambda_2$, it suffices to take $\lambda = 1/\Lambda_1$, $P \in \mathrm{GL}_2(\mathbb C)$ satisfying $P B = \begin{pmatrix}0 & 1 \\\end{pmatrix}^\transp$, and, for $j \in \{1, 2\}$, $K_j \in \mathcal M_ {1, 2}(\mathbb C)$ such that $-K_j P^{-1}$ is equal to the second row of $P A_j P^{-1}$, and in this case $P(A_1 + B K_1)P^{-1}$, $P(A_2 + B K_2)P^{-1}$, $P B$, and $(\lambda \Lambda_1, \lambda \Lambda_2)$ are under the form \eqref{NormalForm}. 

Notice that $a_{j2} = -\det\mathcal C(A_j, B)$ for $j \in \{1, 2\}$, which implies that $a_{j2} =0$ if and only if $(A_j,B)$ is not controllable. 
Moreover, if $(A_j,B)$ for $j\in\{1,2\}$ is controllable, then $P \in \mathrm{GL}_2(\mathbb C)$ and $K_j \in \mathcal M_{1, 2}(\mathbb C)$ can be taken so that, in addition,  
 $P (A_j + B K_j) P^{-1}$ is  under the form 
 \begin{equation*}
A_j = \begin{pmatrix}
0 & 1 \\
0 & 0 \\
\end{pmatrix} \end{equation*}
 (see, e.g., \cite[Definition 5.1.5]{Sontag1998Mathematical}). Clearly, if both $(A_1,B)$ and $(A_2,B)$ are controllable, in general only one of the two matrices $A_1$ and $A_2$ can be put in such a normal form. 

We will thus prove Theorem~\ref{MainTheo22} for $(A_1, A_2, B, \Lambda_1, \Lambda_2) $ in one of the following normal forms:
\begin{equation}
\label{ReducGa1}
A_j = \begin{pmatrix}
a_{j1} & 0 \\
0 & 0 \\
\end{pmatrix} \text{ for } j \in \{1, 2\}, \qquad B = \begin{pmatrix}
0 \\
1 \\
\end{pmatrix}, \qquad (\Lambda_1, \Lambda_2) = (1, L),
\end{equation}
%\begin{equation}
%\label{ReducGa2}
%A_1 = \begin{pmatrix}
%a_{11} & 0 \\
%0 & 0 \\
%\end{pmatrix}, \qquad A_2 = \begin{pmatrix}
%0 & 1 \\
%0 & 0 \\
%\end{pmatrix}, \qquad B = \begin{pmatrix}
%0 \\
%1 \\
%\end{pmatrix}, \qquad (\Lambda_1, \Lambda_2) = (1, L),
%\end{equation}
%corresponding to part \ref{MainTheoA1BNonContr} of the theorem, and the normal forms
\begin{equation}
\label{ReducGb}
A_1 = \begin{pmatrix}
0 & 1 \\
0 & 0 \\
\end{pmatrix}, \qquad A_2 = \begin{pmatrix}
a_{21} & 0 \\
0 & 0 \\
\end{pmatrix}, \qquad B = \begin{pmatrix}
0 \\
1 \\
\end{pmatrix}, \qquad (\Lambda_1, \Lambda_2) = (1, L),
\end{equation}
and
\begin{equation}
\label{ReducGc}
A_1 = \begin{pmatrix}
a_{11} & a_{12} \\
0 & 0 \\
\end{pmatrix}, \qquad A_2 = \begin{pmatrix}
0 & 1 \\
0 & 0 \\
\end{pmatrix}, \qquad B = \begin{pmatrix}
0 \\
1 \\
\end{pmatrix}, \qquad (\Lambda_1, \Lambda_2) = (1, L).
\end{equation}
Part \ref{MainTheoA1BNonContr} in the statement of Theorem \ref{MainTheo22} corresponds to the normal forms \eqref{ReducGa1} and \eqref{ReducGc} in the case $a_{11} = a_{12} = 0$, \ref{MainTheoA2BNonContr} corresponds to \eqref{ReducGb} and \eqref{ReducGc} with $a_{11} \neq 0$ and $a_{12} = 0$, and \ref{MainTheoA1A2Contr} corresponds to \eqref{ReducGc} with $a_{12} \neq 0$. In the latter case, by a straightforward computation, one has $\alpha = a_{11}$ and $\beta = a_{12}$.
\end{remk}

\subsection{Proof of Theorem~\ref{MainTheo22}\ref{MainTheoA1BNonContr}}

\begin{prfsection}
In order to prove \ref{MainTheoA1BNonContr}, suppose first that $(A_1, B)$ and $(A_2, B)$ are not controllable. According to Remark \ref{RemkReduction}, we can assume that $A_1$, $A_2$, $B$, and $(\Lambda_1, \Lambda_2)$ are under the form \eqref{ReducGa1}. Hence one immediately computes
\[
\Xi_{\mathbf n} B =
\begin{dcases*}
B & if $\mathbf n = 0$, \\
0 & otherwise.
\end{dcases*}
\]
Then, for every $u \in \mathsf Y_T$ and $t \in (-1, 0)$, one has $(E(T) u)(t) = B u(T + t)$ if $T + t \geq 0$ and $(E(T) u)(t)\allowbreak = 0$ if $T + t < 0$. In particular, the range of $E(T)$ is contained in the set $L^2((-1, 0),\allowbreak \range B)$, which is not dense in $\mathsf X$. Hence the system is neither approximately nor exactly controllable in any time $T > 0$.

Consider now the case where $\range A_1 \subset \range B$. In particular, $(A_1, B)$ is not controllable, and one is left to consider the case where $(A_2, B)$ is controllable. In this case, the system can be brought down to the normal form \eqref{ReducGc} with $a_{11} = a_{12} = 0$. Hence
\begin{equation}
\label{CoeffsCase2}
\Xi_{\mathbf n} B =
\begin{dcases*}
\begin{pmatrix}0 \\ 1 \\\end{pmatrix} & if $\mathbf n = 0$, \\
\begin{pmatrix}1 \\ 0 \\\end{pmatrix} & if $\mathbf n = (0, 1)$, \\
0 & otherwise.
\end{dcases*}
\end{equation}
Then, for every $u \in \mathsf Y_T$, one has
\begin{equation}
\label{Case3ExplicitE}
(E(T) u)(t) = 
\begin{dcases*}
0 & if $-1 \leq T + t < 0$, \\
\begin{pmatrix}0 \\ u(T + t)\end{pmatrix} & if $0 \leq T + t < L$, \\
\begin{pmatrix} u(T + t - L) \\ u(T + t)\end{pmatrix} & if $T + t \geq L$.
\end{dcases*}
\end{equation}
If $T < 1 + L$, then, for every $u \in \mathsf Y_T$, the first component of $E(T) u$ vanishes in the non-empty interval $(-1, L - T)$, and hence the range of $E(T)$ is not dense in $\mathsf X$, which shows that the system is neither approximately nor exactly controllable in time $T < 1 + L$. If $T \geq 1 + L$, then, for every $u \in \mathsf Y_T$, if $x = E(T) u = (x_1, x_2)$, we have $x_1(t) = u(T + t - L)$ and $x_2(t) = u(T + t)$ for every $t \in (-1, 0)$, which implies that $x_2(t) = x_1(t + L)$ for $t \in (-1, -L)$. Hence the range of $E(T)$ is not dense in $\mathsf X$, which shows that the system is neither approximately nor exactly controllable in time $T \geq 1 + L$ either. This concludes the proof of \ref{MainTheoA1BNonContr}.

\subsection{Proof of Theorem~\ref{MainTheo22}\ref{MainTheoA2BNonContr}}

Concerning \ref{MainTheoA2BNonContr}, assume first that $(A_1, B)$ is controllable and $(A_2, B)$ is not controllable. According to Remark \ref{RemkReduction}, we can assume that $A_1$, $A_2$, $B$, and $(\Lambda_1, \Lambda_2)$ are under the form \eqref{ReducGb}. In this case, one has
\[
\Xi_{\mathbf n} B =
\begin{dcases*}
\begin{pmatrix}0 \\ 1 \\\end{pmatrix} & if $\mathbf n = 0$, \\
a_{21}^k \begin{pmatrix}1 \\ 0 \\\end{pmatrix} & if $\mathbf n = (1, k)$ and $k \in \mathbb N$, \\
0 & otherwise.
\end{dcases*}
\]
Then, for every $u \in \mathsf Y_T$, one has
\begin{equation}
\label{ETu-NewCase}
(E(T) u)(t) = 
\begin{dcases*}
0 & if $-1 \leq T + t < 0$, \\
\begin{pmatrix}0 \\ u(T + t)\end{pmatrix} & if $0 \leq T + t < 1$, \\
\begin{pmatrix} \displaystyle\sum_{k=0}^{\floor{\frac{T + t - 1}{L}}} a_{21}^k u(T + t - 1 - k L) \\ u(T + t)\end{pmatrix} & if $T + t \geq 1$.
\end{dcases*}
\end{equation}
If $T < 2$, then, for every $u \in \mathsf Y_T$, the first component of $E(T) u$ vanishes in the non-empty interval $(-1, 1-T)$, and hence the range of $E(T)$ is not dense in $\mathsf X$, which shows that the system is neither approximately nor exactly controllable in time $T < 2$. If $T \geq 2$, the system is exactly controllable. Indeed, take $x \in \mathsf X$ and write $x = (x_1, x_2)$. Define $u \in \mathsf Y_T$ by
\[
u(t) = 
\begin{dcases*}
x_2(t - T), & if $T - 1 \leq t < T$, \\
x_1(t - T + 1) - a_{21} x_1(t - T + 1 - L), & if $T - 2 + L \leq t < T - 1$, \\
x_1(t - T + 1), & if $T - 2 \leq t < T - 2 + L$, \\
0, & otherwise.
\end{dcases*}
\]
Then, for $t \in (-1, 0)$, one has $u(T + t) = x_2(t)$ and, for $k \in \left\llbracket 0, \floor{\frac{T + t - 1}{L}}\right\rrbracket$,
\[
u(T + t - 1 - k L) = 
\begin{dcases*}
x_1(t - k L) - a_{21} x_1(t - (k+1)L), & if $k \leq \frac{t+1}{L} - 1$, \\
x_1(t - k L), & if $k = \floor{\frac{t + 1}{L}}$, \\
0, & otherwise.
\end{dcases*}
\]
By \eqref{ETu-NewCase}, one immediately checks that $E(T) u = x$. Hence $E(T)$ is surjective, and thus the system is exactly controllable.

Assume now that $\range A_1 \not\subset \range B$, $(A_1, B)$ is not controllable, and $(A_2, B)$ is controllable. Thanks to Remark \ref{RemkReduction}, we can then assume that $A_1$, $A_2$, $B$, and $(\Lambda_1, \Lambda_2)$ are under the form \eqref{ReducGc} with $a_{11} \neq 0$ and $a_{12} = 0$. Hence
\[
\Xi_{\mathbf n} B =
\begin{dcases*}
\begin{pmatrix}0 \\ 1 \\\end{pmatrix} & if $\mathbf n = 0$, \\
a_{11}^k \begin{pmatrix}1 \\ 0 \\\end{pmatrix} & if $\mathbf n = (k, 1)$ and $k \in \mathbb N$, \\
0 & otherwise.
\end{dcases*}
\]
Then, for every $u \in \mathsf Y_T$, one has
\begin{equation}
\label{ETu-OtherNewCase}
(E(T) u)(t) = 
\begin{dcases*}
0 & if $-1 \leq T + t < 0$, \\
\begin{pmatrix}0 \\ u(T + t)\end{pmatrix} & if $0 \leq T + t < L$, \\
\begin{pmatrix} \displaystyle\sum_{k=0}^{\floor{T + t - L}} a_{11}^k u(T + t - k - L) \\ u(T + t)\end{pmatrix} & if $T + t \geq L$.
\end{dcases*}
\end{equation}
If $T < 1 + L$, \eqref{ETu-OtherNewCase} reduces to \eqref{Case3ExplicitE}, and the non-controllability of \eqref{Control22} follows as in \ref{MainTheoA1BNonContr}. If $1 + L \leq T < 2$, then, for every $u \in \mathsf Y_T$, if $x = E(T) u = (x_1, x_2)$, we have $x_1(t) = u(T + t - L)$ for $t \in (-1, 1 + L - T)$ and $x_2(t) = u(T + t)$ for $t \in (-1, 0)$, which implies that $x_2(t) = x_1(t + L)$ for $t \in (-1, 1 - T)$. As in the proof of \ref{MainTheoA1BNonContr}, the range of $E(T)$ is not dense in $\mathsf X$ and \eqref{Control22} is not controllable. To prove that \eqref{Control22} is exactly controllable when $T \geq 2$, take $x \in \mathsf X$ and write $x = (x_1, x_2)$. Define $u \in \mathsf Y_T$ by
\[
u(t) = 
\begin{dcases*}
x_2(t - T), & if $T - 1 \leq t < T$, \\
x_1(t - T + L), & if $T - 1 - L \leq t < T - 1$, \\
a_{11}^{-1} \left[x_1(t - T + 1 + L) - x_2(t - T + 1)\right] & if $T - 2 \leq t < T - 1 - L$, \\
0, & otherwise.
\end{dcases*}
\]
Then, for $t \in (-1, 0)$, one has $u(T + t) = x_2(t)$ and, for $k \in \left\llbracket 0, \floor{T + t - L}\right\rrbracket$,
\[
u(T + t - k - L) = 
\begin{dcases*}
x_2(t - k - L), & if $k = \floor{t - L + 1}$, \\
x_1(t - k), & if $t - L + 1 < k \leq t + 1$, \\
a_{11}^{-1} \left[x_1(t + 1 - k) - x_2(t + 1 - k - L)\right], & if $t + 1 < k \leq t - L + 2$, \\
0, & if $k > t - L + 2$.
\end{dcases*}
\]
If $t \in [-1, L - 1)$, then $\floor{t - L + 1} = -1$, $(t - L + 1, t + 1] \cap \mathbb N = \{0\}$, $(t + 1, t - L + 2] \cap \mathbb N = \emptyset$, and thus $\displaystyle\sum_{k=0}^{\floor{T + t - L}} a_{11}^k u(T + t - k - L) = x_1(t)$. If $t \in [L - 1, 0)$, then $\floor{t - L + 1} = 0$, $(t - L + 1, t + 1] \cap \mathbb N = \emptyset$, $(t + 1, t - L + 2] \cap \mathbb N = \{1\}$, and thus $\displaystyle\sum_{k=0}^{\floor{T + t - L}} a_{11}^k u(T + t - k - L) = x_2(t - L) + a_{11} a_{11}^{-1} \left[x_1(t) - x_2(t - L)\right] = x_1(t)$. It follows that $E(T) u = x$, proving that $E(T)$ is surjective and yielding the exact controllability of \eqref{Control22}.
\end{prfsection}

\subsection{Proof of Theorem~\ref{MainTheo22}\ref{MainTheoA1A2Contr}}
\label{SecProofCaseC}

In order to prove \ref{MainTheoA1A2Contr}, let us first provide explicit expressions for $E(T)$ and $E(T)^\ast$ when $A_1$, $A_2$, $B$, and $(\Lambda_1, \Lambda_2)$ are under the form \eqref{ReducGc}. In this case, one obtains, by a straightforward computation, that
\begin{equation}
\label{CoeffsExplicit}
\Xi_{(n, m)} B =
\begin{dcases*}
\begin{pmatrix}0 \\ 1 \\\end{pmatrix} & if $n = m = 0$, \\
\begin{pmatrix}\alpha^{n-1} \beta \\ 0 \\\end{pmatrix} & if $m = 0$, $n \geq 1$, \\
\begin{pmatrix}\alpha^{n} \\ 0 \\\end{pmatrix} & if $m = 1$, \\
0 & if $m \geq 2$,
\end{dcases*}
\end{equation}
where one uses that $\alpha = a_{11}$ and $\beta = a_{12}$. Hence, for every $u \in \mathsf Y_T$, $(E(T) u)(t) = 0$ for $T + t \in (-1, 0)$ and, for $T + t \geq 0$,
\begin{equation}
\label{Case4ExplicitE}
(E(T) u)(t) = \begin{pmatrix}\displaystyle \sum_{n=0}^{\floor{T + t - 1}} \alpha^{n} \beta u(T + t - n - 1) + \sum_{n=0}^{\floor{T + t - L}} \alpha^n u(T + t - n - L) \\ u(T + t)\end{pmatrix}.
\end{equation}
Moreover, for every $x = (x_1, x_2) \in \mathsf X$ and $t \in (-T, 0)$, one computes from \eqref{EqETAst} that 
\begin{equation}
\label{ETAstUtPlusT}
(E(T)^\ast x)(t + T) = \begin{dcases*}
x_2(t), & if $-L < t < 0$, \\
x_2(t) + x_1(t + L), & if $-1 < t < -L$, \\
\overline\alpha^{-\floor{t} - 2} \overline\beta x_1(\{t\} - 1) + \overline\alpha^{-\floor{t + L} - 1} x_1(\{t + L\} - 1), & if $t < -1$,
\end{dcases*}
\end{equation}
where we recall that $\{\xi\} = \xi - \floor{\xi}$ for $\xi \in \mathbb R$.

\subsubsection{\texorpdfstring{Case $T < 2 \Lambda_1$}{Case T < 2 Lambda1}}

\begin{prfsection}
Assume that $(A_1, B)$ and $(A_2, B)$ are controllable, in which case, according to Remark \ref{RemkReduction}, we can assume that $A_1$, $A_2$, $B$, and $(\Lambda_1, \Lambda_2)$ are under the form \eqref{ReducGc}, and thus $E(T)$ and $E(T)^\ast$ are given by \eqref{Case4ExplicitE} and \eqref{ETAstUtPlusT}, respectively.

If $T < 1 + L$, it follows from \eqref{Case4ExplicitE} that, for every $u \in \mathsf Y_T$, the first component of $E(T) u$ vanishes in the non-empty interval $(-1, L - T)$, and hence the system is neither approximately nor exactly controllable in time $T < 1 + L$.

For $1 + L \leq T < 2$, we will show that approximate controllability does not hold (and hence that exact controllability does not hold either) by showing that $E(T)^\ast$ is not injective. For $x = (x_1, x_2) \in \mathsf X$, it follows from \eqref{ETAstUtPlusT} that $E(T)^\ast x = 0$ in $\mathsf Y_T$ if and only if
\begin{equation}
\label{SystemAdjoint}
\left\{
\begin{aligned}
x_2(t) & = 0, & -L & < t < 0, \\
x_2(t) + x_1(t + L) & = 0, & -1 & < t < -L, \\
\overline\beta x_1(t + 1 - L) + x_1(t) & = 0, & -1 & < t < -1 + L, \\
\overline\beta x_1(t - L) + \overline\alpha x_1(t) & = 0, & \quad 1 + L - T & < t < 0.
\end{aligned}
\right.
\end{equation}
Since the first two equations of \eqref{SystemAdjoint} define $x_2$ uniquely in terms of $x_1$, showing that $E(T)^\ast x = 0$ for some nonzero function $x \in \mathsf X$ amounts to showing that there exists a nonzero function $y \in L^2((-1, 0), \mathbb C)$ such that
\begin{equation}
\label{EquationV}
\left\{
\begin{aligned}
\overline\beta y(t + 1 - L) + y(t) & = 0, & -1 & < t < -1 + L, \\
\overline\beta y(t - L) + \overline\alpha y(t) & = 0, & \quad 1 + L - T & < t < 0.
\end{aligned}
\right.
\end{equation}

Define $f: [-1, 0) \to [-1, 0)$ by $f(t) = t + 1 - L$ if $-1 \leq t < L - 1$ and $f(t) = t - L$ if $L - 1 \leq t < 0$; notice that $f$ is a translation by $1 - L$ modulo $1$. For $n \in \mathbb N$, set $t_n = f^n(-1)$ and let $K = \min\{n \in \mathbb N \suchthat f^{n+1}(-1) \in [-1, 1 - T)\}$. $K$ is clearly well-defined: if $L$ is rational, all orbits of $f$ are periodic and hence $K+1$ is upper bounded by the period of the orbit starting at $-1$, and, if $L$ is irrational, all orbits of $f$ are dense in $[-1, 0)$ and hence they intersect $[-1, 1-T)$ infinitely many times. Moreover, all the points $t_0, \dotsc, t_K$ are distinct. For $n \in \llbracket 0, K\rrbracket$, we define $\gamma_n \in \mathbb C$ inductively as follows. We set $\gamma_0 = 1$ and, for $n \in \llbracket 1, K\rrbracket$, we set $\gamma_n = -\frac{\gamma_{n-1}}{\beta}$ if $-1 \leq t_{n-1} < L - 1$ and $\gamma_n = -\frac{\alpha \gamma_{n-1}}{\beta}$ if $L - 1 \leq t_{n-1} < 0$.

Take $\delta > 0$ small enough such that all the intervals $(t_n, t_n + \delta)$, $n \in \llbracket 0, K\rrbracket$, are pairwise disjoint, contained in $(-1, 0)$, and do not contain any of the points $1 - T$, $L - 1$, $1 + L - T$, and $-L$ (these points may possibly be an extremity of the interval). Let $y \in L^2((-1, 0), \mathbb C)$ be defined by
\begin{equation}
\label{DefV}
y(t) = \sum_{n=0}^K \overline{\gamma_n} \chi_{(t_n, t_n + \delta)}(t).
\end{equation}
We claim that $y$ satisfies \eqref{EquationV}. Consider first the case $t \in (1 + L - T, 0)$, in which we have $f(t) = t - L$ since $(1 + L - T, 0) \subset [L - 1, 0)$. Since $f(1 + L - T) = 1 - T$ and $t_0 = -1$, it follows by construction of $\delta$ that $f(t) \notin (t_0, t_0 + \delta)$. If $t \notin \bigcup_{n=0}^K (t_n, t_n + \delta)$, then $f(t) \notin \bigcup_{n=0}^K (t_n, t_n + \delta)$; indeed, $f(t) \in (t_n, t_n + \delta)$ for some $n \in \llbracket 1, K\rrbracket$ implies immediately, by construction of $f$ and $\delta$, that $t \in (t_{n-1}, t_{n-1} + \delta)$. Hence, if $t \in (1 + L - T, 0) \setminus \bigcup_{n=0}^K (t_n, t_n + \delta)$, one immediately has that $y(t) = y(t - L) = 0$ and hence the second equation of \eqref{EquationV} is satisfied for such a $t$. Notice that $f(t_K) = t_{K+1} < 1 - T$, so that $t_K < 1 + L - T$, and thus, by construction of $\delta$, $(t_K, t_K + \delta) \cap (1 + L - T, 0) = \emptyset$. If $t \in (t_n, t_n + \delta)$ for some $n \in \llbracket 0, K - 1\rrbracket$, one has $t_n \in (1 + L - T, 0) \subset [L - 1, 0)$ by construction of $\delta$ and $f(t) \in (t_{n+1}, t_{n+1} + \delta)$, which shows, by the construction of $(\gamma_n)_{n=0}^K$, that
\[
\overline\alpha y(t) + \overline\beta y(t - L) = \overline\alpha \overline{\gamma_n} + \overline\beta \overline{\gamma_{n+1}} = 0.
\]
Hence the second equation of \eqref{EquationV} is satisfied for every $t \in (1 + L - T, 0)$.

Consider now the case $t \in (-1, L - 1)$, in which we have $f(t) = t + 1 - L$. Since $f^{-1}(t_0, t_0 + \delta) = (L-1, L-1+\delta)$, one has $f(t) \notin (t_0, t_0 + \delta)$. Again, the same argument as before shows that, if $t \notin \bigcup_{n=0}^K (t_n, t_n + \delta)$, then $f(t) \notin \bigcup_{n=0}^K (t_n, t_n + \delta)$, and thus, for such a $t$, $y(t) = y(t + 1 - L) = 0$ and the first equation of \eqref{EquationV} is satisfied. Since $f(t_K) = t_{K + 1} \in [-1, 1 - T)$, one has $t_K \in [L - 1, 1 + L - T)$, and hence $(t_K, t_K + \delta) \cap (-1, L-1) = \emptyset$. If $t \in (t_n, t_n + \delta) \cap (-1, L-1)$ for some $n \in \llbracket 0, K-1\rrbracket$, one has $t_n \in (-1, L-1)$ and $f(t) \in (t_{n+1}, t_{n+1} + \delta)$, which shows, by the construction of $(\gamma_n)_{n=0}^K$, that
\[
\overline\beta y(t + 1 - L) + y(t) = \overline\beta \overline{\gamma_{n+1}} + \overline{\gamma_{n}} = 0.
\]
Hence the first equation of \eqref{EquationV} is satisfied for every $t \in (-1, L - 1)$. Thus $E(T)^\ast$ is not injective, yielding that approximate controllability does not hold.
\end{prfsection}

\begin{remk}
One can modify the above construction to obtain a smooth function $x \in \mathcal C^\infty_0([-1,\allowbreak 0),\allowbreak \mathbb C^2)$ in the kernel of $E(T)^\ast$, simply by replacing the characteristic functions $\chi_{(t_n, t_n + \delta)}$ in \eqref{DefV} by $\varphi(\cdot - t_n)$ for a certain $\mathcal C^\infty$ function $\varphi$ compactly supported in $(0, \delta)$.
\end{remk}

\subsubsection{\texorpdfstring{Case $T \geq 2 \Lambda_1$}{Case T >= 2 Lambda1}}

The next lemma shows that one can reduce the proof of Theorem~\ref{MainTheo22}\ref{MainTheoA1A2Contr} in the case $T \geq 2 \Lambda_1$ to the case $T = 2 \Lambda_1$.

\begin{lemm}
\label{LemmTime2}
Let $A_1, A_2 \in \mathcal M_2(\mathbb C)$, $B \in \mathcal M_{2, 1}(\mathbb C)$, and $(\Lambda_1, \Lambda_2) \in (0, +\infty)^2$ with $\Lambda_1 > \Lambda_2$, and assume that $(A_1, B)$ and $(A_2, B)$ are controllable. Then the following assertions hold.
\begin{enumerate}
\item System \eqref{Control22} is approximately controllable in some time $T \geq 2 \Lambda_1$ if and only if it is approximately controllable in time $T = 2 \Lambda_1$.
\item System \eqref{Control22} is exactly controllable in some time $T \geq 2 \Lambda_1$ if and only if it is exactly controllable in time $T = 2 \Lambda_1$.
\end{enumerate}
\end{lemm}

\begin{prf}
Thanks to Remark \ref{RemkReduction}, it suffices to consider the case where $A_1$, $A_2$, $B$, and $(\Lambda_1, \Lambda_2)$ are given by \eqref{ReducGc}, in which case $E(T)^\ast$ is given by \eqref{ETAstUtPlusT}.

It is trivial that approximate controllability in $T = 2$ implies approximate controllability for larger time. To prove the converse, suppose that the system is approximately controllable in time $T \geq 2$ and take $x \in \mathsf X$ such that $E(2)^\ast x = 0$ in $\mathsf Y_2$. Thanks to \eqref{ETAstUtPlusT}, this means that, for almost every $t \in (-2, 0)$,
\[
\left\{
\begin{aligned}
x_2(t) & = 0, & \quad & \text{ if } -L < t < 0, \\
x_2(t) + x_1(t + L) & = 0, & & \text{ if } -1 < t < -L, \\
\overline\alpha^{-\floor{t} - 2} \overline\beta x_1(\{t\} - 1) + \overline\alpha^{-\floor{t + L} - 1} x_1(\{t + L\} - 1) & = 0, & & \text{ if } -2 < t < -1.
\end{aligned}
\right.
\]
Multiplying the last equation by $\overline\alpha^k$ for $k \in \mathbb N^\ast$ shows that, for almost every $t \in (-\infty, 0)$,
\[
\left\{
\begin{aligned}
x_2(t) & = 0, & \quad & \text{ if } -L < t < 0, \\
x_2(t) + x_1(t + L) & = 0, & & \text{ if } -1 < t < -L, \\
\overline\alpha^{-\floor{t} - 2} \overline\beta x_1(\{t\} - 1) + \overline\alpha^{-\floor{t + L} - 1} x_1(\{t + L\} - 1) & = 0, & & \text{ if } t < -1.
\end{aligned}
\right.
\]
In particular, $E(T)^\ast x = 0$ in $\mathsf Y_T$, and thus $x = 0$ in $\mathsf X$, which shows the approximate controllability in time $2$.

Concerning exact controllability, it is trivial that exact controllability in $T = 2$ implies exact controllability for larger time. To prove the converse, it suffices to show that, for every $T \geq 2$, there exists $C_T > 0$ such that, for every $x \in \mathsf X$,
\[
\norm{E(T)^\ast x}_{\mathsf Y_T}^2 \leq C_T \norm{E(2)^\ast x}_{\mathsf Y_2}^2.
\]
Let $T \geq 2$, $x = (x_1, x_2) \in \mathsf X$. Since the right-hand side of \eqref{ETAstUtPlusT} does not depend on $T$, one obtains that, for $t \in (-2, 0)$, $(E(T)^\ast x)(t + T) = (E(2)^\ast x)(t + 2)$. Hence
\begin{align*}
\norm{E(T)^\ast x}_{\mathsf Y_T}^2 & = \int_0^T \abs{(E(T)^\ast x)(t)}^2 dt = \int_{-T}^0 \abs{(E(T)^\ast x)(t + T)}^2 dt \displaybreak[0] \\
 & = \norm{E(2)^\ast x}_{\mathsf Y_2}^2 + \int_{-T}^{-2} \abs{(E(T)^\ast x)(t + T)}^2 dt \displaybreak[0] \\
 & = \norm{E(2)^\ast x}_{\mathsf Y_2}^2 + \int_{-T}^{-2} \abs{\overline\alpha^{-\floor{t} - 2} \overline\beta x_1(\{t\} - 1) + \overline\alpha^{-\floor{t + L} - 1} x_1(\{t + L\} - 1)}^2 dt \displaybreak[0] \\
 & \leq \norm{E(2)^\ast x}_{\mathsf Y_2}^2 \\
 & \hphantom{\leq} {} + \sum_{k=1}^{\ceil{T} - 2} \int_{-(k+2)}^{-(k+1)} \abs{\overline\alpha^{-\floor{t} - 2} \overline\beta x_1(\{t\} - 1) + \overline\alpha^{-\floor{t + L} - 1} x_1(\{t + L\} - 1)}^2 dt \displaybreak[0] \\
 & = \norm{E(2)^\ast x}_{\mathsf Y_2}^2 \\
 & \hphantom{=} {} + \sum_{k=1}^{\ceil{T} - 2} \abs{\alpha}^k \int_{-2}^{-1} \abs{\overline\alpha^{-\floor{t} - 2} \overline\beta x_1(\{t\} - 1) + \overline\alpha^{-\floor{t + L} - 1} x_1(\{t + L\} - 1)}^2 dt \displaybreak[0] \\
 & \leq \norm{E(2)^\ast x}_{\mathsf Y_2}^2 \sum_{k=0}^{\ceil{T} - 2} \abs{\alpha}^k,
\end{align*}
and one can thus conclude the proof by taking $C_T = \sum_{k=0}^{\ceil{T} - 2} \abs{\alpha}^k$.
\end{prf}

In order to study the controllability of \eqref{Control22} in the case $T = 2 \Lambda_1$, we introduce the following operator.

\begin{defi}
\label{DefinitionS}
We define the Hilbert space $\mathsf Z$ by $\mathsf Z = L^2((-1, 0), \mathbb C)$. Let $\alpha, \beta \in \mathbb C$. We define the bounded linear operator $S \in \mathcal L(\mathsf Z)$ by
\begin{equation}
\label{DefS}
S x(t) = 
\begin{dcases*}
\beta x(t) + \alpha x(t - L), & if $L - 1 < t < 0$, \\
\beta x(t) + x(t - L + 1), & if $-1 < t < L - 1$.
\end{dcases*}
\end{equation}
\end{defi}

By a straightforward computation, one obtains that the adjoint operator $S^\ast \in \mathcal L(\mathsf Z)$ is given, for $x \in \mathsf Z$, by
\begin{equation}
\label{SAst}
S^\ast x(t) = 
\begin{dcases}
\overline\beta x(t) + x(t + L - 1) & \text{ if } -L < t < 0, \\
\overline\beta x(t) + \overline\alpha x(t + L) & \text{ if } -1 < t < -L.
\end{dcases}
\end{equation}
The operators $S$ and $S^\ast$ allow one to characterize approximate and exact controllability for \eqref{Control22}, as shown in the next lemma.

\begin{lemm}
\label{LemmS}
Let $A_1, A_2 \in \mathcal M_2(\mathbb C)$, $B \in \mathcal M_{2, 1}(\mathbb C)$, and $(\Lambda_1, \Lambda_2) \in (0, +\infty)^2$ with $\Lambda_1 > \Lambda_2$, and assume that $(A_1, B)$ and $(A_2, B)$ are controllable. Then the following assertions hold.
\begin{enumerate}
\item\label{Syst22SInj} System \eqref{Control22} is approximately controllable in some time $T \geq 2 \Lambda_1$ if and only if $S^\ast$ is injective.

\item\label{Syst22SAstSurj} System \eqref{Control22} is exactly controllable in some time $T \geq 2 \Lambda_1$ if and only if $S$ is surjective or, equivalently, if there exists $c > 0$ such that $\norm{S^\ast x}_{\mathsf Z} \geq c \norm{x}_{\mathsf Z}$ for every $x \in \mathsf Z$.
\end{enumerate}
\end{lemm}

\begin{prf}
Thanks to Remark \ref{RemkReduction}, we can assume that $A_1$, $A_2$, $B$, and $(\Lambda_1, \Lambda_2)$ are under the form \eqref{ReducGc}, in which case $E(T)$ and $E(T)^\ast$ are given respectively by \eqref{Case4ExplicitE} and \eqref{ETAstUtPlusT}.

Let us first prove \ref{Syst22SInj}. Combining Lemma~\ref{LemmTime2} and Proposition~\ref{PropControl}, one obtains that \eqref{Control22} is approximately controllable in some time $T \geq 2$ if and only if $E(2)^\ast$ is injective. Thanks to \eqref{ETAstUtPlusT} and \eqref{SAst}, $x = (x_1, x_2) \in \mathsf X$ satisfies $E(2)^\ast x = 0$ if and only if
\begin{equation}
\label{SystToShowControllability}
\left\{
\begin{aligned}
x_2(t) & = 0, & \quad & \text{ if } -L < t < 0, \\
x_2(t) & = -x_1(t + L) , & & \text{ if } -1 < t < -L, \\
S^\ast x_1(t) & = 0, & & \text{ if } -1 < t < 0.
\end{aligned}
\right.
\end{equation}
Assume that $E(2)^\ast$ is injective and let $w \in \mathsf Z$ be such that $S^\ast w = 0$. Defining $x = (x_1, x_2) \in \mathsf X$ by $x_1 = w$, $x_2(t) = 0$ for $t \in (-L, 0)$, and $x_2(t) = -w(t + L)$ for $t \in (-1, -L)$, one obtains from \eqref{SystToShowControllability} that $E(2)^\ast x = 0$, which implies that $x = 0$ and hence $w = 0$, yielding the injectivity of $S^\ast$. Assume now that $S^\ast$ is injective and let $x = (x_1, x_2) \in \mathsf X$ be such that $E(2)^\ast x = 0$. Then, by the third equation of \eqref{SystToShowControllability}, one has $S^\ast x_1 = 0$, which shows that $x_1 = 0$, and thus the first two equations of \eqref{SystToShowControllability} show that $x_2 = 0$, yielding the injectivity of $E(2)^\ast$. Hence the injectivity of $E(2)^\ast$ is equivalent to that of $S^\ast$.

Let us now prove \ref{Syst22SAstSurj}. Combining Lemma~\ref{LemmTime2} and Proposition~\ref{PropETControl}, one obtains that \eqref{Control22} is exactly controllable in some time $T \geq 2$ if and only if $E(2)$ is surjective. Thanks to \eqref{Case4ExplicitE}, one has, for $u \in \mathsf Y_2$,
\begin{equation}
\label{E2u}
(E(2) u)(t) = 
\begin{dcases*}
\begin{pmatrix}
\beta u(t + 1) + \alpha u(t + 1 - L) + u(t + 2 - L) \\
u(t + 2) \\
\end{pmatrix}, & if $L - 1 < t < 0$, \\
\begin{pmatrix}
\beta u(t + 1) + u(t + 2 - L) \\
u(t + 2) \\
\end{pmatrix}, & if $-1 < t < L - 1$.
\end{dcases*}
\end{equation}
Assume that $E(2)$ is surjective and take $w \in \mathsf Z$. Let $x = (w, 0) \in \mathsf X$ and take $u \in \mathsf Y_2$ such that $E(2) u = x$. Hence, by \eqref{E2u}, one has that $u(t + 2) = 0$ for $t \in (-1, 0)$, i.e., $u(t) = 0$ for $t \in (1, 2)$. Thus $u(t + 2 - L) = 0$ for $L - 1 < t < 0$, and one obtains from \eqref{E2u} that
\[
\left\{
\begin{aligned}
\beta u(t + 1) + \alpha u(t + 1 - L) & = w(t), & \qquad & \text{ if } L - 1 < t < 0, \\
\beta u(t + 1) + u(t + 2 - L) & = w(t), & & \text{ if } -1 < t < L - 1.
\end{aligned}
\right.
\]
This shows that $S u(\cdot + 1) = w$, and thus $S$ is surjective. Assume now that $S$ is surjective and take $x = (x_1, x_2) \in \mathsf X$. Let $\widetilde u \in \mathsf Z$ be such that
\begin{equation}
\label{DefVTilde}
S \widetilde u(t) = 
\begin{dcases*}
x_1(t) - x_2(t - L), & if $L - 1 < t < 0$, \\
x_1(t), & if $-1 < t < L - 1$,
\end{dcases*}
\end{equation}
and define $u \in \mathsf Y_2$ by $u(t) = \widetilde u(t - 1)$ if $0 < t < 1$ and $u(t) = x_2(t - 2)$ if $1 < t < 2$. Then, combining \eqref{DefS}, \eqref{E2u}, and \eqref{DefVTilde}, one obtains that $E(2) u = x$, which yields the surjectivity of $E(2)$. Hence the surjectivity of $E(2)$ is equivalent to that of $S$. The fact that the latter is equivalent to the existence of $c > 0$ such that $\norm{S^\ast x}_{\mathsf Z} \geq c \norm{x}_{\mathsf Z}$ for every $x \in \mathsf Z$ is a classical result in functional analysis (see, e.g., \cite[Theorem~4.13]{Rudin1991Functional}).
\end{prf}

\begin{remk}
\label{RemkGraphicSSAst}
As in Remark \ref{RemkGraphic}, one can provide a graphical representation for the operators $S$ and $S^\ast$. Notice first that, for $A_1$, $A_2$, $B$, and $(\Lambda_1, \Lambda_2)$ under the form \eqref{ReducGc}, the only line segments $\sigma_{\mathbf n}$ from Remark \ref{RemkGraphic} lying inside the domain $[0, 2) \times [-1, 0)$ and associated with non-zero coefficients are $\sigma_{(0, 0)}$, $\sigma_{(0, 1)}$, $\sigma_{(1, 0)}$, and $\sigma_{(1, 1)}$, which are associated respectively with the coefficients $\begin{pmatrix}0 \\ 1 \\\end{pmatrix}$, $\begin{pmatrix}1 \\ 0 \\\end{pmatrix}$, $\begin{pmatrix}\beta \\ 0 \\\end{pmatrix}$, and $\begin{pmatrix}\alpha \\ 0 \\\end{pmatrix}$.

\newlength{\OneLine}
\setlength{\OneLine}{\baselineskip}

\begin{figure}[ht]
\centering
\begin{tabular}{@{} c @{} c @{}}
\raisebox{-\height+\OneLine}{\begin{tikzpicture}[scale=3.75]

\draw[thick, arrows={-Stealth[scale=1.5]}] (-0.0625, 0) -- (2.12500, 0) node[right] {$\xi$};
\draw[thick, arrows={-Stealth[scale=1.5]}] (0, -1.06250) -- (0, 0.125) node[left] {$\zeta$};

\draw[thick, dashed, color=darkgreen] (2.00000, -1.06250) -- (2.00000, 0.0625) node[above, color=black] {$2$};
\draw[thick, dashed, color=darkgreen] (2.06250, -1.00000) -- (-0.0625, -1.00000) node[left, color=black] {$-1$};

% Droites
\draw (2.00000, 0) -- node[pos=0.19000, below, rotate=45]{\scriptsize $\begin{pmatrix}0 \\ 1\\\end{pmatrix}$} (1.00000, -1);
\draw (1.00000, 0) -- node[pos=0.19000, below, rotate=45]{\scriptsize $\begin{pmatrix}\beta \\ 0\\\end{pmatrix}$} (0, -1.00000);
\draw (1.61850, 0) -- node[pos=0.19000, below, rotate=45]{\scriptsize $\begin{pmatrix}1 \\ 0\\\end{pmatrix}$} (0.61850, -1);
\draw (0.61850, 0) -- node[pos=0.30720, below, rotate=45]{\scriptsize $\begin{pmatrix}\alpha \\ 0\\\end{pmatrix}$} (0, -0.61850);

% Droites verticales
\draw[color=darkgreen, dashed, thick] (1.00000, -1.0625) -- (1.00000, 0.0625) node[above, color=black] {$1$};
\draw[color=darkgreen, dashed, thick] (1.61850, -1.0625) -- (1.61850, 0.0625) node[above, color=black] {$2-L$};
\draw[color=darkgreen, dashed, thick] (0.61850, -1.0625) -- (0.61850, 0.0625) node[above, color=black] {$1-L$};

% Droites horizontales
\draw[color=darkgreen, dashed, thick] (2.06250, -0.38150) -- (-0.0625, -0.38150) node[left, color=black] {$-L$};
\draw[color=darkgreen, dashed, thick] (2.06250, -0.61850) -- (-0.0625, -0.61850) node[left, color=black] {$L - 1$};

\draw[arrows={|[scale=1.5]-|[scale=1.5]}] (0, -1.2) -- node[midway, below] {$\mathcal E_2$} (1, -1.2);
\draw[arrows={|[scale=1.5]-|[scale=1.5]}] (1, -1.2) -- node[midway, below] {$\mathcal E_1$} (2, -1.2);

\end{tikzpicture}} & \raisebox{-\height + \OneLine}{\begin{tikzpicture}[scale=3.75]

\draw[thick, arrows={-Stealth[scale=1.5]}] (-1.0625, 0) -- (0.125, 0) node[right] {$\xi$};
\draw[thick, arrows={-Stealth[scale=1.5]}] (0, -1.0625) -- (0, 0.125) node[left] {$\zeta$};

\draw[thick, dashed, color=darkgreen] (-1, -1.0625) -- (-1, 0.0625) node[above, color=black] {$-1$};
\draw[thick, dashed, color=darkgreen] (-1.0625, -1) -- (0.0625, -1) node[right, color=black] {$-1$};

% Droites
\draw (0, 0) -- node[pos=0.5, below, rotate=45]{\footnotesize $\beta$} (-1, -1);
\draw (0, -0.61850) -- node[pos=0.5, below, rotate=45]{\footnotesize $1$} (-0.38150, -1);
\draw (-0.38150, 0) -- node[pos=0.5, below, rotate=45]{\footnotesize $\alpha$} (-1, -0.61850);

% Droites verticales
\draw[color=darkgreen, dashed, thick] (-0.38150, -1.0625) -- (-0.38150, 0.0625) node[above, color=black] {$-L$};

% Droites horizontales
\draw[color=darkgreen, dashed, thick] (-1.0625, -0.61850) -- (0.0625, -0.61850) node[right, color=black] {$L - 1$};

\end{tikzpicture}} \tabularnewline
(a) & (b) \tabularnewline
\end{tabular}
\caption{Graphical representations of the operators (a) $E(2)$ and $E(2)^\ast$, and (b) $S$ and $S^\ast$.}
\label{FigET2-S}
\end{figure}
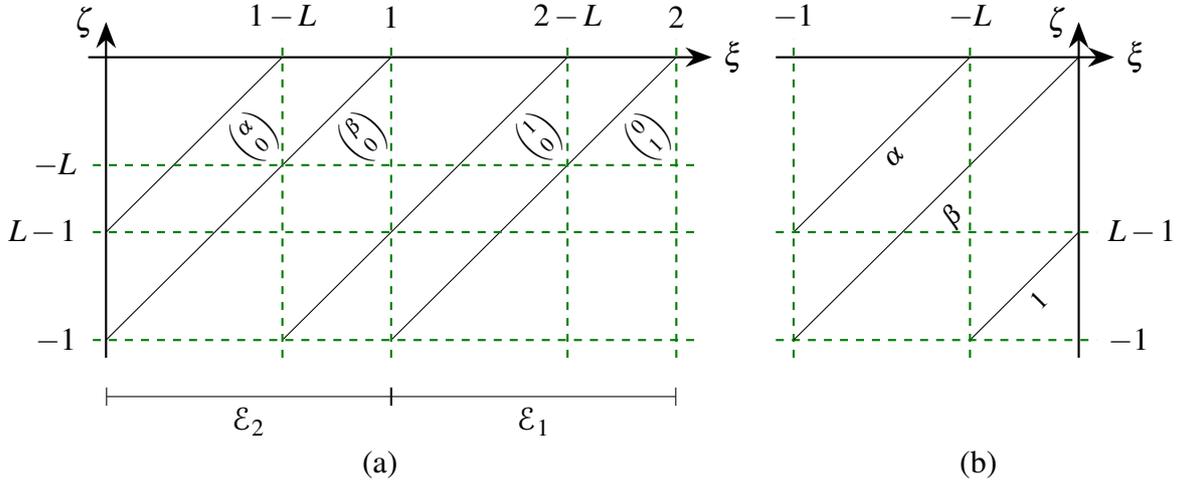

Figure \ref{FigET2-S}(a) provides the graphical representation for $E(2)$ and $E(2)^\ast$ given in Remark \ref{RemkGraphic}. One can decompose the domain $[0, 2) \times [-1, 0)$ in two parts, $\mathcal E_1 = [1, 2) \times [-1, 0)$ and $\mathcal E_2 = [0, 1) \times [-1, 0)$. The value of $E(2)^\ast x(t)$ for $t \in [0, 1)$, which corresponds to the region $\mathcal E_2$, only depends on $x_1$, and $S^\ast$ is defined as the operator that, to each $x_1$, associates the value of $E(2)^\ast x(t)$ for $t \in (0, 1)$, translated by $1$ in order to obtain as a result a function defined in $(-1, 0)$. Hence $S^\ast$ can be seen as the part of $E(2)^\ast$ corresponding to the region $\mathcal E_2$, which is represented in Figure \ref{FigET2-S}(b). It turns out that this part of $E(2)^\ast$ is enough to characterize its injectivity and the surjectivity of its adjoint, as shown in Lemma~\ref{LemmS}.

In the case of commensurable delays, i.e., $L = \frac{p}{q}$ with $p, q \in \mathbb N^\ast$ coprime and $p < q$, one can associate with $S^\ast$ a Toeplitz matrix $M = (m_{ij})_{i, j \in \llbracket 1, q\rrbracket} \in \mathcal M_q(\mathbb C)$, similar to the construction of $C$ and $\Epsilon$ from $E(T)$ performed in Remark \ref{RemkGraphCommens}, and defined by
\begin{equation}
\label{NewDefMij}
m_{ij} = 
\begin{dcases*}
\overline\beta, & if $j = i$, \\
\overline\alpha, & if $j = i - p$, \\
1, & if $j = i + q - p$, \\
0, & otherwise.
\end{dcases*}
\end{equation}
A graphical way to represent $M$ goes as follows. We decompose $(-1, 0)^2$ into squares $\mathsf S_{ij} = \left(-\frac{i}{q}, -\frac{i - 1}{q}\right) \times \left(-\frac{j}{q}, -\frac{j-1}{q}\right)$ for $i, j \in \llbracket 1, q\rrbracket$. Remark that the intersection between one of the line segments representing $S^\ast$ and the square $\mathsf S_{ij}$ is either empty, and in this case $m_{ij} = 0$, or equal to the diagonal of the square from its bottom left corner to its top right corner, in which case $m_{ij}$ is the conjugate of the coefficient corresponding to the intersecting line. Figure \ref{FigSRational} illustrates such a construction in the case $L = \frac{3}{7}$. The link between $M$ and $S^\ast$ is made more explicit in \eqref{ReductionToMatrix}.

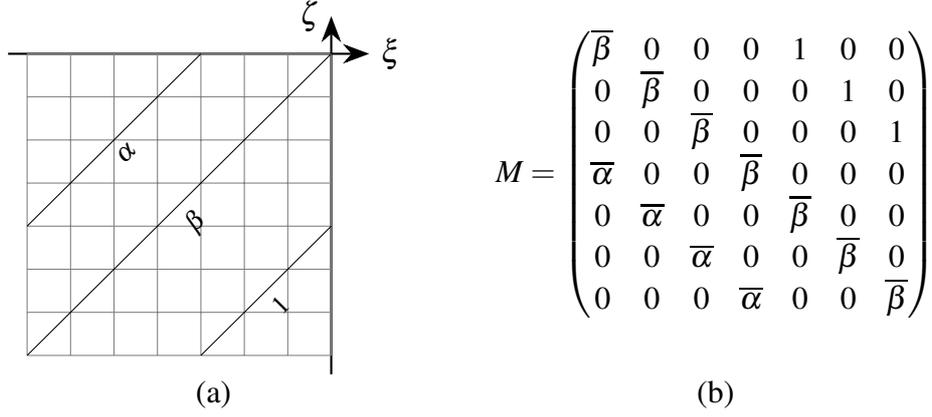
\begin{figure}[ht]
\centering
\begin{tabular}{@{} c @{\hspace*{1cm}} c @{}}
\raisebox{-0.5\height}{\begin{tikzpicture}[scale=4]

\draw[thick, arrows={-Stealth[scale=1.5]}] (-1.0625, 0) -- (0.125, 0) node[right] {$\xi$};
\draw[thick, arrows={-Stealth[scale=1.5]}] (0, -1.0625) -- (0, 0.125) node[left] {$\zeta$};

\draw[step=0.14286cm, color=gray] (0, 0) grid (-1, -1);% Droites
\draw (0, 0) -- node[pos=0.5, below, rotate=45]{\footnotesize $\beta$} (-1, -1);
\draw (0, -0.57143) -- node[pos=0.5, below, rotate=45]{\footnotesize $1$} (-0.42857, -1);
\draw (-0.42857, 0) -- node[pos=0.5, below, rotate=45]{\footnotesize $\alpha$} (-1, -0.57143);

\end{tikzpicture}} & $M = \begin{pmatrix}
\overline\beta & 0 & 0 & 0 & 1 & 0 & 0 \\
0 & \overline\beta & 0 & 0 & 0 & 1 & 0 \\
0 & 0 & \overline\beta & 0 & 0 & 0 & 1 \\
\overline\alpha & 0 & 0 & \overline\beta & 0 & 0 & 0 \\
0 & \overline\alpha & 0 & 0 & \overline\beta & 0 & 0 \\
0 & 0 & \overline\alpha & 0 & 0 & \overline\beta & 0 \\
0 & 0 & 0 & \overline\alpha & 0 & 0 & \overline\beta \\
\end{pmatrix}$ \tabularnewline
(a) & (b) \tabularnewline
\end{tabular}
\caption{Construction of the matrix $M$ from $S^\ast$ in the case $L = \frac{3}{7}$.}
\label{FigSRational}
\end{figure}
\end{remk}

\paragraph{Proof of Theorem~\ref{MainTheo22}\ref{MainTheoA1A2Contr}\ref{MainTheoApprox}}

\begin{prfsection}
Assume that $(A_1, B)$ and $(A_2, B)$ are controllable, in which case, according to Remark \ref{RemkReduction}, we can assume that $A_1$, $A_2$, $B$, and $(\Lambda_1, \Lambda_2)$ are under the form \eqref{ReducGc}. It has already been proved that approximate controllability does not hold for $T < 2$. Thanks to Lemma~\ref{LemmS}, one is left to show that the operator $S^\ast$ from \eqref{SAst} is injective if and only if $0 \notin \mathcal S$. We write in this proof $\alpha = \abs{\alpha} e^{i\theta}$ for some $\theta \in (-\pi, \pi]$.

Consider first the case $L \in (0, 1) \cap \mathbb Q$ and write $L = \frac{p}{q}$ for $p, q \in \mathbb N^\ast$ coprime. Define the operator $R \in \mathcal L\left(\mathsf Z, L^2\left(\left(-{1}/{q}, 0\right), \mathbb C^q\right)\right)$ by
\[
(R x(t))_n = x\left(t - \frac{n-1}{q}\right), \qquad -\frac{1}{q} < t < 0, \; n \in \llbracket 1, q\rrbracket.
\]
One immediately verifies from its definition that $R$ is a unitary transformation and that, for every $x \in L^2\left(\left(-{1}/{q}, 0\right), \mathbb C^q\right)$,
\begin{equation}
\label{ReductionToMatrix}
(R S^\ast R^{-1} x)(t) = M x(t),
\end{equation}
where $M$ is the matrix defined in \eqref{NewDefMij}. One has
\begin{equation}
\label{CRational}
\mathcal S = \left\{\beta + \abs{\alpha}^{1 - \frac{p}{q}} e^{i (\theta + 2 k \pi)\left(1 - \frac{p}{q}\right)} \midsuchthat k \in \llbracket 0, q - 1\rrbracket\right\}.
\end{equation}
Notice that $0 \in \mathcal S$ if and only if $\det M = 0$. Indeed, by Proposition~\ref{AppPropM}\ref{AppCharactPolyM} in the appendix, one has $\det M = 0$ if and only if $(-\beta)^q = \alpha^{q - p}$, i.e., if and only if $-\beta$ is a $q$-th root of $\alpha^{q - p}$, which means that $-\beta = \abs{\alpha}^{\frac{q - p}{q}} e^{i(\theta + 2 k \pi)\frac{q - p}{q}}$ for some $k \in \llbracket 0, q-1\rrbracket$, this being equivalent to $0 \in \mathcal S$. Since $R$ is a unitary transformation, one obtains in particular that the injectivity of $S^\ast$ is equivalent to that of $R S^\ast R^{-1}$, which, thanks to \eqref{ReductionToMatrix}, is equivalent to that of $M$. Since $M$ is injective if and only if $\det M \not = 0$, one concludes that $S^\ast$ is injective if and only if $0 \notin \mathcal S$, as required.

Assume now that $L \in (0, 1) \setminus \mathbb Q$. Let $x \in \mathsf Z$ be such that $S^\ast x = 0$, i.e.,
\[
x(t) = 
\begin{dcases*}
-\frac{1}{\overline\beta} x(t + L - 1), & if $-L < t < 0$, \\
-\frac{\overline\alpha}{\overline\beta} x(t + L), & if $-1 < t < -L$.
\end{dcases*}
\]
Let $\varphi: [-1, 0) \to [-1, 0)$ be the translation by $L$ modulo $1$ on the interval $[-1, 0)$, i.e., $\varphi(t) = t + L$ if $t \in [-1, -L)$ and $\varphi(t) = t + L - 1$ if $t \in [-L, 0)$. Since $L$ is irrational, $\varphi$ is ergodic with respect to the Lebesgue measure in $[-1, 0)$ (see, e.g., \cite[Chapter II, Theorem~3.2]{Mane1987Ergodic}). We have
\[
x(t) = -\frac{\overline\alpha \chi_{(-1, -L)}(t) + \chi_{(-L, 0)}(t)}{\overline\beta} x \circ \varphi(t) \qquad \text{ for } -1 < t < 0.
\]

Choose $\gamma \in \mathbb C$ such that $e^{\gamma (1 - L)} = -\overline\beta$. If $0 \in \mathcal S$, we next show that $\gamma$ can be chosen so that $e^{\gamma} = \overline\alpha$ and that such a choice is unique. Indeed, since $0 \in \mathcal S$, one has $\alpha \neq 0$, for otherwise $\beta = 0$, which contradicts the controllability of $(A_1, B)$. Hence the set of solutions with respect to $\gamma$ of the equation $e^{\gamma} = \overline\alpha$ is equal to $\{\log \abs{\alpha} - i (\theta + 2 m \pi) \mid m \in \mathbb Z\}$. The condition $0 \in \mathcal S$ means that there exists $k \in \mathbb Z$ such that $\beta + \abs{\alpha}^{1 - L} e^{i (\theta + 2 k \pi)(1 - L)} = 0$, and thus $\gamma = \log \abs{\alpha} - i (\theta + 2 k \pi)$ satisfies both equations. As regards uniqueness, consider $\gamma^\prime \in \mathbb C$ satisfying $e^{\gamma^\prime (1 - L)} = -\overline\beta$ and $e^{\gamma^\prime} = \overline\alpha$. Then there exists an integer $k^\prime \in \mathbb Z$ such that $\gamma^\prime = \log \abs{\alpha} - i (\theta + 2 k^\prime \pi)$ and $\beta + \abs{\alpha}^{1 - L} e^{i (\theta + 2 k \pi)(1 - L)} = \beta + \abs{\alpha}^{1 - L} e^{i (\theta + 2 k^\prime \pi)(1 - L)} = 0$. Hence $(k - k^\prime)(1 - L)$ is an integer, which implies $k = k^\prime$ since $L \notin \mathbb Q$.

Let $y \in \mathsf Z$ be defined by $y(t) = e^{\gamma t} x(t)$, i.e., $y$ is the function satisfying
\begin{equation}
\label{FixedPointV}
y(t) = \left(\overline\alpha e^{-\gamma} \chi_{(-1, -L)}(t) + \chi_{(-L, 0)}(t)\right) y \circ \varphi(t) \qquad \text{ for } -1 < t < 0.
\end{equation}
If $0 \in \mathcal S$, then $\overline\alpha e^{-\gamma} = 1$, and thus $y$ satisfies $y = y \circ \varphi$. Since $\varphi$ is ergodic with respect to the Lebesgue measure in $[-1, 0)$, the set of functions $y \in \mathsf Z$ satisfying $y = y \circ \varphi$ is the set of functions constant almost everywhere (see, e.g., \cite[Chapter II, Proposition~2.1]{Mane1987Ergodic}). Hence 
\begin{equation}
\label{EqKerS}
\Ker S^\ast = \{t \mapsto c e^{-\gamma t} \mid c \in \mathbb C\},
\end{equation}
where $\gamma = \log\abs{\alpha} - i (\theta + 2 k \pi)$ for some integer $k$ and $e^{\gamma(1 - L)} = - \overline\beta$. Since such a $\gamma \in \mathbb C$ (i.e., integer $k$) is unique, $\Ker S^\ast$ is of dimension $1$. In particular, $S^\ast$ is not injective, as required.

If $0 \notin \mathcal S$, notice that, from \eqref{FixedPointV},
\[
\norm{y}_{\mathsf Z}^2 = \abs{\overline\alpha e^{-\gamma}}^2 \int_{L-1}^{0} \abs{y(t)}^2 dt + \int_{-1}^{L - 1} \abs{y(t)}^2 dt,
\]
which shows that
\[
\left(1 - \abs{\overline\alpha e^{-\gamma}}^2\right) \int_{L-1}^{0} \abs{y(t)}^2 dt = 0.
\]
Let us prove that $y$ vanishes in the interval $(L - 1, 0)$. If $\abs{\overline\alpha e^{-\gamma}} \not = 1$, this follows immediately from the above equality. If $\abs{\overline\alpha e^{-\gamma}} = 1$, write $\overline\alpha e^{-\gamma} = e^{i \frac{2\pi\eta L}{1 - L}}$ for some $\eta \in \left[0, \frac{1-L}{L}\right)$. Notice that, for every $n \in \mathbb Z$, one has $e^{i \frac{2\pi(\eta - n) L}{1 - L}} \not = 1$; indeed, one has $\overline\alpha = e^{\gamma + i \frac{2\pi\eta L}{1 - L}}$ and hence the possible complex values of $\overline\alpha^{1-L}$ are
\begin{equation}
\label{ComplexAlphaOneMinusL}
e^{\gamma (1 - L) + i (2\pi\eta L + 2\pi k (1 - L))} = -\overline\beta e^{2 i \pi L (\eta - k)}, \qquad k \in \mathbb Z.
\end{equation}
If $e^{i \frac{2\pi(\eta - n) L}{1 - L}} = 1$ for some $n \in \mathbb Z$, then $\eta \equiv n \mod \frac{1-L}{L}$ and, since $\frac{1-L}{L} = \frac{1}{L} - 1$, we conclude that there exists $k \in \mathbb Z$ such that $\eta \equiv k \mod \frac{1}{L}$. Then $e^{2 i \pi L (\eta - k)} = 1$, which is not possible due to \eqref{ComplexAlphaOneMinusL} since we are in the case $0 \notin \mathcal S$. Hence, for every $n \in \mathbb Z$, one has $e^{i \frac{2\pi(\eta - n) L}{1 - L}} \not = 1$. The function $y$ satisfies
\[
y(t) = \left(e^{i \frac{2\pi\eta L}{1 - L}} \chi_{(-1, -L)}(t) + \chi_{(-L, 0)}(t)\right) y \circ \varphi(t) \qquad \text{ for } -1 < t < 0.
\]
Thus, for every $n \in \mathbb Z$,
\begin{align*}
\int_{-1}^0 y(t) e^{i \frac{2 \pi n}{1-L} t} dt & = e^{i \frac{2\pi\eta L}{1 - L}} \int_{-1}^{-L} y(t + L) e^{i \frac{2 \pi n}{1-L} t} dt + \int_{-L}^0 y(t + L - 1) e^{i \frac{2 \pi n}{1-L} t} dt \\
 & = e^{i \frac{2\pi(\eta - n) L}{1 - L}} \int_{L-1}^{0} y(t) e^{i \frac{2 \pi n}{1-L} t} dt + \int_{-1}^{L - 1} y(t) e^{i \frac{2 \pi n}{1-L} t} dt,
\end{align*}
which implies that
\[
\left(1 - e^{i \frac{2\pi(\eta - n) L}{1 - L}}\right) \int_{L-1}^{0} y(t) e^{i \frac{2 \pi n}{1-L} t} dt = 0, \qquad \forall n \in \mathbb Z.
\]
Since $e^{i \frac{2\pi(\eta - n) L}{1 - L}} \not = 1$ for every $n \in \mathbb Z$, we conclude that
\[
\int_{L-1}^{0} y(t) e^{i \frac{2 \pi n}{1-L} t} dt = 0, \qquad \forall n \in \mathbb Z,
\]
which shows that all the Fourier coefficients of $y|_{(L-1, 0)}$ vanish. Thus $y$ is zero in the interval $(L-1, 0)$.

Since $y$ vanishes in $(L - 1, 0)$, it follows from \eqref{FixedPointV} and an immediate inductive argument that $y$ is zero in $\varphi^{-k}(L-1, 0)$ for every $k \in \mathbb N$, which shows that $y = 0$ in $(-1, 0)$ since $\varphi$ is ergodic (see, e.g., \cite[Theorem~1.5]{Walters1982Introduction}). Hence $x = 0$ is the unique solution of $S^\ast x = 0$, proving that $S^\ast$ is injective, as required.
\end{prfsection}

\begin{remk}
\label{RemkControlDelayRational}
One can also obtain from the previous proof that, if $L = \frac{p}{q}$ for some $p, q \in \mathbb N^\ast$ coprime, then approximate and exact controllability in time $T \geq 2$ are equivalent for \eqref{Control22}. Indeed, notice that, when \eqref{Control22} is approximately controllable in time $T \geq 2$, then $0 \notin \mathcal S$, $M$ is invertible, and hence, by \eqref{ReductionToMatrix}, one has $\norm{R S^\ast R^{-1} x}_{L^2\left(\left(-{1}/{q}, 0\right), \mathbb C^q\right)} \geq \abs{M^{-1}}_{2}^{-1} \norm{x}_{L^2\left(\left(-{1}/{q}, 0\right), \mathbb C^q\right)}$ for every $x \in L^2\left(\left(-{1}/{q}, 0\right), \mathbb C^q\right)$, which shows that $\norm{S^\ast x}_{\mathsf Z} \geq \abs{M^{-1}}_{2}^{-1} \norm{x}_{\mathsf Z}$ for every $x \in \mathsf Z$, thus giving the exact controllability of \eqref{Control22} in time $T \geq 2$ thanks to Lemma~\ref{LemmS}. This agrees with the general result of Proposition~\ref{PropCommET} for commensurable delays. Moreover, one obtains from \eqref{CRational} that the set $\mathcal S$ is finite, which shows that $\overline{\mathcal S} = \mathcal S$ and hence conditions $0 \notin \mathcal S$ and $0 \notin \overline{\mathcal S}$ are equivalent. This proves Theorem~\ref{MainTheo22}\ref{MainTheoA1A2Contr}\ref{MainTheoExact} in the case where $\Lambda_1$ and $\Lambda_2$ are commensurable, i.e., $\frac{\Lambda_2}{\Lambda_1} \in \mathbb Q$.
\end{remk}

\begin{remk}
When $0 \in \mathcal S$ and $L \notin \mathbb Q$, this proof also shows that the kernel of $S^\ast$ is the vector space spanned by the function $x(t) = e^{\gamma t}$ with $\gamma \in \mathbb C$ chosen as in the proof of the theorem. Thanks to \eqref{ETAstUtPlusT}, this means that the kernel of $E(2)^\ast$ is the vector space spanned by the function
\[
x(t) = 
\begin{pmatrix}
e^{-\gamma t} \\
-e^{-\gamma(t + L)} \chi_{(-1, -L)}(t)
\end{pmatrix}.
\]
\end{remk}

\begin{remk}
When $0 \in \mathcal S$, $L \notin \mathbb Q$, and $\alpha, \beta \in \mathbb R$, one has that $\gamma \in \mathbb R$, obtaining thus a real-valued nonzero solution to $S^\ast x = 0$, and hence to $E(2)^\ast x = 0$. Indeed, notice first that one can only have $0 \in \mathcal S$ with $\alpha, \beta \in \mathbb R$ if $\alpha > 0$ (in which case $\beta < 0$), since $\alpha = 0$ implies $\beta = 0$, which is not possible, and, for $\alpha < 0$, the equality $\beta + \alpha^{1 - L} = 0$ for some complex value of $\alpha^{1 - L}$ implies that $-\beta = \alpha^{1-L} = \abs{\alpha}^{1-L} e^{i (\pi + 2 n \pi)(1 - L)}$ for some $n \in \mathbb Z$, but such an expression cannot be real for any $n \in \mathbb Z$ since $L \notin \mathbb Q$. Now, when $\alpha > 0$, then $\gamma = \log\alpha \in \mathbb R$.
\end{remk}

\paragraph{Proof of Theorem~\ref{MainTheo22}\ref{MainTheoA1A2Contr}\ref{MainTheoExact}}

\begin{prfsection}
Assume that $(A_1, B)$ and $(A_2, B)$ are controllable, in which case, according to Remark \ref{RemkReduction}, we can assume that $A_1$, $A_2$, $B$, and $(\Lambda_1, \Lambda_2)$ are under the form \eqref{ReducGc}. Since one has already proved that exact controllability does not hold for $T < 2$, it suffices to show that, for $T \geq 2$, the system is exactly controllable if and only if $0 \notin \overline{\mathcal S}$. Remark \ref{RemkControlDelayRational} has already shown the result when $L \in (0, 1) \cap \mathbb Q$, and thus one is left to prove only the case $L \in (0, 1) \setminus \mathbb Q$. Thanks to Lemma~\ref{LemmS}, one is left to show that $0 \notin \overline{\mathcal S}$ if and only if the operator $S$ defined in \eqref{DefS} is surjective or, equivalently, if there exists $c > 0$ such that $S^\ast$ satisfies $\norm{S^\ast x}_{\mathsf Z} \geq c \norm{x}_{\mathsf Z}$ for every $x \in \mathsf Z$. We write in this proof $\alpha = \abs{\alpha} e^{i\theta}$ for some $\theta \in (-\pi, \pi]$.

Take $L \in (0, 1) \setminus \mathbb Q$. Notice first that $0 \in \overline{\mathcal S}$ if and only if $\abs{\beta} = \abs{\alpha}^{1-L}$. Indeed, one has
\[
\mathcal S = \left\{\beta + \abs{\alpha}^{1 - L} e^{i (\theta + 2 k \pi)\left(1 - L\right)} \midsuchthat k \in \mathbb Z\right\},
\]
and, since $L$ is irrational, $\overline{\mathcal S}$ is the circle in $\mathbb C$ of center $\beta$ and radius $\abs{\alpha}^{1-L}$. 

Let us first treat the case $\alpha = 0$. Since $\beta \not = 0$ due to the controllability of $(A_1, B)$, one has $0 \notin \overline{\mathcal S}$ in this case. We will prove the exact controllability of \eqref{Control22} by showing the surjectivity of $S$. Take $x \in \mathsf Z$ and define $u \in \mathsf Z$ by
\[u(t) = \sum_{k=0}^{\floor{\frac{t}{L - 1}}} \frac{(-1)^k}{\beta^{k+1}} x(t + k(1 - L)).\]
Then, for $L - 1 < t < 0$, one has $S u(t) = \beta u(t) = x(t)$ and, for $-1 < t < L-1$, one has
\begin{align*}
S u(t) & = \beta u(t) + u(t - L + 1) \displaybreak[0] \\
 & = \sum_{k=0}^{\floor{\frac{t}{L - 1}}} \frac{(-1)^k}{\beta^{k}} x(t + k(1 - L)) + \sum_{k=0}^{\floor{\frac{t - L + 1}{L - 1}}} \frac{(-1)^k}{\beta^{k+1}} x(t - L + 1 + k(1 - L)) \displaybreak[0] \\
 & = \sum_{k=0}^{\floor{\frac{t}{L - 1}}} \frac{(-1)^k}{\beta^{k}} x(t + k(1 - L)) + \sum_{k=1}^{\floor{\frac{t}{L - 1}}} \frac{(-1)^{k-1}}{\beta^{k}} x(t + k(1 - L)) = x(t),
\end{align*}
which shows that $S u = x$ and thus $S$ is surjective.

Consider now the case $\alpha \not = 0$. Suppose that $0 \notin \overline{\mathcal S}$, which means that $\abs{\beta} \not = \abs{\alpha}^{1-L}$. Let $(p_n), (q_n)$ be two sequences of positive integers such that $p_n$ and $q_n$ are coprime for every $n \in \mathbb N$ and $\frac{p_n}{q_n} \to L$ as $n \to \infty$. Let $r_n = q_n - p_n$. Up to eliminating a finite number of terms in the sequence, we can assume that $\abs{\beta} \not = \abs{\alpha}^\frac{r_n}{q_n}$ for every $n \in \mathbb N$. Let $S_n \in \mathcal L(\mathsf Z)$ be the operator whose adjoint $S_n^\ast$ is given by
\begin{equation*}
S_n^\ast x(t) = 
\begin{dcases}
\overline\beta x(t) + x\left(t + \frac{p_n}{q_n} - 1\right) & \text{ if } -\frac{p_n}{q_n} < t < 0, \\
\overline\beta x(t) + \overline\alpha x\left(t + \frac{p_n}{q_n}\right) & \text{ if } -1 < t < -\frac{p_n}{q_n}.
\end{dcases}
\end{equation*}
One easily verifies (using, e.g., \cite[Theorem~9.5]{Rudin1987Real}) that, for every $x \in \mathsf Z$, one has $S_n^\ast x \to S^\ast x$ as $n \to \infty$. Since $\abs{\beta} \not = \abs{\alpha}^{\frac{r_n}{q_n}}$ for every $n \in \mathbb N$, we obtain, from Remark \ref{RemkControlDelayRational}, that
\[\norm{S_n^\ast x}_{\mathsf Z} \geq \abs{M_n^{-1}}_{2}^{-1} \norm{x}_{\mathsf Z},\]
where $M_n$ is given by \eqref{NewDefMij} with $p$ and $q$ replaced respectively by $p_n$ and $q_n$. Hence, by Proposition~\ref{AppPropM}\ref{AppInverseM},
\[\norm{S_n^\ast x}_{\mathsf Z} \geq \frac{\abs{\abs{\beta} - \abs{\alpha}^{\frac{r_n}{q_n}}}}{\max\left(\abs{\alpha}, \abs{\alpha}^{-1}\right)} \norm{x}_{\mathsf Z},\]
and, letting $n \to \infty$,
\[\norm{S^\ast x}_{\mathsf Z} \geq \frac{\abs{\abs{\beta} - \abs{\alpha}^{1-L}}}{\max\left(\abs{\alpha}, \abs{\alpha}^{-1}\right)} \norm{x}_{\mathsf Z},\]
which proves the surjectivity of $S$.

For $a, b \in \mathbb C$ with $a \neq 0$, let $S_{a, b}^\ast \in \mathcal L(\mathsf Z)$ be defined by
\[
S_{a, b}^\ast x(t) = 
\begin{dcases}
\overline b x(t) + x(t + L - 1) & \text{ if } -L < t < 0, \\
\overline b x(t) + \overline a x(t + L) & \text{ if } -1 < t < -L.
\end{dcases}
\]
In particular, for every $\lambda \in \mathbb C$, one has $S_{a, b}^\ast - \lambda = S_{a, b - \overline\lambda}^\ast$. Let $\sigma_{\mathrm p}(S_{a, b}^\ast)$ denote the set of eigenvalues of $S_{a, b}^\ast$. Thus $\lambda \in \sigma_{\mathrm p}(S_{a, b}^\ast)$ if and only if $0 \in \sigma_{\mathrm p}(S_{a, b - \overline\lambda}^\ast)$, which, by the proof of Theorem~\ref{MainTheo22}\ref{MainTheoA1A2Contr}\ref{MainTheoApprox}, is the case if and only if $\overline b - \lambda + \overline a^{1-L} = 0$ for some complex value of $\overline a^{1-L}$. Hence $\sigma_{\mathrm p}(S^\ast)$ is the set of all possible values of $\overline\beta + \overline\alpha^{1-L}$.

Suppose now that $0 \in \overline{\mathcal S}$, i.e., that $\abs{\beta} = \abs{\alpha}^{1-L}$. Since $L$ is irrational, we conclude that $0 \in \overline{\sigma_{\mathrm p}(S^\ast)}$. Hence there exists a sequence $(\lambda_n)_{n \in \mathbb N}$ in $\sigma_{\mathrm p}(S^\ast)$ such that $\lambda_n \to 0$ as $n \to \infty$. For $n \in \mathbb N$, let $x_n$ be an eigenfunction of $S^\ast$ associated with the eigenvalue $\lambda_n$ and with $\norm{x_n}_{\mathsf Z} = 1$. Hence $S^\ast x_n = \lambda_n x_n \to 0$ as $n \to +\infty$, which shows that there does not exist $c > 0$ such that $\norm{S^\ast x}_{\mathsf Z} \geq c \norm{x}_{\mathsf Z}$ for every $x \in \mathsf Z$, and thus $S$ is not surjective.
\end{prfsection}

\begin{remk}
\label{RemkEigenpairsS}
It follows from the above proof and \eqref{EqKerS} that, for $L \in (0, 1) \setminus \mathbb Q$ and $\alpha \neq 0$, one has a complete description of the eigenvalues and eigenfunctions of $S^\ast$. The set of eigenvalues of $S^\ast$ is $\{\overline\beta + \abs{\alpha}^{1 - L} e^{- i (\theta + 2 k \pi)(1 - L)} \mid k \in \mathbb Z\}$, where $\theta \in \mathbb R$ is an argument of $\alpha$. In addition, every eigenvalue $\lambda$ is simple, with corresponding eigenfunction $x(t) = e^{- \gamma t}$, where $\gamma$ is the unique solution of
\[
\left\{
\begin{aligned}
e^{\gamma} & = \overline\alpha, \\
e^{\gamma(1 - L)} & = - (\overline\beta - \lambda).
\end{aligned}
\right.
\]
\end{remk}

\section{Controllability to constants}
\label{SecContrConst}

The notions of controllability provided in Definition \ref{DefiContrExactApprox} require the possibility of steering the state $x_t$ of \eqref{MainSyst} towards (or arbitrarily close to) an arbitrary state of the infinite-dimensional space $\mathsf X$. We show in this section the equivalence between such controllability notions and notions which are in appearance much weaker, since they involve only target states belonging to a finite-dimensional space.

\begin{defi}
\label{DefiControlToCst}
Let $T \in (0, +\infty)$.
Define $\mathsf K$ by
\begin{equation}
\label{DefMathsfK}
\mathsf K = \left\{x \in \mathsf X \midsuchthat x: (-\Lambda_{\max}, 0) \to \mathbb C^d \text{ is a constant function}\right\}.
\end{equation}

\begin{enumerate}
\item We say that \eqref{MainSyst} is \emph{approximately controllable to constants in time $T$} if $\overline{\range E(T)} \supset \mathsf K$, i.e., for every  $y \in \mathsf K$ and $\varepsilon > 0$, there exists $u \in \mathsf Y_T$ such that the solution $x$ of \eqref{MainSyst} with initial condition $0$ and control $u$ satisfies $\norm{x_T - y}_{\mathsf X} < \varepsilon$.

\item We say that \eqref{MainSyst} is \emph{exactly controllable to constants in time $T$} if $\range E(T) \supset \mathsf K$, i.e., for every $y \in \mathsf K$, there exists $u \in \mathsf Y_T$ such that the solution $x$ of \eqref{MainSyst} with initial condition $0$ and control $u$ satisfies $x_T = y$.
\end{enumerate}
\end{defi}

 As we have proved in Lemma~\ref{LemmControlDelayInvariant} for approximate and exact controllability, approximate and exact controllability to constants are also preserved under linear change of coordinates, linear feedback, and changes of the time scale.

\begin{lemm}
\label{LemmInvariantConstants}
Let $T > 0$, $\lambda > 0$, $K_j \in \mathcal M_{m, d}(\mathbb C)$ for $j \in \llbracket 1, N\rrbracket$, $P \in \mathrm{GL}_d(\mathbb C)$, and consider System \eqref{MainSystChange}. Then
\begin{enumerate}
\item\label{ApproxContrlConstIff} \eqref{MainSyst} is approximately controllable to constants in time $T$ if and only if \eqref{MainSystChange} is approximately controllable to constants in time $\frac{T}{\lambda}$;

\item\label{ExactContrlConstIff} \eqref{MainSyst} is exactly controllable to constants in time $T$ if and only if \eqref{MainSystChange} is exactly controllable to constants in time $\frac{T}{\lambda}$.
\end{enumerate}
\end{lemm}

The following analogue of Proposition~\ref{PropControl} will also be of use in the sequel.

\begin{prop}
\label{PropInegObservConst}
Let $T \in (0, +\infty)$. System \eqref{MainSyst} is exactly controllable to constants in time $T$ if and only if there exists $c > 0$ such that, for every $x \in \mathsf X$,
\[
\norm{E(T)^\ast x}_{\mathsf Y_T}^2 \geq c \abs{\int_{-\Lambda_{\max}}^0 x(s) ds}_2^2.
\]
\end{prop}

\begin{prf}
Let $\kappa \in \mathcal L(\mathbb C^d, \mathsf X)$ be the canonical injection of $\mathbb C^d$ into $\mathsf X$, i.e., for $v \in \mathbb C^d$, $\kappa v$ is the constant function identically equal to $v$. Then clearly $\range \kappa = \mathsf K$, where $\mathsf K$ is defined by \eqref{DefMathsfK}, and thus \eqref{MainSyst} is exactly controllable to constants in time $T$ if and only if $\range \kappa \subset \range E(T)$. By classical results on functional analysis (see, e.g., \cite[Lemma~2.48]{Coron2007Control}), the latter condition is equivalent to the existence of $c > 0$ such that, for every $x \in \mathsf X$,
\[
\norm{E(T)^\ast x}_{\mathsf Y_T}^2 \geq c \abs{\kappa^\ast x}_2^2.
\]
This concludes the proof, since $\kappa^\ast x = \int_{-\Lambda_{\max}}^0 x(s) ds$, as one can verify by a straightforward computation.
\end{prf}

\subsection{Approximate controllability to constants}
\label{SecApproxToConst}

The main result of this section, Theorem~\ref{TheoApproxConstant}, states that approximate controllability and approximate controllability to constants are equivalent. Its proof relies on the following lemma, inspired by \cite[Theorem~5.1]{Gohberg1989Inversion}, which provides a link between the operator $E(T)$ and some suitable integration operators.

\begin{lemm}
\label{LemmPETETQF}
Let $T \in (0, +\infty)$. Define the bounded linear operators $P \in \mathcal L(\mathsf X)$, $Q \in \mathcal L(\mathsf Y_T)$, and $F \in \mathcal L(\mathsf Y_T, \mathbb C^d)$ by
\begin{align*}
(P x)(t) & = \int_{-\Lambda_{\max}}^t x(s) ds, & & x \in \mathsf X,\; t \in (-\Lambda_{\max}, 0), \displaybreak[0] \\
(Q u)(t) & = \int_0^t u(s) ds, & & u \in \mathsf Y_T,\; t \in (0, T), \displaybreak[0] \\
F u & = \sum_{\substack{\mathbf n \in \mathbb N^N \\ \Lambda \cdot \mathbf n \leq T - \Lambda_{\max}}} \Xi_{\mathbf n} B \int_0^{T - \Lambda_{\max} - \Lambda \cdot \mathbf n} u(s) ds, & & u \in \mathsf Y_T.
\end{align*}
Then $\norm{P}_{\mathcal L(X)} \leq \frac{\sqrt{2}\Lambda_{\max}}{2}$ and
\begin{equation}
\label{EqPETETQF}
P E(T) = E(T) Q - F.
\end{equation}
\end{lemm}

\begin{prf}
For $x \in \mathsf X$, one has
\begin{align*}
\norm{P x}_{\mathsf X}^2 & = \int_{-\Lambda_{\max}}^0 \abs{\int_{-\Lambda_{\max}}^t x(s) ds}_2^2 dt \leq \int_{-\Lambda_{\max}}^0 \left(\int_{-\Lambda_{\max}}^t \abs{x(s)}_2^2 ds\right) (t + \Lambda_{\max}) dt \displaybreak[0] \\
& \leq \norm{x}_{\mathsf X}^2 \int_{-\Lambda_{\max}}^0 (t + \Lambda_{\max}) dt = \frac{\Lambda_{\max}^2}{2} \norm{x}_{\mathsf X}^2,
\end{align*}
and thus $\norm{P}_{\mathcal L(X)} \leq \frac{\sqrt{2}\Lambda_{\max}}{2}$.

Let $u \in \mathsf Y_T$ and extend $u$ by zero in the interval $(-\infty, 0)$. Then, for almost every $t \in (-\Lambda_{\max}, 0)$,
\begin{align*}
(P E(T) u)(t) & = \int_{-\Lambda_{\max}}^t \sum_{\mathbf n \in \mathbb N^N} \Xi_{\mathbf n} B u(T + s - \Lambda \cdot \mathbf n) ds \displaybreak[0] \\
& = \sum_{\mathbf n \in \mathbb N^N} \Xi_{\mathbf n} B \int_{T - \Lambda_{\max} - \Lambda \cdot \mathbf n}^{T + t - \Lambda \cdot \mathbf n} u(s) ds \displaybreak[0] \\
& = \sum_{\mathbf n \in \mathbb N^N} \Xi_{\mathbf n} B \left[\int_0^{T + t - \Lambda \cdot \mathbf n} u(s) ds - \int_0^{T - \Lambda_{\max} - \Lambda \cdot \mathbf n} u(s) ds\right] \displaybreak[0] \\
& = (E(T) Q u)(t) - F u,
\end{align*}
where we use that the above infinite sums have only finitely many non-zero terms.
\end{prf}

As a consequence of Lemma~\ref{LemmPETETQF}, one obtains that approximate controllability to constants implies approximate controllability to polynomials.

\begin{lemm}
\label{LemmApproxToPoly}
Let $T \in (0, +\infty)$ and assume that \eqref{MainSyst} is approximately controllable to constants in time $T$. Then, for every polynomial $p: (-\Lambda_{\max}, 0) \to \mathbb C^d$ and $\varepsilon > 0$, there exists $u \in \mathsf Y_T$ such that $\norm{E(T) u - p}_{\mathsf X} < \varepsilon$.
\end{lemm}

\begin{prf}
Let $P$, $Q$, and $F$ be as in Lemma~\ref{LemmPETETQF}. We prove the result by induction on the degree of the polynomial. The result is true for polynomials of degree at most $0$ since this is precisely the definition of approximate controllability to constants.

Assume that $r \in \mathbb N$ is such that, for every polynomial $p: (-\Lambda_{\max}, 0) \to \mathbb C^d$ of degree at most $r$ and $\varepsilon > 0$, there exists $u \in \mathsf Y_T$ such that $\norm{E(T) u - p}_{\mathsf X} < \varepsilon$. Let $q: (-\Lambda_{\max}, 0) \to \mathbb C^d$ be a polynomial of degree $r + 1$ and take $\varepsilon > 0$. Let $a_0, \dotsc, a_{r+1} \in \mathbb C^d$ be such that
\[q(t) = \sum_{n=0}^{r + 1} a_n (t + \Lambda_{\max})^n, \qquad \forall t \in (-\Lambda_{\max}, 0).\]

Since $t \mapsto a_{r+1} (r+1) (t + \Lambda_{\max})^r$ is a polynomial of degree $r$, thanks to the induction hypothesis, there exists $u_0 \in \mathsf Y_T$ such that
\[
\norm{E(T) u_0 - a_{r+1} (r+1) (\cdot + \Lambda_{\max})^r}_{\mathsf X} < \frac{\sqrt{2}\varepsilon}{3 \Lambda_{\max}}.
\]
Hence, since $P \left[a_{r+1} (r+1) (\cdot + \Lambda_{\max})^r\right] = a_{r+1} (\cdot + \Lambda_{\max})^{r+1}$ and $\norm{P}_{\mathcal L(\mathsf X)} < \frac{\sqrt{2} \Lambda_{\max}}{2}$, one obtains that
\[
\norm{P E(T) u_0 - a_{r+1} (\cdot + \Lambda_{\max})^{r+1}}_{\mathsf X} < \frac{\varepsilon}{3},
\]
which yields, thanks to \eqref{EqPETETQF},
\begin{equation}
\label{EstimU0}
\norm{E(T) Q u_0 - F u_0 - a_{r+1} (\cdot + \Lambda_{\max})^{r+1}}_{\mathsf X} < \frac{\varepsilon}{3}.
\end{equation}

Since $F u_0$ is a constant vector, there exists $u_1 \in \mathsf Y_T$ such that
\begin{equation}
\label{EstimU1}
\norm{E(T) u_1 + F u_0}_{\mathsf X} < \frac{\varepsilon}{3}.
\end{equation}
Since $t \mapsto \sum_{n=0}^r a_n (t + \Lambda_{\max})^r$ is a polynomial of degree at most $r$, there exists $u_2 \in \mathsf Y_T$ such that
\begin{equation}
\label{EstimU2}
\norm{E(T) u_2 - \sum_{n=0}^r a_n (\cdot + \Lambda_{\max})^r}_{\mathsf X} < \frac{\varepsilon}{3}.
\end{equation}
Let $u = Q u_0 + u_1 + u_2 \in \mathsf Y_T$. Combining \eqref{EstimU0}, \eqref{EstimU1}, and \eqref{EstimU2}, one finally obtains that
\[
\norm{E(T) u - q}_{\mathsf X} < \varepsilon,
\]
which concludes the inductive argument.
\end{prf}

Since the set of all $\mathbb C^d$-valued polynomials defined on $(-\Lambda_{\max}, 0)$ is dense in $\mathsf X$, one obtains as an immediate consequence the main result of this section.

\begin{theo}
\label{TheoApproxConstant}
Let $T \in (0, +\infty)$. Then \eqref{MainSyst} is approximately controllable in time $T$ if and only if it is approximately controllable to constants in time $T$.
\end{theo}

\subsection{Exact controllability to constants}

In this section, we are interested in the relation between exact controllability and exact controllability to constants. The technique used in Section \ref{SecApproxToConst} to prove Theorem~\ref{TheoApproxConstant} does not seem well adapted to treat such a question, since, even though one can easily adapt Lemma~\ref{LemmApproxToPoly} to prove that exact controllability to constants implies exact controllability to polynomials, this is not sufficient to decide whether exact controllability holds.

We rely instead in the characterization of exact controllability to constants from Proposition~\ref{PropInegObservConst}. We are only able to treat the case of two-dimensional systems with two delays and a scalar control, since, in that case, the tools from Section \ref{SecTwoTwo}, and in particular the spectral decomposition of the operator $S^\ast$ from \eqref{SAst}, are available. The general case remains an open problem.

Let us then consider System \eqref{Control22}, i.e.,
\begin{equation*}
\tag{\ref{Control22}}
x(t) = A_1 x(t - \Lambda_1) + A_2 x(t - \Lambda_2) + B u(t),
\end{equation*}
where $x(t) \in \mathbb C^2$, $u(t) \in \mathbb C$, $A_1, A_2 \in \mathcal M_2(\mathbb C)$, and $B \in \mathbb C^2$, and we still assume, without loss of generality, that $\Lambda_1 > \Lambda_2$. We start by proving that the analogue of Lemma~\ref{LemmTime2} for exact controllability to constants also holds.

\begin{lemm}
\label{LemmContConstTIff2}
Let $A_1, A_2 \in \mathcal M_2(\mathbb C)$, $B \in \mathcal M_{2, 1}(\mathbb C)$, and $(\Lambda_1, \Lambda_2) \in (0, +\infty)^2$ with $\Lambda_1 > \Lambda_2$, and assume that $(A_1, B)$ and $(A_2, B)$ are controllable. Then \eqref{Control22} is exactly controllable to constants in some time $T \geq 2 \Lambda_1$ if and only if it is exactly controllable to constants in time $T = 2 \Lambda_1$.
\end{lemm}

\begin{prf}
Thanks to Lemma~\ref{LemmInvariantConstants}, one can proceed as in Remark \ref{RemkReduction} and assume with no loss of generality that $A_1, A_2, B$, and $(\Lambda_1, \Lambda_2)$ are given by \eqref{ReducGc}, in which case $E(T)^\ast$ is given by \eqref{ETAstUtPlusT}.

Notice that, for every $T \geq 2$, there exists $C_T > 0$ such that, for every $x \in \mathsf X$,
\begin{equation}
\label{EquivE2ET}
\norm{E(2)^\ast x}_{\mathsf Y_T}^2 \leq \norm{E(T)^\ast x}_{\mathsf Y_T}^2 \leq C_T \norm{E(2)^\ast x}_{\mathsf Y_2}^2.
\end{equation}
Indeed, the first inequality is trivial since, by \eqref{ETAstUtPlusT}, $(E(2)^\ast x)(t + 2) = (E(T)^\ast x)(t + T)$ for every $t \in (-2, 0)$, and the second inequality has been shown in the proof of Lemma~\ref{LemmTime2}. The conclusion of the lemma now follows from Proposition~\ref{PropInegObservConst}.
\end{prf}

In order to prove an analogue of Lemma~\ref{LemmS} for exact controllability to constants, we first introduce the space $\mathsf K_{\mathrm r}(L)$ defined for $L \in (0, 1)$ by
\[
\mathsf K_{\mathrm r}(L) = \left\{x \in \mathsf Z \midsuchthat x \text{ is constant on the intervals } (-1, L-1) \text{ and } (L-1, 0)\right\}.
\]

\begin{lemm}
\label{LemmIneqObservS}
Let $A_1, A_2 \in \mathcal M_2(\mathbb C)$, $B \in \mathcal M_{2, 1}(\mathbb C)$, $(\Lambda_1, \Lambda_2) \in (0, +\infty)^2$ with $\Lambda_1 > \Lambda_2$, and $L = \Lambda_2 / \Lambda_1$. Assume that $(A_1, B)$ and $(A_2, B)$ are controllable. Let $S \in \mathcal L(\mathsf Z)$ be the operator defined in \eqref{DefS}. Then \eqref{Control22} is exactly controllable to constants in some time $T \geq 2 \Lambda_1$ if and only if $\range S \supset \mathsf K_{\mathrm r}(L)$, or, equivalently, if there exists $c > 0$ such that, for every $x \in \mathsf Z$,
\begin{equation}
\label{IneqObservS}
\norm{S^\ast x}_{\mathsf Z}^2 \geq c \left(\abs{\int_{-1}^{L-1} x(t) dt}^2 + \abs{\int_{L-1}^{0} x(t) dt}^2\right).
\end{equation}
\end{lemm}

\begin{prf}
As in the proof of Lemma~\ref{LemmContConstTIff2}, we assume, with no loss of generality, that $A_1, A_2, B$, and $(\Lambda_1, \Lambda_2)$ are given by \eqref{ReducGc}. By Lemma~\ref{LemmContConstTIff2},  \eqref{Control22} is exactly controllable to constants in some time $T \geq 2$ if and only if $\range E(2) \supset \mathsf K$.

Assume that \eqref{Control22} is exactly controllable to constants in some time $T \geq 2$ and take $y \in \mathsf K_{\mathrm r}(L)$. Let $(a, b) \in \mathbb C^2$ be such that
\[
y(t) = \begin{dcases*}
a, & if $L - 1 < t < 0$, \\
b, & if $-1 < t < L - 1$.
\end{dcases*}
\]
Consider the function $z \in \mathsf K$ given by $z(t) = (b, b - a)$ for every $t \in (-1, 0)$. Since $\range E(2) \supset \mathsf K$, there exists $u \in \mathsf Y_2$ such that $E(2) u = z$, i.e.,
\[
\left\{
\begin{aligned}
\begin{pmatrix}
\beta u(t + 1) + \alpha u(t + 1 - L) + u(t + 2 - L) \\
u(t + 2) \\
\end{pmatrix} & = \begin{pmatrix}
b \\
b - a \\
\end{pmatrix}, & & \text{ if } L - 1 < t < 0, \\
\begin{pmatrix}
\beta u(t + 1) + u(t + 2 - L) \\
u(t + 2) \\
\end{pmatrix} & = \begin{pmatrix}
b \\
b - a \\
\end{pmatrix}, & & \text{ if } -1 < t < L - 1,
\end{aligned}
\right.
\]
where we use the explicit expression of $E(2)$ from \eqref{E2u}. Hence
\[
\left\{
\begin{aligned}
u(t) & = b - a, & & \text{ if } 1 < t < 2, \\
\beta u(t + 1) + \alpha u(t + 1 - L) + u(t + 2 - L) & = b, & & \text{ if } L - 1 < t < 0, \\
\beta u(t + 1) + u(t + 2 - L) & = b, & & \text{ if } -1 < t < L - 1,
\end{aligned}
\right.
\]
and, since $t + 2 - L \in (1, 2)$ for $L - 1 < t < 0$, one obtains that
\begin{equation}
\label{SystU}
\left\{
\begin{aligned}
\beta u(t + 1) + \alpha u(t + 1 - L) & = a, & & \text{ if } L - 1 < t < 0, \\
\beta u(t + 1) + u(t + 2 - L) & = b, & & \text{ if } -1 < t < L - 1.
\end{aligned}
\right.
\end{equation}
Let $x \in \mathsf Z$ be defined by $x(t) = u(t + 1)$ for $-1 < t < 0$. Then \eqref{SystU} means precisely that $S x = y$, and thus $\mathsf K_{\mathrm r}(L) \subset \range S$.

Assume now that $\mathsf K_{\mathrm r}(L) \subset \range S$ and take $x \in \mathsf K$. Let $(a, b) \in \mathbb C^2$ be such that $x(t) = (a, b)$ for $t \in (-1, 0)$. Let $y \in \mathsf Z$ be given for $t \in (-1, 0)$ by
\[
y(t) =  \begin{dcases*}
a - b, & if $L - 1 < t < 0$, \\
a, & if $-1 < t < L - 1$.
\end{dcases*}
\]
Hence $y \in \mathsf K_{\mathrm r}(L)$, and thus there exists $z \in \mathsf Z$ such that $S z = y$, i.e., for $t \in (-1, 0)$,
\[
\left\{
\begin{aligned}
\beta z(t) + \alpha z(t - L) & = a - b, & & \text{ if } L - 1 < t < 0, \\
\beta z(t) + z(t + 1 - L) & = a, & & \text{ if } -1 < t < L - 1.
\end{aligned}
\right.
\]
Let $u \in \mathsf Y_2$ be defined by
\[
u(t) = 
\begin{dcases*}
z(t - 1), & if $0 < t < 1$, \\
b, & if $1 < t < 2$.
\end{dcases*}
\]
Then, for $t \in (-1, 0)$,
\[
\left\{
\begin{aligned}
u(t + 2) & = b, & & \text{ if } -1 < t < 0, \\
\beta u(t + 1) + \alpha u(t + 1 - L) + u(t + 2 - L) & = a, & & \text{ if } L - 1 < t < 0, \\
\beta u(t + 1) + u(t + 2 - L) & = a, & & \text{ if } -1 < t < L - 1,
\end{aligned}
\right.
\]
and, using the explicit expression \eqref{E2u} of $E(2)$, one obtains that $E(2) u = x$. Then $\mathsf K \subset \range E(2)$, and thus \eqref{Control22} is exactly controllable to constants in time $T \geq 2$.

Finally, let $\kappa_{\mathrm r} \in \mathcal L(\mathbb C^2, \mathsf Z)$ be the bounded linear operator defined for $(a, b) \in \mathbb C^2$ by
\[
(\kappa_{\mathrm r} (a, b)) (t) = 
\begin{dcases*}
a, & if $L - 1 < t < 0$, \\
b, & if $-1 < t < L - 1$.
\end{dcases*}
\]
Then $\range \kappa_{\mathrm r} = \mathsf K_{\mathrm r}(L)$, which means that \eqref{Control22} is exactly controllable to constants in time $T \geq 2$ if and only if $\range \kappa_{\mathrm r} \subset \range S$. By classical results on functional analysis (see, e.g., \cite[Lemma~2.48]{Coron2007Control}), the latter condition is equivalent to the existence of $c > 0$ such that, for every $x \in \mathsf X$,
\begin{equation}
\label{IneqObservSKappaR}
\norm{S^\ast x}_{\mathsf Z}^2 \geq c \abs{\kappa_{\mathrm r}^\ast x}_2^2.
\end{equation}
By a straightforward computation, one obtains that
\[
\kappa_{\mathrm r}^\ast x = 
\left(
\begin{aligned}
\displaystyle\int_{L-1}^0 x(t) dt \\
\displaystyle\int_{-1}^{L-1} x(t) dt
\end{aligned}
\right),
\]
and thus \eqref{IneqObservSKappaR} is the same as \eqref{IneqObservS}.
\end{prf}

We can now state the main result of this section.

\begin{theo}
\label{TheoExactConstant}
Let $T \in (0, +\infty)$. Then \eqref{Control22} is exactly controllable in time $T$ if and only if it is exactly controllable to constants in time $T$.
\end{theo}

\begin{prf}
Notice that exact controllability in time $T$ implies exact controllability to constants in time $T$, which in turn implies approximate controllability to constants in time $T$, the latter being equivalent, thanks to Theorem~\ref{TheoApproxConstant}, to approximate controllability in time $T$. Hence, equivalence between exact controllability to constants in time $T$ and exact controllability in time $T$ is true in particular when approximate and exact controllability in time $T$ are equivalent. Thanks to Theorem~\ref{MainTheo22}, this is the case if at least one of the following conditions holds.
\begin{itemize}
\item $T < 2 \Lambda_1$;
\item $(A_1, B)$ or $(A_2, B)$ is not controllable;
\item $(A_1, B)$ and $(A_2, B)$ are controllable and $0 \notin \overline{\mathcal S} \setminus \mathcal S$, where $\mathcal S \subset \mathbb C$ is as in the statement of Theorem~\ref{MainTheo22}\ref{MainTheoA1A2Contr}.
\end{itemize}
Hence Theorem~\ref{TheoExactConstant} is proved in such situations, and one is left to consider the case where $T \geq 2 \Lambda_1$, $(A_1, B)$ and $(A_2, B)$ are controllable, and $0 \in \overline{\mathcal S} \setminus \mathcal S$.

Assume that $T \geq 2 \Lambda_1$, $(A_1, B)$ and $(A_2, B)$ are controllable, and $0 \in \overline{\mathcal S} \setminus \mathcal S$. Notice that, due to the definition of $\mathcal S$, one has $\Lambda_2 / \Lambda_1 \notin \mathbb Q$ in this case. Thanks to Theorem~\ref{MainTheo22}\ref{MainTheoA1A2Contr}, \eqref{Control22} is not exactly controllable in time $T$, and thus the proposition is proved if one shows that \eqref{Control22} is not exactly controllable to constants in time $T$.

As in Lemmas \ref{LemmContConstTIff2} and \ref{LemmIneqObservS}, we assume, with no loss of generality, that $A_1, A_2, B$, and $(\Lambda_1, \Lambda_2)$ are given by \eqref{ReducGc}, with $L = \Lambda_2 / \Lambda_1$. Let $\alpha, \beta \in \mathbb C$ be as in the statement of Theorem~\ref{MainTheo22}\ref{MainTheoA1A2Contr}, $\theta_\alpha, \theta_\beta \in \mathbb R$ be the arguments of $\alpha$ and $\beta$, respectively, and $S$ be the operator defined in \eqref{DefS}. Notice that, since $0 \in \overline{\mathcal S}\setminus \mathcal S$, one has $\alpha \neq 0$. Define the operators $M_\alpha \in \mathcal L(\mathsf Z)$ and $\widehat S \in \mathcal L(\mathsf Z)$ for $x \in \mathsf Z$ by
\[M_\alpha x(t) = e^{-(\log\abs{\alpha} - i \theta_\alpha) t} x(t) \quad \text{ and } \quad \widehat S = M_{\alpha}^{-1} S^\ast M_\alpha.\]
According to Remark \ref{RemkEigenpairsS}, the eigenvalues of $\widehat S$ are $\lambda_k = \overline\beta + \abs{\alpha}^{1 - L} e^{- i (\theta_\alpha + 2 k \pi)(1 - L)}$ for $k \in \mathbb Z$, with corresponding eigenfunctions $e_k$ given, for $t \in (-1, 0)$, by $e_k(t) = e^{2 i k \pi t}$.

Notice that
\[
\min(1, \abs{\alpha}) \leq \norm{M_\alpha}_{\mathcal L(\mathsf Z)} \leq \max(1, \abs{\alpha}),
\]
and thus, for every $x \in \mathsf Z$,
\[
\min\left(\abs{\alpha}^2, \abs{\alpha}^{-2}\right) \norm{S^\ast x}_{\mathsf Z}^2 \leq \norm{\widehat S x}_{\mathsf Z}^2 \leq \max\left(\abs{\alpha}^2, \abs{\alpha}^{-2}\right) \norm{S^\ast x}_{\mathsf Z}^2.
\]
Hence, thanks to Lemma~\ref{LemmIneqObservS}, \eqref{Control22} is exactly controllable to constants in time $T$ if and only if there exists $c > 0$ such that, for every $x \in \mathsf Z$,
\begin{equation}
\label{IneqObservSHat}
\norm{\widehat S x}_{\mathsf Z}^2 \geq c \left(\abs{\int_{-1}^{L-1} x(t) dt}^2 + \abs{\int_{L-1}^{0} x(t) dt}^2\right).
\end{equation}

Assume, to obtain a contradiction, that \eqref{Control22} is exactly controllable to constants in time $T$, and let $c > 0$ be such that \eqref{IneqObservSHat} holds for every $x \in \mathsf Z$. Notice that 
\begin{equation}\label{eq:mn}
\tfrac{\pi + \theta_\beta - \theta_\alpha(1 - L)}{2 \pi} \text{ is not of the form } m (1 - L) + n \text{ for } m, n \in \mathbb Z.
\end{equation}
Indeed, if it were the case, one would have $\pi + \theta_\beta \equiv (\theta_\alpha - 2 \pi m) (1 - L) \mod 2 \pi$; since $0 \in \overline{\mathcal S} \setminus \mathcal S$, one has $\abs{\beta} = \abs{\alpha}^{1 - L}$, and thus $-\overline\beta = \abs{\alpha}^{1 - L} e^{- i (\theta_\alpha - 2 \pi m) (1 - L)}$, which contradicts the fact that $0 \notin \mathcal S$. Hence, using the Inhomogeneous Diophantine Approximation Theorem (see, e.g., \cite[Chapter III, Theorem~II A]{Cassels1957Introduction}), there exist two sequences $(p_n)_{n \in \mathbb N}$ and $(q_n)_{n \in \mathbb N}$ in $\mathbb Z$ with $\abs{q_n} \to \infty$ as $n \to \infty$ such that, for every $n \in \mathbb N$, one has $q_n \neq 0$ and
\begin{equation}
\label{DiophantineEq}
\abs{2 \pi q_n (1 - L) - (\pi + \theta_\beta - \theta_\alpha (1 - L)) - 2 \pi p_n} < \frac{\pi}{2 \abs{q_n}}.
\end{equation}
Recalling that $\abs{\beta} = \abs{\alpha}^{1 - L}$, one obtains that, for every $n \in \mathbb N$, the eigenvalue $\lambda_{q_n}$ of $\widehat S$ satisfies
\begin{equation}
\label{LambdaN}
\begin{aligned}
\lambda_{q_n} & = \overline\beta + \abs{\alpha}^{1 - L} e^{- i (\theta_\alpha + 2 \pi q_n)(1 - L)} \\
& = \overline\beta\left[1 + e^{-i(\theta_\alpha + 2 \pi q_n)(1 - L)} e^{i \theta_\beta}\right] \\
& = \overline\beta\left[1 - e^{i (2 \pi p_n + \pi + \theta_\beta - \theta_\alpha(1 - L) - 2 \pi q_n (1 - L))}\right].
\end{aligned}
\end{equation}
Notice that, if $z \in \mathbb C$ is such that $\abs{z} \leq 1$, then $\abs{1 - e^z} \leq 2 \abs{z}$. By \eqref{DiophantineEq}, one has $\lvert 2 \pi p_n + \pi + \theta_\beta - \theta_\alpha(1 - L) - 2 \pi q_n (1 - L) \rvert < \frac{\pi}{2 \abs{q_n}} \leq 1$ for $n$ large enough, and thus, for every $n$ large enough,
\begin{equation}
\label{EstimLambdaN}
\abs{\lambda_{q_n}} \leq \abs{\beta} \frac{\pi}{\abs{q_n}}.
\end{equation}
In particular, one has $\lambda_{q_n} \to 0$ as $n \to \infty$, and, by \eqref{LambdaN}, this also proves that $e^{- 2 i \pi q_n (1 - L)} \to e^{-i(\pi + \theta_\beta - \theta_\alpha(1 - L))}$ as $n \to \infty$. Notice that $e^{-i(\pi + \theta_\beta - \theta_\alpha(1 - L))} \neq 1$, as it follows from \eqref{eq:mn}. Hence there exists $C > 0$ such that, for every $n$ large enough,
\begin{equation}
\label{LowerBoundExp2PiQn1MinusLMinus1}
\abs{e^{- 2 i \pi q_n (1 - L)} - 1} \geq \frac{1}{C}.
\end{equation}
Fix $n_0 \in \mathbb N$ such that \eqref{EstimLambdaN} and \eqref{LowerBoundExp2PiQn1MinusLMinus1} hold for every $n \geq n_0$.

For $n \geq n_0$, define $x_n \in \mathsf Z$ by
\[
x_n = \frac{1}{n - n_0 + 1} \sum_{j = n_0}^n \frac{2 i \pi q_j}{e^{- 2 i \pi q_j (1 - L)} - 1} e_{q_j}.
\]
Then, since $\{e_k \mid k \in \mathbb Z\}$ is an orthonormal basis of $\mathsf Z$ made of the eigenfunctions of $\widehat S$, one has
\begin{equation}
\label{NormSXn}
\begin{aligned}
\norm{\widehat S x_n}_{\mathsf Z}^2 & = \sum_{j = n_0}^n \abs{\frac{2 i \pi q_j}{(e^{- 2 i \pi q_j (1 - L)} - 1)(n - n_0 + 1)}}^2 \abs{\lambda_{q_j}}^2 \\
& \leq \frac{4 \pi^2 C^2}{(n - n_0 + 1)^2} \sum_{j = n_0}^n \abs{q_j}^2 \abs{\lambda_{q_j}}^2 \\
& \leq \frac{4 \abs{\beta}^2 \pi^4 C^2}{n - n_0 + 1}.
\end{aligned}
\end{equation}
On the other hand, one computes
\begin{equation}
\label{IntXn}
\int_{-1}^{L - 1} x_n(t) dt = \frac{1}{n - n_0 + 1} \sum_{j = n_0}^n \frac{2 i \pi q_j}{e^{- 2 i \pi q_j (1 - L)} - 1} \frac{e^{- 2 i \pi q_j (1 - L)} - 1}{2 i \pi q_j} = 1.
\end{equation}
Hence, inserting \eqref{NormSXn} and \eqref{IntXn} into \eqref{IneqObservSHat}, one obtains that, for every $n \geq n_0$,
\[
\frac{4 \abs{\beta}^2 \pi^4 C^2}{n - n_0 + 1} \geq c,
\]
which implies, by taking the limit as $n \to \infty$, that $c \leq 0$, contradicting the fact that $c > 0$. This contradiction proves that \eqref{Control22} is not exactly controllable to constants in time $T$, as required.
\end{prf}

\section{Conclusion and open problems}

This paper has provided new results on the approximate and exact controllability of \eqref{MainSyst} in the function space $L^2((-\Lambda_{\max}, 0), \mathbb C^d)$. The case of commensurable delays has been completely characterized in Section \ref{SecCommensurable}, using both the classical augmented state space technique in Proposition \ref{PropCommensurable} and the explicit expression of the end-point operator $E(T)$ in Proposition \ref{PropCommET}, with a comparison between such techniques provided in Theorem \ref{TheoCommens}. In particular, approximate and exact controllability are equivalent in this context and can be characterized by the Kalman criterion from Proposition \ref{PropCommensurable}\ref{CtrlCommensurableKalman}.

A complete characterization of approximate and exact controllability has been provided in the first non-trivial case of \eqref{MainSyst} where incommensurable delays appear, namely the case $N = d = 2$ and $m = 1$. This complete characterization, provided in Theorem \ref{MainTheo22}, has been proved using several tools, the first one being a reduction to normal forms carried out in Section \ref{SecReduction}. The easy cases from Theorem \ref{MainTheo22}\ref{MainTheoA1BNonContr} and \ref{MainTheoA2BNonContr}, in which approximate and exact controllability are equivalent, were then studied using the expression of the end-point operator $E(T)$, with explicit constructions of controls in the cases where controllability holds.

The interesting and more subtle case from Theorem \ref{MainTheo22}\ref{MainTheoA1A2Contr} has been tackled using different tools, including classical characterizations of approximate and exact controllability in terms of the dual notions of unique continuation property and observability inequality, the ergodicity of translations by $L$ modulo $1$ when $L$ is irrational, and rational approximation of the delays combined with a fine spectral analysis of a sequence of Toeplitz matrices whose sizes tend to infinity.

We have also considered the notions of approximate and exact controllability to constants in Section \ref{SecContrConst}, proving in Theorem \ref{TheoApproxConstant} that approximate controllability and approximate controllability to constants are equivalent. The main tool in the proof of this result is \eqref{EqPETETQF}, which essentially means that a (sort of) commutator between integration and the end-point operator $E(T)$ is given by the operator $F$, which takes values in constant states. Exact controllability to constants has been proved to be equivalent to exact controllability in Theorem \ref{TheoExactConstant} in the case $N = d = 2$ and $m = 1$, whose proof is built upon the spectral analysis of $S^\ast$ from Remark \ref{RemkEigenpairsS} and uses an inhomogeneous Diophantine approximation result to bound the absolute value of a subsequence of the eigenvalues of an operator related to $S^\ast$.

We next propose two open problems that we believe to be interesting and challenging.

\begin{enumerate}[leftmargin=0pt, itemindent=*, listparindent=\parindent]
\item \emph{Is it possible to provide approximate and exact controllability criteria for \eqref{MainSyst} similar to Theorem \ref{MainTheo22} in higher dimensions and with more delays and control inputs?}

The most interesting case seems to be the analogue of Theorem \ref{MainTheo22}\ref{MainTheoA1A2Contr}, in which approximate and exact controllability are not equivalent and can be characterized in terms of the position of $0$ with respect to some set $\mathcal S$ constructed from the parameters of the system. It is not clear how the assumptions of \ref{MainTheoA1A2Contr} should be generalized to more than two delays, and many subtleties might appear depending on the ranks of the controllability matrices $\mathcal C(A_j, B)$ for $j \in \llbracket 1, N\rrbracket$. An important starting point would be to consider the case where all pairs $(A_j, B)$, $j \in \llbracket 1, N\rrbracket$, are controllable.

If one tries to follow the ideas of the proof of Theorem \ref{MainTheo22}, a first difficulty comes from the reduction to normal forms from Section \ref{SecReduction}. Even though similar reductions are still possible in higher dimensions and with more delays, explicit computations of $\Xi_{\mathbf n}$ and $E(T)$ used in Section \ref{SecTwoTwo} are much more tricky to handle. In particular, it is not immediate what should be a suitable generalization for the operator $S$ defined in \eqref{DefS}.

Concerning the main tools used in Section \ref{SecProofCaseC}, we expect the translations by $L$ modulo $1$ used in the analysis of approximate controllability to be replaced by more general interval exchange maps, on which ergodicity results are available (see, e.g., \cite{Viana2006Ergodic}). However, it is not clear how to transform approximate controllability into an interval exchange problem similar to \eqref{FixedPointV} in the general case. As regards the spectral analysis of Toeplitz matrices of sizes tending to infinity, it seems that reasonable generalizations of the operator $S$ would yield matrices that are only Toeplitz by blocks, whose spectral analysis seems intractable. We then expect a general characterization of exact controllability to rely on different techniques.

\item \emph{Are exact controllability and exact controllability to constants equivalent in general?}

The proof of Theorem \ref{TheoExactConstant} relies on spectral properties of $S$, and so we expect any generalization of this result using similar techniques to face the same difficulties as the general characterization of exact controllability.

%Another possible technique is to try an approach based on Lemma \ref{LemmPETETQF}, which allows one to prove that, similarly to Lemma \ref{LemmApproxToPoly}, exact controllability to constants implies exact controllability to polynomials. Then, given a target state $x \in \mathsf X$, one can construct a sequence $(p_n)_{n \in \mathbb N}$ of polynomials converging to $x$ and a sequence $(u_n)_{n \in \mathbb Y_T}$ of controls such that $E(T) u_n = p_n$. One is thus confronted with the problem of providing estimates on $u_n$ guaranteeing its convergence, up to the extraction of a subsequence, to some control $u$, which would then satisfy $E(T) u = x$. However, the lack of a suitable explicit characterization of $u_n$ in terms of $p_n$ renders this approach intractable.

\end{enumerate}

\appendix

\section{Appendix}
\label{AppM}

\begin{prop}
\label{AppPropM}
Let $\alpha, \beta \in \mathbb C$ and $p, q \in \mathbb N^\ast$ with $p, q$ coprime and $p < q$. Define the matrix $M = (m_{ij})_{i, j \in \llbracket 1, q\rrbracket} \in \mathcal M_{q}(\mathbb C)$ by
\begin{equation}
\label{DefMij}
m_{ij} = 
\begin{dcases*}
\overline\beta, & if $j = i$, \\
\overline\alpha, & if $j = i - p$, \\
1, & if $j = i + q - p$, \\
0, & otherwise.
\end{dcases*}
\end{equation}
Then the following holds.

\begin{enumerate}
\item\label{AppCharactPolyM} The characteristic polynomial and the determinant of $M$ are given by $P(\lambda) = \left(\lambda - \overline\beta\right)^q - \overline\alpha^{q - p}$ and $\det M = \overline\beta^q - (-1)^q \overline\alpha^{q - p}$, respectively.

\item\label{AppEigenvectorsM} Assume that $\alpha \not = 0$ and write $\alpha = \abs{\alpha} e^{i\theta}$ for some $\theta \in (-\pi, \pi]$. The eigenvalues of the matrix $M$ are
\begin{equation}
\label{EigenvaluesM}
\lambda_j = \overline\beta + \abs{\alpha}^{\frac{q - p}{q}} e^{-i\frac{\theta(q - p)}{q}} e^{i \frac{2\pi j (q - p)}{q}}, \qquad j \in \llbracket 1, q\rrbracket.
\end{equation}
For $j \in \llbracket 1, q\rrbracket$, a right eigenvector $v_j \in \mathbb C^q \simeq \mathcal M_{q, 1}(\mathbb C)$ of $M$ associated with $\lambda_j$ is
\[
v_j = \left(\abs{\alpha}^{\frac{k}{q}} e^{-i \frac{\theta k}{q}} e^{i \frac{2\pi j k}{q}}\right)_{k=1}^q
\]
and a left eigenvector $w_j \in \mathcal M_{1, q}(\mathbb C)$ of $M$ associated with $\lambda_j$ is
\[
w_j = \frac{1}{q} \left(\abs{\alpha}^{-\frac{k}{q}} e^{i \frac{\theta k}{q}} e^{-i \frac{2\pi j k}{q}}\right)_{k=1}^q.
\]
Moreover, for every $j, k \in \llbracket 1, q\rrbracket$, we have $w_k v_j = \delta_{jk}$, where $\delta_{jk}$ denotes the Kronecker delta, i.e., $\delta_{jk} = 1$ if $j = k$ and $\delta_{jk} = 0$ otherwise.

\item\label{AppInverseM} If $\alpha \not = 0$ and $\abs{\beta} \not = \abs{\alpha}^{\frac{q - p}{q}}$, then $M$ is invertible and
\[\abs{M^{-1}}_{2} \leq \frac{\max\left(\abs{\alpha}, \abs{\alpha}^{-1}\right)}{\abs{\abs{\beta} - \abs{\alpha}^{\frac{q - p}{q}}}}.\]
\end{enumerate}
\end{prop}

\begin{prf}
We start by proving \ref{AppCharactPolyM}. Set $M_\lambda = \lambda \id_q - M$ and notice that $P(\lambda) = \det M_\lambda$. Let $\mathfrak S_q$ denote the group of permutations of $\llbracket 1, q\rrbracket$ and $\epsilon(\sigma)$ denote the signature of an element $\sigma \in \mathfrak S_q$. Leibniz formula for the determinant gives
\begin{equation}
\label{DetLeibniz}
P(\lambda) = \det M_\lambda = \sum_{\sigma \in \mathfrak S_q} \epsilon(\sigma) \prod_{i=1}^q m_{i \sigma(i)}^{(\lambda)}.
\end{equation}
Thanks to \eqref{DefMij}, the product $\prod_{i=1}^q m_{i \sigma(i)}^{(\lambda)}$ is nonzero only if $\sigma \in \mathfrak S_q$ satisfies, for every $i \in \llbracket 1, q\rrbracket$,
\begin{equation}
\label{PermutationNonzero}
\sigma(i) \in
\begin{dcases*}
\{i, i + q - p\}, & if $i \in \llbracket 1, p\rrbracket$, \\
\{i, i - p\}, & if $i \in \llbracket p+1, q\rrbracket$.
\end{dcases*}
\end{equation}
Let $\tau \in \mathfrak S_q$ be the translation by $-1$ modulo $q$, i.e., $\tau(i) = i - 1$ if $i \in \llbracket 2, q\rrbracket$ and $\tau(1) = q$. We have $\epsilon(\tau) = (-1)^{q-1}$, and thus $\epsilon(\tau^p) = (-1)^{(q-1)p}$. Since $p, q$ are coprime, one has $pq \equiv p + q + 1 \mod 2$ and thus $(q-1)p \equiv q+1 \mod 2$, which gives $\epsilon(\tau^p) = (-1)^{q+1}$. Notice, moreover, that \eqref{PermutationNonzero} can be written as $\sigma(i) \in \{i, \tau^p(i)\}$ for every $i \in \llbracket 1, q\rrbracket$.

One immediately verifies that the only permutations $\sigma \in \mathfrak S_q$ satisfying \eqref{PermutationNonzero} are $\id_{\mathfrak S_q}$ and $\tau^p$. Then, it follows from \eqref{DetLeibniz} that
\[
P(\lambda) = \prod_{i=1}^q m_{ii}^{(\lambda)} + (-1)^{q+1} \prod_{i=1}^q m_{i \tau^p(i)}^{(\lambda)} = \left(\lambda - \overline\beta\right)^q + (-1)^{q+1} (-1)^q \overline\alpha^{q - p} = \left(\lambda - \overline\beta\right)^q - \overline\alpha^{q - p}.
\]
Moreover, $\det M = (-1)^q \det(-M) = (-1)^q P(0) = \overline\beta^q - (-1)^q \overline\alpha^{q - p}$.

We now turn to the proof of \ref{AppEigenvectorsM}. Formula \eqref{EigenvaluesM} for the eigenvalues of $M$ follows immediately from the expression of the characteristic polynomial of $M$. Let $j \in \llbracket 1, q\rrbracket$. For $k \in \llbracket 1, p\rrbracket$,
\begin{align*}
(M v_j)_k & = \overline\beta \abs{\alpha}^{\frac{k}{q}} e^{-i \frac{\theta k}{q}} e^{i \frac{2\pi j k}{q}} + \abs{\alpha}^{\frac{k+q-p}{q}} e^{-i \frac{\theta (k+q-p)}{q}} e^{i \frac{2\pi j (k+q-p)}{q}} \\
 & = \abs{\alpha}^{\frac{k}{q}} e^{-i \frac{\theta k}{q}} e^{i \frac{2\pi j k}{q}} \left(\overline\beta + \abs{\alpha}^{\frac{q-p}{q}} e^{-i \frac{\theta (q-p)}{q}} e^{i \frac{2\pi j (q-p)}{q}}\right) = \lambda_j (v_j)_k,
\end{align*}
and, for $k \in \llbracket p+1, q\rrbracket$,
\begin{align*}
(M v_j)_k & = \overline\beta \abs{\alpha}^{\frac{k}{q}} e^{-i \frac{\theta k}{q}} e^{i \frac{2\pi j k}{q}} + \overline\alpha \abs{\alpha}^{\frac{k-p}{q}} e^{-i \frac{\theta (k-p)}{q}} e^{i \frac{2\pi j (k-p)}{q}} \\
 & = \abs{\alpha}^{\frac{k}{q}} e^{-i \frac{\theta k}{q}} e^{i \frac{2\pi j k}{q}} \left(\overline\beta + \abs{\alpha}^{\frac{q-p}{q}} e^{-i \frac{\theta (q-p)}{q}} e^{i \frac{2\pi j (q-p)}{q}}\right) = \lambda_j (v_j)_k,
\end{align*}
which shows that $M v_j = \lambda_j v_j$, and hence $v_j$ is a right eigenvector of $M$ associated with $\lambda_j$. Now, for $k \in \llbracket 1, q-p\rrbracket$,
\begin{align*}
(w_j M)_k & = \frac{1}{q} \overline\beta \abs{\alpha}^{-\frac{k}{q}} e^{i \frac{\theta k}{q}} e^{-i \frac{2\pi j k}{q}} + \frac{1}{q} \overline\alpha \abs{\alpha}^{-\frac{k+p}{q}} e^{i \frac{\theta (k+p)}{q}} e^{-i \frac{2\pi j (k+p)}{q}} \\
 & = \frac{1}{q} \abs{\alpha}^{-\frac{k}{q}} e^{i \frac{\theta k}{q}} e^{-i \frac{2\pi j k}{q}} \left(\overline\beta + \abs{\alpha}^{\frac{q-p}{q}} e^{-i \frac{\theta (q-p)}{q}} e^{i \frac{2\pi j (q-p)}{q}}\right) = \lambda_j (w_j)_k,
\end{align*}
and, for $k \in \llbracket q-p+1, q\rrbracket$,
\begin{align*}
(w_j M)_k & = \frac{1}{q} \overline\beta \abs{\alpha}^{-\frac{k}{q}} e^{i \frac{\theta k}{q}} e^{-i \frac{2\pi j k}{q}} + \frac{1}{q} \abs{\alpha}^{-\frac{k+p-q}{q}} e^{i \frac{\theta (k+p-q)}{q}} e^{-i \frac{2\pi j (k+p-q)}{q}} \\
 & = \frac{1}{q} \abs{\alpha}^{-\frac{k}{q}} e^{i \frac{\theta k}{q}} e^{-i \frac{2\pi j k}{q}} \left(\overline\beta + \abs{\alpha}^{\frac{q-p}{q}} e^{-i \frac{\theta (q-p)}{q}} e^{i \frac{2\pi j (q-p)}{q}}\right) = \lambda_j (w_j)_k,
\end{align*}
which shows that $w_j M = \lambda_j w_j$, and hence $w_j$ is a left eigenvector of $M$ associated with $\lambda_j$. For $j, k \in \llbracket 1, q\rrbracket$, one evaluates immediately $w_k v_j = \frac{1}{q} \sum_{\ell=1}^q e^{i \frac{2\pi (j - k) \ell}{q}} = \delta_{jk}$.

To prove \ref{AppInverseM}, we first consider the matrices $V, W, D \in \mathcal M_{q}(\mathbb C)$ defined by
\[
V = (V_{jk})_{j, k \in \llbracket 1, q\rrbracket}, \qquad W = (W_{jk})_{j, k \in \llbracket 1, q\rrbracket}, \qquad D = (D_{jk})_{j, k \in \llbracket 1, q\rrbracket},
\]
with, for $j, k \in \llbracket 1, q\rrbracket$
\[
V_{jk} = (v_k)_j, \qquad W_{jk} = (w_j)_k, \qquad D_{jk} = \lambda_j \delta_{jk}.
\]
It follows from \ref{AppEigenvectorsM} that
\[
M = V D W \qquad \text{ and } \qquad V = W^{-1}.
\]

For simplicity, we set $r = q - p$. By \ref{AppCharactPolyM}, $M$ is invertible if and only if $\overline\beta^q - (-1)^q \overline\alpha^r \not = 0$, which is the case if $\alpha \not = 0$ and $\abs{\beta} \not = \abs{\alpha}^{\frac{r}{q}}$. In this case, $M^{-1} = V D^{-1} W$ and thus, for $j, k \in \llbracket 1, q\rrbracket$,
\begin{equation}
\label{DevelopMMinus1}
\begin{aligned}
\left(M^{-1}\right)_{jk} & = \sum_{\ell=1}^q (v_\ell)_j \lambda_\ell^{-1} (w_\ell)_k = \frac{\abs{\alpha}^{\frac{j-k}{q}} e^{-i\theta\frac{j-k}{q}}}{q} \sum_{\ell=1}^q \lambda_\ell^{-1} e^{i \frac{2\pi \ell (j - k)}{q}} = \\
 & = \frac{\abs{\alpha}^{\frac{j-k}{q}} e^{-i\theta\frac{j-k}{q}}}{q} \sum_{\ell=1}^q \frac{e^{i \frac{2\pi \ell (j - k)}{q}}}{\overline\beta + \abs{\alpha}^{\frac{r}{q}} e^{-i\frac{\theta r}{q}} e^{i \frac{2\pi \ell r}{q}}} = \frac{\abs{\alpha}^{\frac{j-k}{q}} e^{-i\theta\frac{j-k}{q}}}{q \overline\beta} \sum_{\ell=1}^q \frac{e^{i \frac{2\pi \ell (j - k)}{q}}}{1 + \frac{\abs{\alpha}^{\frac{r}{q}} e^{-i\frac{\theta r}{q}}}{\overline\beta} e^{i \frac{2\pi \ell r}{q}}}.
\end{aligned}
\end{equation}

We claim that, for every $z \in \mathbb C$ such that $z^q \not = 1$, we have
\begin{equation}
\label{CalculSommeMeromorphe}
\sum_{\ell=1}^q \frac{e^{i \frac{2\pi \ell (j - k)}{q}}}{1 - z e^{i \frac{2 \pi \ell r}{q}}} = \frac{q z^{d_{j, k}}}{1 - z^q},
\end{equation}
where $d_{j, k}$ is the unique integer in $\llbracket 0, q-1\rrbracket$ such that $r d_{j, k} + j - k \equiv 0 \mod q$, which is well-defined since $q$ and $r$ are coprime.

To show that \eqref{CalculSommeMeromorphe} holds for every $z \in \mathbb C$ such that $z^q \not = 1$, it suffices to show that it holds for $z \in \mathbb C$ with $\abs{z} < 1$, since both left- and right-hand sides of \eqref{CalculSommeMeromorphe} are meromorphic functions with simple poles at the $q$ roots of $z^q = 1$. If $z \in \mathbb C$ is such that $\abs{z} < 1$, then
\[
\sum_{\ell=1}^q \frac{e^{i \frac{2\pi \ell (j - k)}{q}}}{1 - z e^{i \frac{2 \pi \ell r}{q}}} = \sum_{\ell=1}^q e^{i \frac{2\pi \ell (j - k)}{q}} \sum_{s=0}^\infty z^s e^{i \frac{2 \pi \ell r s}{q}} = \sum_{s=0}^\infty z^s \sum_{\ell=1}^q e^{i \frac{2\pi \ell (r s + j - k)}{q}} = q z^{d_{j, k}} \sum_{t = 0}^\infty z^{t q} = \frac{q z^{d_{j, k}}}{1 - z^q},
\]
where we use that $\sum_{\ell=1}^q e^{i \frac{2\pi \ell (r s + j - k)}{q}} = q$ if $r s + j - k \equiv 0 \mod q$ and is equal to zero otherwise, and that $\{s \in \mathbb N \suchthat r s + j - k \equiv 0 \mod q\} = \{d_{j, k} + t q \suchthat t \in \mathbb N\}$. Hence \eqref{CalculSommeMeromorphe} is proved.

Since $\abs{\beta} \not = \abs{\alpha}^{\frac{r}{q}}$ implies $\overline\beta^q \not = (-1)^q \overline\alpha^r$, we have $\left(-\frac{\abs{\alpha}^{\frac{r}{q}} e^{-i\frac{\theta r}{q}}}{\overline\beta}\right)^q \not = 1$. Hence, combining \eqref{DevelopMMinus1} and \eqref{CalculSommeMeromorphe}, we obtain that
\[
\left(M^{-1}\right)_{jk} = \frac{\abs{\alpha}^{\frac{j-k}{q}} e^{-i\theta\frac{j-k}{q}}}{q \overline\beta} \frac{q \left(-\frac{\abs{\alpha}^{\frac{r}{q}} e^{-i\frac{\theta r}{q}}}{\overline\beta}\right)^{d_{j, k}}}{1 - \left(-\frac{\abs{\alpha}^{\frac{r}{q}} e^{-i\frac{\theta r}{q}}}{\overline\beta}\right)^q} = (-1)^{d_{j, k}} \frac{\overline\alpha^{n_{j, k}}\overline\beta^{q-1-d_{j,k}}}{\overline\beta^q - (-1)^q \overline\alpha^r},
\]
where $n_{j, k} \in \mathbb Z$ is the unique integer satisfying $r d_{j, k} + j - k = n_{j, k} q$; moreover, since $d_{j, k} \in \llbracket 0, q-1\rrbracket$ and $j, k \in \llbracket 1, q\rrbracket$, we have $n_{j, k} \in \llbracket 0, r\rrbracket$.

Notice that, for $j, k \in \llbracket 1, q\rrbracket$, $\frac{r d_{j, k}}{q} = n_{j, k} + \frac{k-j}{q}$, and hence $n_{j, k} = \floor{\frac{r d_{j, k}}{q}} + \delta_{j > k}$, where $\delta_{j > k} = 1$ if $j > k$ and $\delta_{j > k} = 0$ otherwise. Thus, for $k \in \llbracket 1, q\rrbracket$,
\[
\sum_{j=1}^q \abs{\left(M^{-1}\right)_{jk}} = \frac{1}{\abs{\beta^q - (-1)^q \alpha^r}} \sum_{j=1}^q \abs{\alpha}^{\floor{\frac{r d_{j, k}}{q}} + \delta_{j > k}} \abs{\beta}^{q-1-d_{j, k}}.
\]
Since $d_{j, k}$ is defined as the unique integer in $\llbracket 0, q-1\rrbracket$ satisfying $r d_{j, k} + j - k \equiv 0 \mod q$ and $r, q$ are coprime, we obtain that, for fixed $k \in \llbracket 1, q\rrbracket$, the map $j \mapsto d_{j, k}$ is a bijection between $\llbracket 1, q\rrbracket$ and $\llbracket 0, q-1\rrbracket$. Hence, when $\abs{\alpha} \geq 1$,
\begin{align*}
\sum_{j=1}^q \abs{\left(M^{-1}\right)_{jk}} & \leq \frac{\abs{\alpha}}{\abs{\beta^q - (-1)^q \alpha^r}} \sum_{j=0}^{q-1} \abs{\alpha}^{\floor{\frac{r j}{q}}} \abs{\beta}^{q-1-j} \leq \frac{\abs{\alpha} \abs{\beta}^{q-1}}{\abs{\beta^q - (-1)^q \alpha^r}} \sum_{j=0}^{q-1} \abs{\alpha}^{\frac{r j}{q}} \abs{\beta}^{-j} \displaybreak[0] \\
 & = \frac{\abs{\alpha} \abs{\beta}^{q-1}}{\abs{\beta^q - (-1)^q \alpha^r}} \abs{\frac{1 - \abs{\alpha}^r \abs{\beta}^{-q}}{1 - \abs{\alpha}^{\frac{r}{q}} \abs{\beta}^{-1}}} = \frac{\abs{\alpha}}{\abs{\abs{\beta} - \abs{\alpha}^{\frac{r}{q}}}} \frac{\abs{\abs{\beta}^q - \abs{\alpha}^r}}{\abs{\beta^q - (-1)^q \alpha^r}} \displaybreak[0] \\
& \leq \frac{\abs{\alpha}}{\abs{\abs{\beta} - \abs{\alpha}^{\frac{r}{q}}}},
\end{align*}
and, similarly, when $0 < \abs{\alpha} < 1$,
\begin{align*}
\sum_{j=1}^q \abs{\left(M^{-1}\right)_{jk}} & \leq \frac{1}{\abs{\beta^q - (-1)^q \alpha^r}} \sum_{j=0}^{q-1} \abs{\alpha}^{\floor{\frac{r j}{q}}} \abs{\beta}^{q-1-j} \leq \frac{\abs{\alpha}^{-1} \abs{\beta}^{q-1}}{\abs{\beta^q - (-1)^q \alpha^r}} \sum_{j=0}^{q-1} \abs{\alpha}^{\frac{r j}{q}} \abs{\beta}^{-j} \\
 & \leq \frac{\abs{\alpha}^{-1}}{\abs{\abs{\beta} - \abs{\alpha}^{\frac{r}{q}}}},
\end{align*}
which shows that
\[
\abs{M^{-1}}_{1} = \max_{k \in \llbracket 1, q\rrbracket} \sum_{j=1}^q \abs{\left(M^{-1}\right)_{jk}} \leq \frac{\max\left(\abs{\alpha}, \abs{\alpha}^{-1}\right)}{\abs{\abs{\beta} - \abs{\alpha}^{\frac{r}{q}}}}.
\]
A similar argument also shows that
\[
\abs{M^{-1}}_{\infty} = \max_{j \in \llbracket 1, q\rrbracket} \sum_{k=1}^q \abs{\left(M^{-1}\right)_{jk}} \leq \frac{\max\left(\abs{\alpha}, \abs{\alpha}^{-1}\right)}{\abs{\abs{\beta} - \abs{\alpha}^{\frac{r}{q}}}},
\]
and the result follows since $\abs{M^{-1}}_{2} \leq \sqrt{\abs{M^{-1}}_{1} \abs{M^{-1}}_{\infty}}$.
\end{prf}

\hyphenation{ar-Xiv}
\bibliographystyle{abbrv}
\bibliography{Bib}

\begin{thebibliography}{10}

\bibitem{Cassels1957Introduction}
J.~W.~S. Cassels.
\newblock {\em An introduction to {D}iophantine approximation}.
\newblock Cambridge Tracts in Mathematics and Mathematical Physics, No. 45.
  Cambridge University Press, New York, 1957.

\bibitem{Chitour2016Stability}
Y.~Chitour, G.~Mazanti, and M.~Sigalotti.
\newblock Stability of non-autonomous difference equations with applications to
  transport and wave propagation on networks.
\newblock {\em Netw. Heterog. Media}, 11(4):563--601, 2016.

\bibitem{Chitour2017Persistently}
Y.~Chitour, G.~Mazanti, and M.~Sigalotti.
\newblock Persistently damped transport on a network of circles.
\newblock {\em Trans. Amer. Math. Soc.}, 369(6):3841--3881, 2017.

\bibitem{Chyung1970Controllability}
D.~H. Chyung.
\newblock On the controllability of linear systems with delay in control.
\newblock {\em IEEE Trans. Automatic Control}, 15(2):255--257, 1970.

\bibitem{Cooke1968Differential}
K.~L. Cooke and D.~W. Krumme.
\newblock Differential-difference equations and nonlinear initial-boundary
  value problems for linear hyperbolic partial differential equations.
\newblock {\em J. Math. Anal. Appl.}, 24:372--387, 1968.

\bibitem{Coron2007Control}
J.-M. Coron.
\newblock {\em Control and nonlinearity}, volume 136 of {\em Mathematical
  Surveys and Monographs}.
\newblock American Mathematical Society, Providence, RI, 2007.

\bibitem{Coron2008Dissipative}
J.-M. Coron, G.~Bastin, and B.~d'Andr{\'e}a Novel.
\newblock Dissipative boundary conditions for one-dimensional nonlinear
  hyperbolic systems.
\newblock {\em SIAM J. Control Optim.}, 47(3):1460--1498, 2008.

\bibitem{Coron2015Dissipative}
J.-M. Coron and H.-M. Nguyen.
\newblock Dissipative boundary conditions for nonlinear 1-{D} hyperbolic
  systems: sharp conditions through an approach via time-delay systems.
\newblock {\em SIAM J. Math. Anal.}, 47(3):2220--2240, 2015.

\bibitem{Cruz1970Stability}
M.~Cruz, A. and J.~K. Hale.
\newblock Stability of functional differential equations of neutral type.
\newblock {\em J. Differential Equations}, 7:334--355, 1970.

\bibitem{Datko1977Linear}
R.~Datko.
\newblock Linear autonomous neutral differential equations in a {B}anach space.
\newblock {\em J. Diff. Equations}, 25(2):258--274, 1977.

\bibitem{Avellar1980Zeros}
C.~E. de~Avellar and J.~K. Hale.
\newblock On the zeros of exponential polynomials.
\newblock {\em J. Math. Anal. Appl.}, 73(2):434--452, 1980.

\bibitem{Diblik2008Controllability}
J.~Dibl{\'{\i}}k, D.~Y. Khusainov, and M.~Rů{\v{z}}i{\v{c}}kov{\'a}.
\newblock Controllability of linear discrete systems with constant coefficients
  and pure delay.
\newblock {\em SIAM J. Control Optim.}, 47(3):1140--1149, 2008.

\bibitem{Franklin1997Digital}
G.~F. Franklin, J.~D. Powell, and M.~L. Workman.
\newblock {\em Digital Control of Dynamic Systems}.
\newblock Addison-Wesley, 3 edition, 1997.

\bibitem{Fridman2010Bounds}
E.~Fridman, S.~Mondi{\'e}, and B.~Saldivar.
\newblock Bounds on the response of a drilling pipe model.
\newblock {\em IMA J. Math. Control Inform.}, 27(4):513--526, 2010.

\bibitem{Gohberg1989Inversion}
I.~Gohberg and T.~Shalom.
\newblock On inversion of square matrices partitioned into nonsquare blocks.
\newblock {\em Integral Equations Operator Theory}, 12(4):539--566, 1989.

\bibitem{Hale1985Stability}
J.~K. Hale, E.~F. Infante, and F.~S.~P. Tsen.
\newblock Stability in linear delay equations.
\newblock {\em J. Math. Anal. Appl.}, 105(2):533--555, 1985.

\bibitem{Hale1993Introduction}
J.~K. Hale and S.~M. Verduyn~Lunel.
\newblock {\em Introduction to functional-differential equations}, volume~99 of
  {\em Applied Mathematical Sciences}.
\newblock Springer-Verlag, New York, 1993.

\bibitem{Hale2002Strong}
J.~K. Hale and S.~M. Verduyn~Lunel.
\newblock Strong stabilization of neutral functional differential equations.
\newblock {\em IMA J. Math. Control Inform.}, 19(1-2):5--23, 2002.
\newblock Special issue on analysis and design of delay and propagation
  systems.

\bibitem{Henry1974Linear}
D.~Henry.
\newblock Linear autonomous neutral functional differential equations.
\newblock {\em J. Differential Equations}, 15:106--128, 1974.

\bibitem{Kloss2012Flow}
B.~Kl{\"o}ss.
\newblock The flow approach for waves in networks.
\newblock {\em Oper. Matrices}, 6(1):107--128, 2012.

\bibitem{Lions1986Controlabilite}
J.-L. Lions.
\newblock Contr\^olabilit\'e exacte des syst\`emes distribu\'es.
\newblock {\em C. R. Acad. Sci. Paris S\'er. I Math.}, 302(13):471--475, 1986.

\bibitem{Lions1988Exact}
J.-L. Lions.
\newblock Exact controllability, stabilization and perturbations for
  distributed systems.
\newblock {\em SIAM Rev.}, 30(1):1--68, 1988.

\bibitem{Mane1987Ergodic}
R.~Ma{\~n}{\'e}.
\newblock {\em Ergodic theory and differentiable dynamics}, volume~8 of {\em
  Ergebnisse der Mathematik und ihrer Grenzgebiete (3) [Results in Mathematics
  and Related Areas (3)]}.
\newblock Springer-Verlag, Berlin, 1987.
\newblock Translated from the Portuguese by Silvio Levy.

\bibitem{Mazanti2017Relative}
G.~Mazanti.
\newblock Relative controllability of linear difference equations.
\newblock {\em SIAM J. Control Optim.}, 55(5):3132--3153, 2017.

\bibitem{Melvin1974Stability}
W.~R. Melvin.
\newblock Stability properties of functional difference equations.
\newblock {\em J. Math. Anal. Appl.}, 48:749--763, 1974.

\bibitem{Michiels2009Strong}
W.~Michiels, T.~Vyhl{\'{\i}}dal, P.~Z{\'{\i}}tek, H.~Nijmeijer, and D.~Henrion.
\newblock Strong stability of neutral equations with an arbitrary delay
  dependency structure.
\newblock {\em SIAM J. Control Optim.}, 48(2):763--786, 2009.

\bibitem{Ngoc2015Exponential}
P.~H.~A. Ngoc and N.~D. Huy.
\newblock Exponential stability of linear delay difference equations with
  continuous time.
\newblock {\em Vietnam J. Math.}, 43(2):195--205, 2015.

\bibitem{OConnor1983Stabilization}
D.~A. O'Connor and T.~J. Tarn.
\newblock On stabilization by state feedback for neutral
  differential-difference equations.
\newblock {\em IEEE Trans. Automat. Control}, 28(5):615--618, 1983.

\bibitem{OConnor1983Function}
D.~A. O'Connor and T.~J. Tarn.
\newblock On the function space controllability of linear neutral systems.
\newblock {\em SIAM J. Control Optim.}, 21(2):306--329, 1983.

\bibitem{Pandolfi1976Stabilization}
L.~Pandolfi.
\newblock Stabilization of neutral functional differential equations.
\newblock {\em J. Optimization Theory Appl.}, 20(2):191--204, 1976.

\bibitem{Pospisil2015Relative}
M.~Posp{\'\i}{\v{s}}il, J.~Dibl{\'\i}k, and M.~Fe{\v{c}}kan.
\newblock On relative controllability of delayed difference equations with
  multiple control functions.
\newblock In {\em Proceedings of the International Conference on Numerical
  Analysis and Applied Mathematics 2014 (ICNAAM-2014)}, volume 1648, page
  130001. AIP Publishing, 2015.

\bibitem{Rudin1987Real}
W.~Rudin.
\newblock {\em Real and complex analysis}.
\newblock McGraw-Hill Book Co., New York, third edition, 1987.

\bibitem{Rudin1991Functional}
W.~Rudin.
\newblock {\em Functional analysis}.
\newblock International Series in Pure and Applied Mathematics. McGraw-Hill,
  Inc., New York, second edition, 1991.

\bibitem{Salamon1984Control}
D.~Salamon.
\newblock {\em Control and observation of neutral systems}, volume~91 of {\em
  Research Notes in Mathematics}.
\newblock Pitman (Advanced Publishing Program), Boston, MA, 1984.

\bibitem{Slemrod1971Nonexistence}
M.~Slemrod.
\newblock Nonexistence of oscillations in a nonlinear distributed network.
\newblock {\em J. Math. Anal. Appl.}, 36:22--40, 1971.

\bibitem{Sontag1998Mathematical}
E.~D. Sontag.
\newblock {\em Mathematical control theory}, volume~6 of {\em Texts in Applied
  Mathematics}.
\newblock Springer-Verlag, New York, second edition, 1998.
\newblock Deterministic finite-dimensional systems.

\bibitem{Viana2006Ergodic}
M.~Viana.
\newblock Ergodic theory of interval exchange maps.
\newblock {\em Rev. Mat. Complut.}, 19(1):7--100, 2006.

\bibitem{Walters1982Introduction}
P.~Walters.
\newblock {\em An introduction to ergodic theory}, volume~79 of {\em Graduate
  Texts in Mathematics}.
\newblock Springer-Verlag, New York-Berlin, 1982.

\end{thebibliography}
\end{document}